\documentclass[12pt]{article}

\usepackage{
	subfiles,
	enumitem,
	times,
	amsmath,
	amssymb,
	amsthm,
	mathtools,
	fancyhdr,
	color,
	thmtools,
	etoolbox,
	xspace,
	lscape,
	thm-restate,
	tikz-cd,
	subcaption,
	xr,
	tikz,
	rank-2-roots,
	appendix,
	dynkin-diagrams,
	gensymb
}

\usepackage[hypertexnames=false]{hyperref}
\usepackage[margin=2cm]{geometry}

\declaretheorem[style=definition,numberwithin=section]{definition}
\declaretheorem[style=definition,qed=$\oslash$,sibling=definition]{example}

\declaretheorem[style=slanted,sibling=definition]{theorem}
\declaretheorem[style=slanted,sibling=definition]{conjecture}
\declaretheorem[style=slanted,sibling=definition]{lemma}
\declaretheorem[style=slanted,sibling=definition]{proposition}
\declaretheorem[style=slanted,sibling=definition]{corollary}

\declaretheorem[style=definition,sibling=example]{remark}

\AtEndEnvironment{proof}{\setcounter{claim}{0}}

\usepackage[
	backend=bibtex,
	giveninits=true,
	style=trad-alpha,
	url=false]{biblatex}
\newcommand{\C}{\ensuremath{\mathbb{C} }}
\newcommand{\R}{\ensuremath{\mathbb{R} }}
\newcommand{\N}{\ensuremath{\mathbb{N} }}
\newcommand{\Z}{\ensuremath{\mathbb{Z} }}

\newcommand{\triv}{\ensuremath{\mathrm{triv}}}
\DeclareMathOperator{\Tr}{Tr}
\DeclareMathOperator{\End}{End}

\DeclareMathOperator{\ad}{ad}

\DeclareMathOperator{\Motimes}{\text{\raisebox{0.25ex}{\scalebox{1.5}{$\otimes$}}}}

\numberwithin{equation}{section}

\title{Mathieu's approach to the Jacobian Conjecture}
\author{Kevin Zwart\footnote{Email address: \href{mailto:kevin.zwart@ru.nl}{kevin.zwart@ru.nl}}}
\date{%
	IMAPP, Radboud University Nijmegen, The Netherlands \\[2ex]%
	}

\begin{document}
	\maketitle
\abstract{In this paper, we give an expository presentation of a paper of Olivier Mathieu \cite{Mathieu}. The paper of Mathieu proves that a certain Lie group-theoretic conjecture implies the Jacobian Conjecture. To give Mathieu's proof, we first discuss the required background on representation theory in an expository way. We continue to prove some results on the irreducible subrepresentations of the tensor algebra of the standard representation of $SU(N)$. The last part of the paper is dedicated to Mathieu's proof.}

\section{Introduction}\label{chapter:Intro}

The Jacobian Conjecture is both a famous and infamous conjecture. The first formulation of the Jacobian Conjecture was by O.H. Keller in 1939 \cite{Keller}, and was known as Keller's Problem for a long time. Nowadays we call it the Jacobian Conjecture, and is formulated as follows:

\begin{conjecture}[The Jacobian Conjecture]\label{conj:Intro_Jacobian_conj}
	Let $n\in\N$, and let $f:\C^n\rightarrow \C^n$ be a polynomial map, that is to say that every $f_i:\C^n\rightarrow \C$ is a polynomial map, in such a way that the Jacobian matrix $J(f):=(\frac{\partial f_i}{\partial z_j})_{1\leq i,j\leq n}$ satisfies $\det(J(f))=1$. Then $f$ is invertible with polynomial inverse.
\end{conjecture}

Curiously Keller proposed this question for functions $f$ with integer coefficients, and verified the conjecture in case $f$ has an rational inverse. Apart from the $n=1$ case, which is trivial, the Jacobian Conjecture remains an open problem to the author's knowledge. Many erroneous proofs have been presented for the case $n=2$, and some have even been published. For a discussion on this, we refer to \cite{Wright_On_the_Jacobian_Conj}.

However, some progress have been made. Especially for the case $n=2$, most work has been done. For example, in the lecture notes of S. S. Abhyankar \cite{Abhyankar}, he discusses the case $n=2$ extensively. In 1976, Y. Nakai and K. Baba \cite{Nakai-Baba} proved the Jacobian Conjecture to be true in the case $n=2$ when the degrees $m=\deg f_1, n=\deg f_2$ are either coprime, whenever $m=4$ or $n=4$, or when $m=4p\geq n$ where $p$ is prime. In 1983, T.T. Moh proved \cite{MohJacobian} that the Jacobian Conjecture is true assuming $f_1,f_2$ have degree $\leq 100$, and a paper by D. Wright \cite{WrightJacobian} proved that $f=(f_1,f_2)$ is invertible if $J(f)$ is a product of elementary matrices in $GL(2,k)$ for any algebraically closed field $k$. 

In a more general setting, S. Wang proved in 1980 \cite{WangJacobian} that the Jacobian Conjecture is true for any $n$, assuming each $f_i$ has at most degree 2. An influential survey by H. Bass, E. H. Connell and D. Wright \cite{Jacobian_overzicht_artikel} showed that the general Jacobian Conjecture can be proven from the assumption that every $f_i$ is of the form $f_i(z)=z_i-h_i$ where $h_i:\C^n\rightarrow \C$ is a homogeneous polynomial of degree $d$,  for all $i$. In fact, they show that it is enough to prove the Jacobian Conjecture just for the $d=3$ case. Other excellent books on the Jacobian Conjecture have been written in the last 30 years, such as \cite{vdEssen,vdEssen_en_de_rest}. Newer work in 2005 by M. de Bondt and A. Van den Essen \cite{BondtVanDenEssen} shows that it is enough to prove the Jacobian Conjecture if $h:=(h_1,\ldots, h_n)$ satisfies $J(h)$ is a homogeneous nilpotent symmetric matrix. The Jacobian Conjecture has been suggested in the paper by S. Smale in 1998 \cite{Smale} to be one of the 18 great problems of the 21st century. 

In 1997, O. Mathieu \cite{Mathieu} published a paper where he proved the Jacobian Conjecture assuming a Lie-theoretic conjecture to hold. His methods differed to previous attempts, for he was using group and representation theory as tools to get an algebraic result. This Lie-theoretic conjecture is nowadays known as the Mathieu Conjecture. Only one year after Mathieu's paper, Duistermaat and van der Kallen \cite{Duistermaat} proved Mathieu's Conjecture in the case the group under consideration is abelian. Apart from this result, the Mathieu Conjecture is still an open problem, both in examples as in its full generality. Mathieu also stated a second conjecture in his paper, which would imply the Mathieu Conjecture. However, a counterexample to this strange conjecture was produced by M.V. Nori \cite{Nori} using algebraic arguments. In recent years, W. Zhao \cite{Zhao_Mathieu_Subspaces} generalized some ideas of Mathieu to a notion called Mathieu-Zhao spaces, but although interesting, these generalities have not succeeded in proving the Mathieu Conjecture so far.

This present paper is written as an extended presentation of the paper of Mathieu. This paper has two parts: in Section \ref{chapter:Representation_theory} we will introduce some representation theory required for the proof Mathieu gave, and apply it to $SU(n)$ and $SL(n,\C)$. In Section \ref{sec:prerequisites_Mathieu} we will analyze the representation theory of $SL(N,\C)$ on the tensor products of $S^k\C^n$, the $k$-th symmetric power of $\C^n$, with its dual $S^k(\C^n)^*$. In Section \ref{chapter:Mathieu_implies_Jacobian} we will apply the tools of representation theory to prove the Jacobian Conjecture following the lines of Mathieu, where we end with the Mathieu Conjecture.

\section{Review of representation theory}\label{chapter:Representation_theory}
\label{sec:short_overview_repr_theory}

In this paper, we aim to give an expository view of the Mathieu Conjecture, a conjecture that relates representation theory of compact Lie groups to the Jacobian Conjecture. Since the focus of this paper will be on the Mathieu Conjecture, we introduce it now. Unfortunately, there is some jargon to go through, and this subsection is dedicated to giving the definitions to understand the conjecture. The entirety of Section \ref{chapter:Representation_theory} is meant for those that are not that familiar with representation theory of (compact) Lie groups. If the reader is familiar with representation theory, they can skip Section \ref{chapter:Representation_theory} and immediately go to Section \ref{sec:prerequisites_Mathieu} where we start to discuss Mathieu's proof.
\begin{conjecture}[The Mathieu Conjecture] \cite[Main Conjecture]{Mathieu}\label{conj:Mathieu}
	Let $K$ be a compact connected Lie group, and let $f$ and $h$ be complex-valued finite-type functions. Assume $$\int_K f^n(k)dk=0$$ for all $n\in \N$. Then $$\int_K f^n(k)h(k)dk=0$$ for all $n$ large enough. Here $dk$ is the Haar measure on $K$.
\end{conjecture}

The focus of this paper is to prove the following theorem:

\begin{theorem}\label{thm:big_theorem_Mathieu_implies_Jacobian}
	Let $N\in\N$ and let $K=SU(N)$. Assume the Mathieu Conjecture for $SU(N)$. Then the Jacobian Conjecture on $\C^N$ is true.
\end{theorem}

The proof Mathieu gave is given in Section \ref{sec:Mahtieu's_proof}. The rest of this subsection is dedicated to explaining the definitions in the conjecture, while the rest of this section will be dedicated to developing tools to prove that the Mathieu Conjecture implies the Jacobian Conjecture. This section is based on \cite{Broecker_tom_Dieck,Bump,Hall,Humphreys,Kirillov,KnappBeyond,KnappRepresentation,Procesi}. 

\begin{example}\label{ex:Intro_groups}
	Throughout this paper, we will be interested in Lie groups. These groups have, in addition to a group structure, also a manifold structure in such a way that the multiplication and inversion maps are smooth maps. For our purposes, we will only need a few Lie groups, which we will introduce now. Let $N\in\N$ and let us denote $M(N,\C)$ as the $N\times N$ matrices with complex entries. For the purpose of this paper, we only need to consider the following matrix groups: 
	\begin{align*}
		GL(N,\C)&:=\{g\in M(N,\C)\mid\det(g)\neq 0\},\\
		SL(N,\C)&:=\{g\in M(N,\C)\mid\det(g)=1\},\\
		U(N)&:=\{g\in M(N,\C)\mid AA^\dagger=A^\dagger A= \mathbf{1}_N\},\\
		SU(N)&:=U(N)\cap SL(N,\C).
	\end{align*}
	where $A^\dagger$ means taking the adjoint (= the complex conjugate and transpose) of $A$, and $\mathbf{1}_N$ is the $N\times N$ identity matrix. All of these groups are Lie groups, with the smooth structure coming from being closed subsets of $M(N,\C)$, which itself is a manifold by noting $M(N,\C)\simeq \C^{N^2}\simeq \R^{2N^2}$. We remark that $U(N)$ and $SU(N)$ are compact, while $GL(N,\C)$ and $SL(N,\C)$ are not. In addition, $SU(N)$ and $SL(N,\C)$ are simply connected.
	
	When discussing Lie groups, one has to discuss Lie algebras as well. Let $G$ be a Lie group. A Lie algebra $\mathfrak{g}$ is just the tangent space at the identity, i.e. $T_eG$, but can also be described by a a vector space with a bilinear product $[\cdot,\cdot]$ that is antisymmetric and satisfies the Jacobi identity. For the purposes of this paper, we only need to consider the Lie algebras corresponding to the above groups:
	\begin{align*}
		\mathfrak{gl}(N,\C)&:=\{X\in M(N,\C)\}=M(N,\C),\\
		\mathfrak{sl}(N,\C)&:=\{X\in M(N,\C)\mid\Tr(X)=0\},\\
		\mathfrak{u}(N)&:=\{X\in M(N,\C)\mid X+X^\dagger=0\},\\
		\mathfrak{su}(N)&:=\mathfrak{u}(N)\cap \mathfrak{sl}(N),
	\end{align*}
	where $[X,Y]=XY-YX$ for all the cases. Note that each of these Lie algebras corresponds to the Lie groups with the same capital Roman letters, i.e. $\mathfrak{u}(N)$ is the Lie algebra of $U(N)$ etc.
\end{example}

For the rest of this section, we shall reserve the letter $K$ for a compact Lie group and $\mathfrak{k}$ for the Lie algebra corresponding to $K$. If we are talking about \emph{any} Lie group, we will primarily refer to $G$ or $H$, with Lie algebra $\mathfrak{g}$ or $\mathfrak{h}$ respectively. For the purposes of this paper, we give the following theorem:
\begin{theorem}\cite[Thm. 4.2]{Hall}
	Let $K$ be a compact Lie group. Then there exist some $n\in\N$ and a closed subgroup $L\subseteq U(n)$ such that $K\simeq L$. 
\end{theorem}
In other words, without loss of generality we can assume that $K$ is a closed subgroup of $U(n)$. Looking at Theorem \ref{thm:big_theorem_Mathieu_implies_Jacobian}, we are mostly interested in $K=SU(n)$, but occasionally the more general case is useful as well.

\begin{definition}
	Let $G$ be a Lie group. Let $V_\pi$ be a non-zero \emph{complex} Hilbert space with Hermitian inner product $\langle \cdot,\cdot\rangle$, and let $\pi:G\rightarrow GL(V_\pi)$ be a homomorphism such that the map $G\times V_\pi\rightarrow V_\pi, (g,v)\mapsto \pi(g)v$ is continuous. Here $GL(V_\pi)$ is the set of bounded operators on $V_\pi$ with bounded inverse. We call the pair $(\pi,V_\pi)$ a \emph{representation} of $G$. If $V_\pi$ is finite-dimensional, we call the representation \emph{finite-dimensional}.
\end{definition}

\begin{example}
	\begin{enumerate}
	\item Let us begin with the trivial example: let $G$ be any Lie group, and consider the space $\C$. Define the map $\pi_\triv:G\rightarrow GL(1,\C)=\C^\times$ by $\pi_\triv(g)z=z$ for any $z\in\C$. In other words, $\pi_\triv(g)=1$ for all $g$. Then $(\pi_\triv,\C)$ is a representation. We call this representation the \emph{trivial representation}.
	\item Let $G$ be a matrix group, e.g. $GL(N,\C)$ or $SU(N)$ for any $n$ as in Example \ref{ex:Intro_groups}. Define the mapping $\pi_{\mathrm{std}}:G\rightarrow GL(N,\C)$ given by $$\pi_\mathrm{std}(g)v=gv$$ for any $v\in \C^N$. This map is a homomorphism and since $\C^N$ is finite-dimensional, it is automatically continuous. So $(\pi_\mathrm{std},\C^n)$ is a representation of $G$. We call this representation the \emph{standard representation of $G$}.
	\item Let $G$ be a matrix group as in previous example. Let us consider the $k$-th tensor product of $\C^N$, denoted $$T^k\C^N:=\C^N\otimes \C^N\otimes\cdots\otimes \C^N.$$ Let us construct the map $\bigotimes^k\pi_{\mathrm{std}}:G\rightarrow GL(T^k\C^N)$ by $$\Motimes^k \pi_{\mathrm{std}}(g):=\pi_{\mathrm{std}}(g)\otimes\pi_{\mathrm{std}}(g)\otimes\cdots\otimes \pi_{\mathrm{std}}(g).$$ Then $(\bigotimes^k\pi_{\mathrm{std}},T^k\C^N)$ is a representation of $G$.
	\item Let $G$ be a matrix group as in previous examples. Consider the $k$-th symmetric product of $\C^N$, denoted $S^k\C^N$. The symmetric product can be realized as $S^k\C^N=T^k\C^N/I$ where $I$ is the two-sided ideal generated by the elements $v\otimes w - w\otimes v$. One can factor $\bigotimes^k\pi_\mathrm{std}$ through to $S^k\C^N$. Let us denote the resulting map by $\pi_k$ for now. Then $(\pi_k,S^k\C^N)$ is a representation of $G$.
	\item Let $G$ be a matrix group again, and let $(\pi,V)$ be some representation of $G$. Let us define the following map $$\pi^*:G\rightarrow GL(V^*),\qquad\qquad \pi^*(g)f=f\circ \pi(g^{-1})$$ where $f\in V^*$ is a linear functional on $V$. This defines a representation $(\pi^*,V^*)$. This representation is called the \emph{dual of $\pi$}, or the \emph{contragradient of $\pi$}.
	\end{enumerate}
\end{example}

To define finite-type functions, we first need the definition of a matrix coefficient of a finite dimensional representation.

\begin{definition}\label{def:matrix_coeff}
	Let $(\pi,V_\pi)$ be a finite dimensional representation of $K$, and let $u,v\in V_\pi$. Let $\langle\cdot,\cdot\rangle$ be the inner product on $V_\pi$. We define a \emph{matrix coefficient of $\pi$} as any function $m_{u,v}^\pi:K\rightarrow\C$ of the form $$m_{u,v}^\pi(g):=\langle \pi(g)v,u\rangle.$$ Let $C(G)_\pi$ be the space of matrix coefficients of $(\pi,V_\pi)$.
\end{definition}	
\begin{definition}\label{def:finite_type_function}
	Let $G$ be a compact Lie group, and let $f:G\rightarrow \C$. We call $f$ a \emph{finite-type function} if $f$ can be written as a finite sum of matrix coefficients, i.e. there exist a finite family of finite dimensional representations $(\pi_j,V_{\pi_j})$ such that $$f(g)=\sum_{j=1}^kc_j m_{u_j,v_j}^{\pi_j}(g),$$ where $c_j\in\C$ and $u_j,v_j\in V_{\pi_j}$.
\end{definition}

\begin{remark}
	Before we continue towards the Mathieu Conjecture, we would like to take a moment to point out that our definition of a finite-type function is not the most commonly used one. For the sake of completeness, we give the other commonly used definition of a finite-type function.
	\begin{definition}
		Let $G$ be a compact Lie group, and let $f:G\rightarrow \C$ be a continuous function. Define $(L(g)f)(x):=f(g^{-1}x)$ for any $g,x\in G$. We call $f$ a \emph{finite-type function} if the functions $\{L(g)f\mid g\in G\}$ span a finite dimensional vector space.
	\end{definition}
	However, it is known that these two notions are the same:
	\begin{theorem}\cite[Thm. 2.1]{Bump}
		Let $K$ be a compact Lie group, let $f:K\rightarrow \C$ be a continuous function, and define $(L(g)f)(x):=f(g^{-1}x)$. Then the functions $L(g)f$ span a finite dimensional vector space if and only if $f$ is a finite sum of matrix coefficients, i.e. a finite-type function.
	\end{theorem}
	We will use Definition \ref{def:finite_type_function} as our definition for finite-type functions to make the connection to representation theory evident.
\end{remark}

The last ingredient we need for the Mathieu Conjecture is the notion of a (normalized) Haar measure. A Haar measure is a specific kind of measure that is compatible with the group structure. 

\begin{definition}
	Let $G$ be a Lie group, and let $\mu$ be any regular Borel measure. We say that $\mu$ is a \emph{left Haar measure} if, for every measurable function $f$, we get
	$$\int_G f(hg)\,d\mu(g)=\int_G f(g)\,d\mu(g)$$
	for all $h\in G$. Similarly, we call $\mu$ a \emph{right Haar measure} if $$\int_G f(gh)\,d\mu(g)=\int_G f(g)\,d\mu(g)$$ for all $h\in G$. If $\mu$ is both a left and right Haar measure, we call $\mu$ a \emph{Haar measure}. To simplify the notation, we will use $\int_G f(g)dg$ instead of $\int_G f(g)\,d\mu(g)$ and say that $dg$ is a left/ right Haar measure.
\end{definition}
It is a well-known fact that whenever $G$ is compact, a Haar measure $dg$ always exists, is unique up to multiplication with a positive constant, and satisfies $dg(G)<\infty$ \cite[Prop. 1.2]{Bump}. Since $dg(G)<\infty$, we can \emph{normalize} the measure, which is to say, we multiply with a constant such that $dg(G)=1$. 
\begin{definition}
	Let $G$ be a compact Lie group. We call a Haar measure such that $dg(G)=1$ a \emph{normalized} Haar measure.
\end{definition} 
For the rest of this paper, whenever we consider a compact Lie group $G$, we always equip it with a normalized Haar measure, which we will denote by $dg$.

Now that we have given all the terms to be able to read the Mathieu Conjecture, we dedicate the rest of the section by developing more representation theory tools. In Section \ref{sec:facts_Jacobian_Conjecture} we return to the Jacobian Conjecture, and show that the Mathieu Conjecture implies the Jacobian Conjecture using the theory developed.
	
\subsection{Some representation theory}	
Now that the Mathieu Conjecture has been introduced, we will quickly review some representation theory we will need for the proof of Theorem \ref{thm:big_theorem_Mathieu_implies_Jacobian}.

To start, we note that there are ways to build representations out of already known representations. 
\begin{definition}
	Let $(\pi,V_\pi)$ and $(\rho,W_\rho)$ be two representations of a Lie group $G$. One can consider the \emph{direct sum representation}, denoted as $(\pi\oplus\rho, V_\pi\oplus W_\rho)$, where $$(\pi\oplus\rho)(g)(v\oplus w)=(\pi(g)v)\oplus (\rho(g)w).$$ In the same way, we can consider the \emph{tensor product representation}, denoted as $(\pi\otimes \rho,V_\pi\otimes W_\rho)$, where $$(\pi\otimes \rho)(g)(v\otimes w)=(\pi(g)v)\otimes (\rho(g)w).$$
\end{definition}

In addition to representations of $G$, one could be interested in representations of the Lie algebra $\mathfrak{g}$.

\begin{definition}
	Let $\mathfrak{g}$ be any Lie algebra. Let $V_\Pi$ be a non-zero complex Hilbert space, and let $\Pi:\mathfrak{g}\rightarrow \End(V_\Pi)$ be a Lie algebra homomorphism, meaning $$\Pi([X,Y])=[\Pi(X),\Pi(Y)]=\Pi(X)\Pi(Y)-\Pi(Y)\Pi(X)$$ for any $X,Y\in\mathfrak{g}$. Here $\End(V_\Pi)$ is the set of endomorphisms on $V_\Pi$. We call the pair $(\Pi,V_\Pi)$ a \emph{representation} of $\mathfrak{g}$. If $V_\Pi$ is finite dimensional, we call the representation \emph{finite dimensional}.
\end{definition}

If $(\pi,V_\pi)$ is a finite dimensional representation of $G$, it is immediately smooth by standard Lie group arguments. Thus one can take the derivative at the identity element $e$ to get the map $$T_e\pi:T_eG=\mathfrak{g}\rightarrow T_\mathbf{1} GL(V_\pi)=\End(V_\pi),\qquad (T_e\pi)(X)=\left.\frac{d}{dt}\right|_{t=0}\pi(\exp(tX)).$$ Then $(T_e\pi,V_\pi)$ is a representation of $\mathfrak{g}$. 
\begin{example}
	\begin{enumerate}
		\item Let $\Pi_{\mathrm{triv}}(X)=0$ for all $X\in \mathfrak{g}$. Then $(\Pi_{\mathrm{triv}},\C)$ is a representation of $\mathfrak{g}$. We call this the \emph{trivial representation}. If $G$ is a Lie group with Lie algebra $\mathfrak{g}$, then $T_e\pi_{\mathrm{triv}}=\Pi_{\mathrm{triv}}$
		\item Let $\mathfrak{g}\subseteq \mathfrak{gl}(N,\C)$, where $\mathfrak{gl}(N,\C)$ is the Lie algebra of $GL(N,\C)$. Then define the map $\Pi_{\mathrm{std}}:\mathfrak{g}\rightarrow M(N,\C)$ by $$\Pi_{\mathrm{std}}(X)v=Xv$$ for any $v\in \C^n$. Then $(\Pi_{\mathrm{std}},\C^n)$ is a representation of $\mathfrak{g}$. We will also call this representation the \emph{standard representation}. If $G\subseteq GL(N,\C)$ is a Lie group with Lie algebra $\mathfrak{g}$, then $T_e\pi_{\mathrm{std}}=\Pi_{\mathrm{std}}$.
		\item Let $\mathfrak{g}\subseteq \mathfrak{gl}(N,\C)$ as in previous example. Consider $T^k\C^N$,  the $k$-th tensor product of $\C^N$. Consider the map $\Pi:\mathfrak{g}\rightarrow \End(T^k\C^N)$, defined by $$\Pi(X):=\Pi_{\mathrm{std}}(X)\otimes \mathrm{id}\otimes\cdots\otimes \mathrm{id}+\mathrm{id}\otimes \Pi_{\mathrm{std}}(X)\cdots \otimes \mathrm{id}+\cdots+\mathrm{id}\otimes \mathrm{id}\otimes\otimes \mathrm{id}\otimes \Pi_\mathrm{std}(X).$$ Then $(\pi,T^k\C^N)$ is a representation of $\mathfrak{g}$. We will denote $\Pi$ by $\bigotimes \Pi_\mathrm{std}$. If $G\subseteq GL(N,\C)$ is a Lie group with Lie algebra $\mathfrak{g}$, then $T_e(\bigotimes^k\pi_\mathrm{std})=\bigotimes^k\Pi_\mathrm{std}.$
		\item Let $\mathfrak{g}\subseteq \mathfrak{gl}(N,\C)$ as in previous examples. Consider the $k$-th symmetric product of $\C^N$, denoted $S^k\C^N$. One can factor $\bigotimes^k\Pi_\mathrm{std}$ through to $S^k\C^N$. Let us denote the resulting map by $\Pi_k$ for now. Then $(\Pi_k,S^k\C^N)$ is a representation of $\mathfrak{g}$.
		\item Let $\mathfrak{g}$ be a complex Lie algebra, and define the map $$\ad:\mathfrak{g}\rightarrow M(n,\mathbb{K}),\qquad \ad(X)Y:=[X,Y],$$ where $n=\dim(\mathfrak{g})$. Then $(\ad,\mathfrak{g})$ is a representation. This is often called the \emph{adjoint representation}.
	\end{enumerate}
\end{example}		

The above examples suggest that each representation of $\mathfrak{g}$ comes from a representation of $G$. This is, unfortunately, not true. However, whenever $G$ is simply connected, we \emph{do} have a bijection between representations of $G$ and $\mathfrak{g}$. 

\begin{lemma}\label{lemma:representation_G_and_Lie_alg_isomorphic_simply_connected}
	Let $G$ be a simply connected Lie group. Then every representation of $G$ induces a unique representation of $\mathfrak{g}$, and every representation of $\mathfrak{g}$ induces a unique representation of $G$. In other words, there is a bijection between representations of $G$ and $\mathfrak{g}$.
\end{lemma}
\begin{proof}
	Let $(\pi, V_\pi)$ be a representation of $G$. Then by the above, $(T_e\pi,V_\pi)$ is a representation of $\mathfrak{g}$. On the other hand, let $(\rho,W_\rho)$ be a representation of $\mathfrak{g}$. Since $G$ is simply connected, every homomorphism $\phi:\mathfrak{g}\rightarrow\mathfrak{h}$ comes from a unique homomorphism $\Phi:G\rightarrow H$ such that $T_e\Phi=\phi$, see for example \cite[Cor. 1.10.5]{Duistermaat_Kolk}. In particular this means that there exists a representation $(\pi,W_\rho)$ of $G$ such that $T_e\pi=\rho.$
\end{proof}
Next up, we quickly discuss some properties of representations.
\begin{definition}
	Let $(\pi,V_\pi)$ be a representation of $G$ or $\mathfrak{g}$. We call a subspace $W\subseteq V_\pi$ an \emph{invariant subspace} if $\pi(g)W\subseteq W$ for all $g\in G$, and say that $(\pi|_{W},W)$ is a \emph{subrepresentation} of $(\pi,V)$. If the only invariant subspaces are $\{0\}$ and $V_\pi$, we say the representation $(\pi,V_\pi)$ is \emph{irreducible}. Otherwise, we will call it reducible.
\end{definition}

\begin{definition}\label{def:intertwiner}
	Let $(\pi,V)$ and $(\rho,W)$ be two representations of $G$. We call a bounded linear map $A:V\rightarrow W$ a $G$-\emph{intertwiner} or a $G$-\emph{equivariant map} if $$\rho(g)A=A\pi(g)\qquad \forall g\in G.$$ We denote the set of $G$-intertwiners $A:V\rightarrow W$ by $\mathrm{Hom}_G(V,W)$. We say $\pi$ and $\rho$ are \emph{equivalent}, if there exists a $G$-intertwiner that is invertible with bounded inverse. We will denote it with $(\pi,V)\simeq (\rho,W)$, but sometimes we will also use $\pi\simeq \rho$, or $V\simeq W$ when the representations are implicit. 
	
	Equivalences of representations of $\mathfrak{g}$ and $\mathfrak{g}$-intertwiners are defined in the same way.
\end{definition}
	
\begin{remark}
	In representation theory, when considering a representation $(\pi,V)$, one often drops either the $\pi$ or the $V$. We will be complete whenever possible, but sometimes we will use the notation that $V$ is a representation of $G$, where $\pi$ is implicit.
\end{remark}

\begin{definition}\label{def:G_hat}
	Let $K$ be a compact Lie group. We define $\widehat{K}$ as the set of equivalence classes of finite dimensional irreducible representations of $K$. That is to say, if $(\pi,V)$ is a finite dimensional irreducible representation of $K$, we say $[(\pi,V)]\in \widehat{K}$. We will often not make the distinction between $(\pi,V)$ and $[(\pi,V)]$.
\end{definition}

One of the most used tools in representation theory is Schur's Lemma, of which there are many versions. There exist more general versions of Schur's Lemma, such as a description for Hilbert spaces (e.g., \cite[Prop. 1.5]{KnappRepresentation}), but we will mainly be looking at finite dimensional representations. Schur's Lemma becomes the following then:
	\begin{proposition}[Schur's Lemma]\cite[Thm. 4.29]{Hall}
		Let $(\pi,V)$ and $(\rho,W)$ be two irreducible finite dimensional representations of $G$ or $\mathfrak{g}$. Let $T:V\rightarrow W$ be a $G$- or $\mathfrak{g}$-intertwiner. Then either $T$ is an isomorphism, or $T=0$.
	\end{proposition}
	\begin{proof}
		The sets $\ker(T)$ and $\mathrm{im}(T)$ are invariant subspaces. Thus, they must either be 0, or the entire space. The possibilities give the results listed.
	\end{proof}
	\begin{corollary}
		Let $(\pi,V)$ be an irreducible finite dimensional representation of $G$ or $\mathfrak{g}$. Then any $G$- or $\mathfrak{g}$-intertwiner $T$ must be of the form $T=\lambda\,\mathrm{id}$ where $\lambda\in\C$.
	\end{corollary}
	
\begin{proposition}\cite[Prop. 2.2]{Bump}\label{cor:compact_group_completely_reducible}
	Let $K$ be a compact Lie group, and let $(\pi,V_\pi)$ be a finite dimensional representation. Then $$(\pi,V_\pi)\simeq\left(\bigoplus_{j=1}^k\pi_j,\bigoplus_{j=1}^k V_j\right),$$ where $(\pi_j,V_j)\in \widehat{K}$ and $\pi_j(g):=\pi(g)|_{V_j}$. If such a decomposition happens, we say the representation $(\pi,V)$ is \emph{completely reducible.}
\end{proposition}

\begin{definition}\label{def:multiplicity}
	Let $(\pi,V)$ be a completely reducible finite dimensional representation of $G$, i.e. $$(\pi,V)\simeq \left(\bigoplus_{i=1}^n\pi_i,\bigoplus_{i=1}^n V_i\right).$$ Let $(\rho,W)\in \widehat{G}$. We say \emph{$(\rho,W)$ occurs in $V$} if $(\rho,W)\simeq(\pi_i,V_i)$ for some $i$ in the above decomposition. We say \emph{$(\pi_i,V_i)$ occurs in $V$ with multiplicity $m$} if $$\dim(\mathrm{Hom}_G(V,V_i))=m.$$ We say $(\pi,V)$ is \emph{multiplicity free} if each irreducible subrepresentation has multiplicity at most 1.
\end{definition}

When we are considering compact Lie groups, one has to mention one of the big achievements in representation theory of compact groups: the Peter-Weyl theorem. We will not state the theorem itself, but only one of its many corollaries:

\begin{lemma}\label{lemma:compact_grp_irreduc_means_finite_dim}
	Let $K$ be a compact Lie group and let $(\pi,V_\pi)$ be an irreducible representation of $K$. Then $V_\pi$ is finite dimensional.
\end{lemma}

Combining Lemma \ref{lemma:compact_grp_irreduc_means_finite_dim} with Proposition \ref{cor:compact_group_completely_reducible} gives the following corollary:

\begin{corollary}
	Let $K$ be a compact Lie group. Then any irreducible representation is completely reducible.
\end{corollary}

We will be interested in decomposing representations into their irreducible parts. Since all irreducible representations are finite dimensional, we will only consider finite dimensional representations. And when considering finite dimensional representations, we can consider matrix coefficients introduced in Definition \ref{def:matrix_coeff}. We sum up some properties of matrix coefficients.

\begin{lemma}[Schur's orthogonality relations]\cite[Cor 1.10]{KnappRepresentation}\label{lemma:Schur_orthogonality}
	Let $K$ be a compact Lie group, and let $(\pi,V)$ and $(\rho,W)$ be two irreducible finite dimensional representations. If $\pi$ and $\rho$ are not equivalent, then $$\int_K m_{v,w}^\pi(k)\overline{m_{v',w'}^\rho(k)}\,dk=0$$ for any $v,w\in V$ and $v',w'\in W$. Here $dk$ is the Haar measure on $K$. In other words, the spaces of matrix coefficients are mutually orthogonal with respect to the $L^2$-norm whenever the representations are not equivalent.
\end{lemma}
Another tool that is useful and will be used later is the notion of a character:
\begin{definition}
	Let $(\pi,V)$ be a finite dimensional representation of $K$. Then we define \emph{the character of $\pi$} as the function $$\chi_\pi:K\rightarrow\C,\qquad \chi_\pi(k):=\Tr(\pi(k)).$$
\end{definition}
\noindent Note that $\chi_\pi$ is conjugation invariant, which is to say that $\chi_\pi(xkx^{-1})=\chi_\pi(k)$ for all $x\in K$. In addition, let $(\pi,V)$ be a finite dimensional representation. If we choose an orthonormal basis $e_1,\ldots,e_n$ on $V$, we get that $$\chi_\pi(k)=\sum_{i=1}^nm_{e_i,e_i}^\pi(k).$$
With this realization, we can use Schur's orthogonality relations and some elementary permutations to derive the following lemma.
\begin{lemma}\label{lemma:Repr_theory_review_character_results}
	\begin{enumerate}[label=(\roman*)]
		\item Let $(\pi_1,V_1)$ and $(\pi_2,V_2)$ be representations of $K$. Then $$\chi_{\pi_1\oplus\pi_2}=\chi_{\pi_1}+\chi_{\pi_2},\qquad\qquad\chi_{\pi_1\otimes\pi_2}=\chi_{\pi_1}\chi_{\pi_2}.$$
		\item Let $(\pi,V)$ be a representation of $K$, and let $(\pi^*,V^*)$ be the dual of $(\pi,V)$. Then $$\chi_{\pi^*}(k)=\overline{\chi_\pi(k)}.$$
		\item Let $\langle\cdot,\cdot\rangle$ be the $L^2(K)$-inner product. Let $(\pi,V)$ and $(\rho,W)$ be finite dimensional irreducible representations of $K$. Then $\chi_\pi=\chi_\rho$ if and only if $\pi\simeq \rho$ and $$\langle \chi_\pi,\chi_\rho\rangle=\begin{cases}
			1 \text{ if } \pi \simeq \rho\\
			0 \text{ if } \pi \not\simeq \rho
		\end{cases}.$$
		\item Let $(\pi,V)$ be any representation of $K$. Because $K$ is compact, $\pi$ is completely reducible, i.e. $$V=\bigoplus_{i}V_i$$ as representations, where $(\pi_i,V_i)$ is an irreducible representation, and $\pi_i(g)=\pi(g)|_{V_i}$. Then, by using the above, $$\langle\chi_\pi,\chi_{\pi_i} \rangle = m_i$$ where $m_i$ is the multiplicity of $(\pi_i,V_i)$ in the decomposition of $(\pi,V)$ into irreducibles.
	\end{enumerate}
\end{lemma}
\noindent We will not use characters until Section \ref{sec:tensor_prod} to study the multiplicities of tensor products of irreducible representations. However, they are useful in many situations to show irreducibility, or to show that two representations are equivalent.

\subsection{Root system of a compact group}\label{sec:Root_system_compact_group}

After having introduced most of the basic representation theory, we move on to discuss some facts about the compact Lie group $K$ itself. One of the significant results of Lie theory is the characterization of complex semisimple Lie algebras using root systems, see for example \cite{Humphreys, Kirillov}. Not only does this give rise to a complete characterization of the Lie algebra, it also defines parabolic subgroups of $G$ (see for example \cite{KnappBeyond,KnappRepresentation}), makes certain decompositions like the Iwasawa decomposition or the Bruhat decomposition possible \cite{KnappBeyond}, and characterizes all finite dimensional irreducible representations (see the Theorem of Highest Weight, see Section \ref{sec:Highest_weight_repres}). For completeness sake, let us give the root decomposition.

Before defining the root system, we briefly want to discuss the fact that any connected compact Lie group is reductive. We recall the definition of a reductive Lie algebra and Lie group.

\begin{definition}
	Let $\mathfrak{g}$ be any Lie algebra. We call $\mathfrak{g}$ \emph{reductive} if the Lie algebra can be written as the direct sum $$\mathfrak{g}=\mathfrak{z}(\mathfrak{g})\oplus[\mathfrak{g},\mathfrak{g}].$$ Here $\mathfrak{z}(\mathfrak{g})$ is the center of $\mathfrak{g}$. We call $\mathfrak{g}$ \emph{semisimple} if $\mathfrak{g}$ is reductive and $\mathfrak{z}(\mathfrak{g})=0$.
	In addition, if $G$ is a Lie group with Lie algebra $\mathfrak{g}$, we call $G$ \emph{reductive} or \emph{semisimple} if $\mathfrak{g}$ is reductive or semisimple respectively.
\end{definition}
\begin{example}
	Consider the case of $SU(N)$ with Lie algebra $\mathfrak{su}(N)$, as described in Example \ref{ex:Intro_groups}. It can be shown that $\mathfrak{su}(N)$ is semisimple.
\end{example}

\begin{theorem}[Root system decomposition for complex semisimple Lie algebras]\cite[Thm. 6.44]{Kirillov}\label{thm:Root_decomp}
	Let $\mathfrak{g}$ be a complex semisimple Lie algebra, and let $\mathfrak{h}\subseteq \mathfrak{g}$ be a maximal abelian subalgebra ($\mathfrak{h}$ is called a Cartan subalgebra). Then
	\begin{enumerate}
		\item The following decomposition holds $$\mathfrak{g}=\mathfrak{h}\oplus\bigoplus_{\alpha\in \Delta(\mathfrak{g},\mathfrak{h})}\mathfrak{g}_\alpha$$ where 
		\begin{align*}
			\mathfrak{g}_\alpha&:=\{X\in \mathfrak{g}\mid \ad(H)X=\alpha(H)X\text{ for all } H\in\mathfrak{h}\},\\
			\Delta(\mathfrak{g},\mathfrak{h})&:=\{\alpha\in\mathfrak{h}^*\setminus \{0\}\mid \mathfrak{g}_\alpha\neq 0\}.
		\end{align*}
		We call $\Delta(\mathfrak{g},\mathfrak{h})$ the \emph{root system} of $\mathfrak{k}$ and the subspaces $\mathfrak{g}_\alpha$ the \emph{root spaces}.
		\item If $\alpha, \beta$ are roots such that $\alpha+\beta$ is a root, then $[\mathfrak{g}_\alpha,\mathfrak{g}_\beta]=\mathfrak{g}_{\alpha+\beta}$.
		\item Let $(\cdot,\cdot)$ be a non-degenerate symmetric bilinear form on $\mathfrak{h}$ that is invariant under $\ad$. Then $$\frac{2(\alpha,\beta)}{(\alpha,\alpha)}\in\Z$$ for all roots $\alpha,\beta$.
		\item For any $\alpha\in\Delta(\mathfrak{g},\mathfrak{h})$ define the map $s_\alpha:\mathfrak{h}^*\rightarrow\mathfrak{h}^*$ by $$s_\alpha(\lambda)=\lambda-\frac{2(\alpha,\lambda)}{(\alpha,\alpha)}\alpha.$$ If $\alpha$ and $\beta$ are roots, then $s_\alpha(\beta)$ is also a root.
	\end{enumerate} 
\end{theorem}

We would also like to do some root decomposition for compact Lie groups. However, the following lemma shows a problem.
\begin{lemma}\label{lemma:compact_grp_complex_is_finite}
	Let $K$ be a compact complex Lie group. Then $K$ is finite.
\end{lemma}
\begin{proof}
	Let $\mathfrak{k}$ be the Lie algebra of $K$. It is equivalent to show that $\mathfrak{k}=0$ since $K$ is compact. Since $K\subseteq U(n)$, so $\mathfrak{k}\subseteq \mathfrak{u}(n)$. In particular, if $X\in\mathfrak{k}$ then $X+X^\dagger=0$. Let $Y\in\mathfrak{k}$. Since $\mathfrak{k}$ is $\C$ linear, both $Y\in\mathfrak{k}$ and $iY\in\mathfrak{k}$. But this gives $Y+Y^\dagger = 0$ and $i(Y-Y^\dagger)=0$. This can only be true if $Y=0$.
\end{proof}

In other words, finding roots on the Lie algebra $\mathfrak{k}$ requires more work. To solve this problem, we recall the complexification of a Lie algebra.

\begin{definition}
	Let $\mathfrak{k}$ be a real Lie algebra. We define the \emph{complexification} of $\mathfrak{k}$ as the complex Lie algebra $\mathfrak{k}_\C:=\mathfrak{k}\otimes \C =\mathfrak{k}\oplus i\mathfrak{k}$ with Lie bracket $[x+iy,x'+iy'] = [x,x']-[y,y']+i([x,y']+[y,x']).$ If $\mathfrak{g}$ has a subalgebra $\mathfrak{k}$ such that $\mathfrak{k}_\C\simeq\mathfrak{g}$, we say that $\mathfrak{k}$ is a \emph{real form} of $\mathfrak{g}$.
\end{definition}

With the previous example in mind, we will focus on a semisimple Lie algebra $\mathfrak{k}$, coming from a compact Lie group $K$. Because of Cartan's semisimplicity criterion, it means that $\mathfrak{k}_\C$ is also semisimple. This allows us to decompose $\mathfrak{k}_\C$ into root spaces. Let us discuss this process in more detail.
\begin{definition}
	Let $K$ be a compact connected Lie group. A subgroup $T$ is called a \emph{torus} if $T$ is connected and abelian.
\end{definition}
\begin{proposition}\cite[Prop. 15.3]{Bump}\label{prop:what_is_a_torus_really?}
	Let $K$ be a connected compact Lie group and let $T$ be a torus. Let $\mathfrak{t}$ be the Lie algebra of $T$. Then $T\simeq (\R/\Z)^r\simeq (S^1)^r$ for some $r$, where $r$ is the dimension of $\mathfrak{t}$.
\end{proposition}
Let us consider a partial order of tori, where the partial order is inclusion. If $T,T'$ are tori with $T\subseteq T'$, it is clear by Proposition \ref{prop:what_is_a_torus_really?} that $\dim(T)\leq \dim(T')$. Since there is a maximum to the dimension, namely $\dim(K)$, any torus is contained in a maximal torus. This maximal torus plays a role that is similar to a `Cartan subalgebra' in the following way:

\begin{proposition}\cite[Prop. 11.7 and Thm. 11.9]{Hall}\label{prop:torus_is_maximal_abelian}
	Let $K$ be a connected compact Lie group. The maximal tori in $K$ are exactly the analytic subgroups $T$ corresponding to the maximal abelian subalgebras $\mathfrak{t}$ of $\mathfrak{k}$, where $\mathfrak{t}$ is the Lie algebra of $T$. In addition, all maximal tori are conjugate. That is to say, if $T'$ is another maximal torus, then there exists $k\in K$ such that $T'=kTk^{-1}$.
\end{proposition}

\begin{proof}
	Let $T$ be a maximal torus, and $\mathfrak{t}$ be its Lie algebra. If $\mathfrak{t}$ is not maximal, then there exists an abelian subalgebra $\mathfrak{t}'$ such that $\mathfrak{t}\subseteq \mathfrak{t}'$. The analytic subgroup $T'$ corresponding to $\mathfrak{t}'$ then is also abelian and satisfies $T\subseteq T'$. So, $T$ is not maximal. Proving the converse goes similarly.
\end{proof}

Now let $T$ be a maximal torus in $K$, and let $\mathfrak{t}$ be its Lie algebra. By Proposition \ref{prop:torus_is_maximal_abelian}, $\mathfrak{t}$ is maximal abelian in $\mathfrak{k}$. This also means that $\mathfrak{t}_\C$ is maximal abelian in $\mathfrak{k}_\C$. This means that we can apply Theorem \ref{thm:Root_decomp} to get the following theorem

\begin{theorem}\label{thm:Root_decomp_compact_grp}
	Let $K$ be a semisimple Lie group, let $\mathfrak{k}$ be its Lie algebra, and let $T$ be a maximal torus in $K$. Let us define $\mathfrak{g}:=\mathfrak{k}_\C$. Then we have the following decomposition for $\mathfrak{g}$ with respect to $\ad(\mathfrak{t}_\C)$:
	\begin{align}\label{eq:root_decomp_compact_grp}
		\mathfrak{k}_\C=\mathfrak{g} = \mathfrak{t}_\C\oplus\bigoplus_{\alpha\in \Delta(\mathfrak{k}_\C,\mathfrak{t}_\C)}\mathfrak{g}_\alpha
	\end{align}
	where
	\begin{align*}
		\mathfrak{g}_\alpha&:=\{X\in\mathfrak{g}\mid\ad(H)X=\alpha(H)X\text{ for all }H\in\mathfrak{t}_\C\}, \qquad \Delta(\mathfrak{k}_\C,\mathfrak{t}_\C):=\{\alpha\in(\mathfrak{t}_\C)^*\setminus \{0\}\mid \mathfrak{g}_\alpha\neq 0\}.
	\end{align*}
	We call $\Delta(\mathfrak{k}_\C,\mathfrak{t}_\C)$ the \emph{root system} of $\mathfrak{k}$ and the subspaces $\mathfrak{g}_\alpha$ the \emph{root spaces}.
\end{theorem}
The usual way textbooks handle the next part, is to go in depth on root systems and give all kinds of interesting features. We, however, will not need these, for we are only interested in the case of $SU(n)$. We only require two definitions.
\begin{definition}
	Let $\Delta$ be a root system, such as $\Delta(\mathfrak{g},\mathfrak{t}_\C)$. We remark that $\Delta\subseteq \R^n$ for some $n$. Let $(\cdot,\cdot)$ be an inner product on $\R^n$. Let $v\in\R^n$ such that for any $\alpha\in\Delta$ we have that $(\alpha,v)\neq 0$. Then we can write $$\Delta=\Delta^+\sqcup \Delta^-$$ where $\Delta^{+}=\{\alpha\in\Delta\mid (\alpha,v)>0\}$ and $\Delta^{-}=\{\alpha\in\Delta\mid (\alpha,v)<0\}$. We call a root $\alpha\in \Delta^+$ a \emph{positive root}.
	In addition, we call a root $\alpha\in\Delta^+$ \emph{simple} if it cannot be written as a sum of two positive roots. Denote $\Phi$ by the set of simple roots.
\end{definition}
\begin{lemma}\cite[Lemma 7.13]{Kirillov}
	Let $\beta$ be a positive root. Then $$\beta=\sum_{\alpha\in\Phi}k_\alpha\alpha$$ where $k_\alpha\in\N_0$.
\end{lemma}

\begin{example}\label{ex:root_system_SU(N)}
	Consider the group $SU(N)$ with $N\geq 2$. This example is especially useful, since this specific example will be necessary to prove the Jacobian Conjecture using the Mathieu Conjecture. Therefore we will discuss this example in detail and will frequently come back to it throughout this paper. The group $SU(N)$ is semisimple as mentioned in Example \ref{ex:Intro_groups}. Note that $\mathfrak{su}(N)_\C=\mathfrak{sl}(N,\C)$. Let $$T:=\left\{\left.\begin{pmatrix}
		e^{i\phi_1}&&&\\
		&e^{i\phi_2}&&\\
		&&\ddots&\\
		&&&e^{i\phi_N}
	\end{pmatrix}\right|\phi_j\in\R\text{ such that }\sum_{j=1}^N\phi_j=0\right\}.$$ Then $$\mathfrak{t}=\left\{\left.\begin{pmatrix}
		i\phi_1&&&\\
		&i\phi_2&&\\
		&&\ddots&\\
		&&&i\phi_N
	\end{pmatrix}\right|\phi_j\in\R\text{ such that }\sum_{j=1}^N\phi_j=0\right\},$$ and its complexification is $$\mathfrak{h}=\mathfrak{t}_\C=\left\{\left.\begin{pmatrix}
		\phi_1&&&\\
		&\phi_2&&\\
		&&\ddots&\\
		&&&\phi_N
	\end{pmatrix}\right|\phi_j\in\C\text{ such that }\sum_{j=1}^N\phi_j=0\right\}\subseteq \mathfrak{sl}(N,\C).$$ To find the roots, we define the functionals $\epsilon_1,\ldots,\epsilon_N:\mathfrak{h}\rightarrow\C$ by $$\epsilon_{i}\left(\mathrm{diag}(\phi_1,\ldots,\phi_N)\right)=\phi_{i}.$$ A calculation shows that the roots are given by $$\Delta(\mathfrak{g},\mathfrak{h})=\{\epsilon_{i}-\epsilon_j\mid i\neq j\},\qquad \Delta^+(\mathfrak{g},\mathfrak{h})=\{\epsilon_i-\epsilon_j\mid i>j\}$$ and the simple roots are given by $$\Phi=\{\alpha_i\mid i=1,\ldots,N-1\},\qquad \alpha_i:=\epsilon_{i+1}-\epsilon_{i}.$$
\end{example}
We will quickly touch upon the notion of the Weyl group. 
\begin{definition}\label{def:analytic_Weyl_group_compact}
	Let $\Delta$ be a root system. We recall $\Delta\subset \R^n$ for some $n$. The \emph{Weyl group} of $\Delta$ is the subgroup generated by the reflections $s_\alpha$, $\alpha\in\Delta$.
\end{definition}
Although there is plenty to say about the Weyl group, we will skip the discussion and refer to excellent sources as \cite{Broecker_tom_Dieck,Hall,KnappBeyond}.

\subsection{Characterizing representations}
\subsubsection{Dominating dominant integral weights}
Leaving compact groups behind for a moment, we would like to discuss the theory of dominant integral weights, and highest weight representations. In this section, we will discuss them for complex semisimple Lie algebras, while in Section \ref{sec:highest_weight_repr_compact_grp} we discuss the relation with compact Lie groups. This section is based on \cite{Bump,Hall,Humphreys,KnappRepresentation,Kirillov}.

For the rest of this section, let $\mathfrak{k}$ be a semisimple Lie algebra over $\R$, as we had before, and let $\mathfrak{g}:=\mathfrak{k}_\C$ be a complex semisimple Lie algebra. Choose any maximal abelian subalgebra $\mathfrak{t}$ of $\mathfrak{k}$, and let $\mathfrak{h}:=\mathfrak{t}_\C$. Taking Theorem \ref{thm:Root_decomp_compact_grp} into account, let $\Delta:=\Delta(\mathfrak{g},\mathfrak{t}_\C)$ be the set of roots, choose a set of positive roots $\Delta^+$ and a set of simple roots $\Phi$. Choose any bilinear non-degenerate symmetric form on $\mathfrak{h}$, such as the Killing form, and this bilinear form induces a bilinear non-degenerate symmetric form on $\mathfrak{h}^*$. Let us denote the latter bilinear form by $(\cdot,\cdot)$.
\begin{definition}
	For every $\alpha\in\Delta$ we define $H_\alpha\in \mathfrak{h}$ as the unique element such that $$\beta(H_\alpha)=2\frac{(\alpha,\beta)}{(\alpha,\alpha)}$$ for every $\beta\in\mathfrak{h}^*$ where $(\cdot,\cdot)$ is the inner product on $\mathfrak{h}^*$. We call this the \emph{coroot} of $\alpha$.
\end{definition}
The set of coroots is interesting in its own right. The coroots also form an abstract root system, with the simple roots given by $H_{\alpha_1},\ldots,H_{\alpha_n}$, where $\Phi=\{\alpha_1,\ldots,\alpha_n\}$. The coroots also satisfy $$\beta(H_\alpha)\in\Z$$ since $\frac{2(\alpha,\beta)}{(\alpha,\alpha)}\in\Z$ by Theorem \ref{thm:Root_decomp}. This leads to the following definition:
\begin{definition}
	Let $\mu\in\mathfrak{h}^*$. We say $\mu$ is \emph{integral} if $$\mu(H_\alpha)=2\frac{(\mu,\alpha)}{(\alpha,\alpha)}\in\Z$$ for all $\alpha\in\Delta$. We say $\mu$ is \emph{dominant} if $$(\mu,\alpha)\geq 0$$ for all $\alpha\in\Phi$. Let $P^+$ be the set of $\mu\in\mathfrak{h}^*$ that are dominant integral. 
\end{definition}
Note that all roots are integral, but not all integral elements are roots. We consider an example.

\begin{example}\label{ex:roots_fundamental_su(3)}
	We consider $\mathfrak{su}(3)$. Then $\mathfrak{su}(3)_\C=\mathfrak{sl}(3,\C)$, and $\Phi=\{\alpha_1,\alpha_2\}$ as in Example \ref{ex:root_system_SU(N)}. The simple roots $\Phi$ span $\R^2$, and if we choose $(\alpha_1,\alpha_1)=2=(\alpha_2,\alpha_2)$, we get that $(\alpha_1,\alpha_2)=-1$. This gives Figure \ref{fig:root_system_A_2}. We colored the area where all the dominant $\mu\in\mathfrak{h}^*$ live in gray. Note that the root system can also be spanned by the dominant elements $\omega_1,\omega_2\in\mathfrak{h}^*$ such that $$\omega_i(H_{\alpha_j})=\delta_{ij}.$$ Solving the equations gives $$\omega_1=\frac{1}{3}(2\alpha_1+\alpha_2),\qquad \omega_2=\frac{1}{3}(\alpha_1+2\alpha_2).$$ The linear functionals $\omega_1,\omega_2$ are also drawn in Figure \ref{fig:root_system_A_2}.
	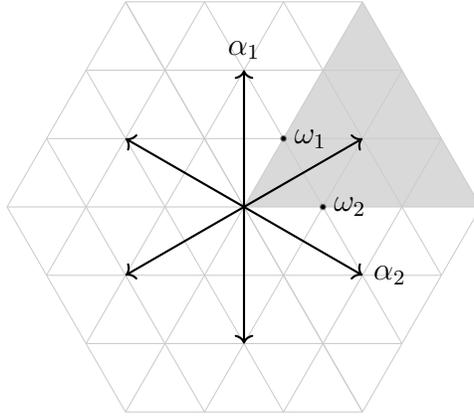
\begin{figure}
		\centering
		\begin{tikzpicture}[scale=0.7]
			\foreach \angle in {30,90,...,330} {
				\draw[->, thick, black] (0,0) -- ({2.6*cos(\angle)}, {2.6*sin(\angle)});
			}
			\begin{rootSystem}[weight length=1.5cm,weight radius=1.5pt]{A}
				\node [above] at \Root {1}{0} {\(\alpha_1\)};
				\node [right] at \Root {0}{1} {\(\alpha_2\)};
				\node [right] at \Root {1}{-1} {\(\)};
				\fundamentalweights
				\node [right] at \weight {1}{0} {\(\omega_2\)};
				\node [right] at \weight {0}{1} {\(\omega_1\)};
				\WeylChamber
			\end{rootSystem}
		\end{tikzpicture}
		\caption{The root system for $A_2$. The grayed-out area are the dominant elements $\mu\in\mathfrak{h}^*$, $\alpha_1,\alpha_2$ are the simple roots, and $\omega_1$ and $\omega_2$ are the fundamental weights.}
		\label{fig:root_system_A_2}
	\end{figure}
\end{example}
We define the $\omega_j$ now, as we will need them later.
\begin{definition}\label{def:fundamental_weight}
	Let $\Delta$ be a root system of $\mathfrak{g}$. Let $\Phi=\{\alpha_1,\ldots,\alpha_n\}$ be the set of simple roots. We define $\omega_i\in\mathfrak{h}^*$ as the linear functional satisfying $$\omega_i(H_{\alpha_j})=\delta_{ij}.$$ We call $\omega_i$ a \emph{fundamental weight}. 
\end{definition}
We note that the fundamental weights are integral. Also note that the fundamental weight $\omega_i$ is dominant, for $$(\omega_i,\alpha)=2(\alpha,\alpha)\omega_i(H_\alpha)\geq 0,$$ so $\omega_i\in P^+$ In addition, $\omega_1,\ldots,\omega_n$ is a basis of $\mathfrak{h}^*$, and it is not hard to see that $\bigoplus_j\Z\alpha_j\subset \bigoplus_j \Z\omega_j$ while the opposite is not true, as can be seen in Example \ref{ex:roots_fundamental_su(3)}. We also note that $P^+=\bigoplus_j\N_0\omega_j.$

\begin{example}\label{ex:root_system_SU(N)_part_2}
	Let us go back to the case of $SU(N)$, and let us continue where we left off in Example \ref{ex:root_system_SU(N)}. The simple roots of $\mathfrak{sl}(N,\C)$ are $\Phi=\{\alpha_1,\ldots,\alpha_{N-1}\}$, and the root system is of type $A_{N-1}$. Let us calculate $H_{\alpha_j}$ and $\omega_j$. We choose an inner product on $\mathfrak{h}^*$ such that $\|\alpha_i\|=2$ for all $i$. The Cartan matrix of $A_{N-1}$ then dictates the value of $(\alpha_i,\alpha_j)$, which is
	\begin{align}
		M=\begin{pmatrix}
			2&-1&&&\\
			-1&2&\ddots&&\\
			&\ddots&\ddots&-1\\
			&&-1&2
		\end{pmatrix},
	\end{align}
	where $M_{ij}=\frac{2(\alpha_j,\alpha_i)}{(\alpha_i,\alpha_i)}$ \cite[Ex. 7.45]{Kirillov}. In this case, $H_\alpha$ are given by $$H_{\epsilon_i-\epsilon_j}=\mathrm{diag}(0,\ldots,0,1,0,\ldots,0 -1,0\ldots 0),$$ where the $1$ is on the $i$-th position, and the $-1$ on the $j$-th position. From the definition of fundamental weights, one can solve the linear equations to find the explicit form of $\omega_j$. For our purpose, we only calculate $\omega_1$ and $\omega_{N-1}$, and they are given by 
	\begin{align}
		\omega_1&=\frac{1}{N}((N-1)\alpha_1+(N-2)\alpha_2+\cdots +\alpha_{N-1}),\label{eq:omega_1_explicit}\\
		\omega_{N-1}&=\frac{1}{N}(\alpha_1+2\alpha_2+\cdots+(N-1)\alpha_{N-1})\label{eq:omega_{n-1}_explicit}.
	\end{align}
	In general, the fundamental weights of the root system of type $A_{N-1}$ are $$\omega_j=\frac{1}{N}\left[(N-j)(\alpha_1+\cdots+(j-1)\alpha_{j-1})+j((N-j)\alpha_j+\cdots+\alpha_{N-1})\right]$$ for $j=1,\ldots,N-1$, see for example \cite[Planche 1]{Bourbaki_ch4_6}.
\end{example}

To end the section, we define a partial order on $\mathfrak{h}^*.$ These will be used to formulate the theorem of highest weight representation, so we will introduce them now.

\begin{definition}
	Let $\Delta$ be a root system of $\mathfrak{g}$, and let $\Phi=\{\alpha_1,\ldots,\alpha_n\}$. Let $\lambda,\mu\in\mathfrak{h}^*$. We say that $\lambda$ is \emph{higher} than $\mu$, denoted as $\mu\preceq\lambda$, if $$\lambda-\mu = \sum_{j=1}^nn_j\alpha_j$$ where $n_j\in\N_0$. 
\end{definition}
Note that the $\preceq$ defines a partial order on $\mathfrak{h}^*$, especially on $P^+$ by restriction. Although the partial ordering is subtle, we will keep the discussion brief. For more information, we recommend \cite[Ch. 8 and 9]{Hall}. We require one more lemma.

\begin{lemma}\cite[Prop. 8.42]{Hall}
	Let $\mu\in\mathfrak{h}^*$ be dominant. Then $w(\mu)\preceq \mu$ for all $w\in W$, where $W$ is the Weyl group.
\end{lemma}

\begin{example}
	Let us go back to the $A_2$ example in Example \ref{ex:roots_fundamental_su(3)}. Let us consider $\lambda$, for example $\lambda=2\omega_1$. In Figure \ref{fig:root_system_A_2_preceq} we denoted all $\mu\in\mathfrak{h}^*$ satisfying $\mu\preceq\lambda$ with a black dot. The grayed-out area represents the convex hull of all $\mu\preceq\lambda$. Of course, not all linear functionals colored like this are applicable; for example, $\mu=\alpha_1$ or $\mu=\alpha_2$ do not satisfy $\mu\preceq\lambda$.
\end{example}	
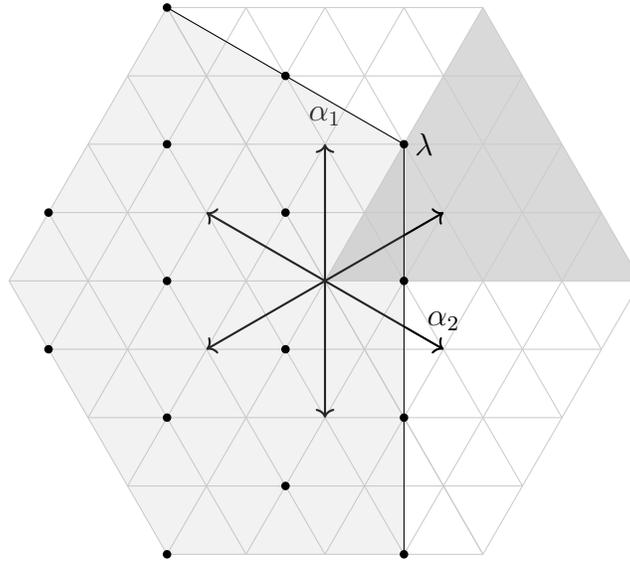
\begin{figure}
	\centering
	\begin{tikzpicture}[scale=0.70]
		\draw[->, thick, black] (0,0) -- ({2.6*cos(90)}, {2.6*sin(90)}) node[above=1mm]{$\alpha_1$};
		\draw[->, thick, black] (0,0) -- ({2.6*cos(330)}, {2.6*sin(330)}) node[above=1mm]{$\alpha_2$};
		\foreach \angle in {30,150,210,270} {
			\draw[->, thick, black] (0,0) -- ({2.6*cos(\angle)}, {2.6*sin(\angle)});
		}
		\begin{rootSystem}[weight length=1.5cm,weight radius=1.5pt]{A}
			\fill[gray!50,opacity=.2] (hex cs:x=3,y=-4) -- (hex cs:x=0,y=2) --
			(hex cs:x=-4,y=4) --
			(hex cs:x=-4,y=0) --
			(hex cs:x=0,y=-4);
			\draw[black, thin] (hex cs:x=3,y=-4) -- (hex cs:x=0,y=2) --
			(hex cs:x=-4,y=4);
			\node [right] at \weight {0}{2} {\(\lambda\)};
			\filldraw (hex cs:x=0,y=2) circle (2pt);
			\filldraw (hex cs:x=1,y=0) circle (2pt);
			\filldraw (hex cs:x=2,y=-2) circle (2pt);
			\filldraw (hex cs:x=3,y=-4) circle (2pt);
			\filldraw (hex cs:x=-1,y=1) circle (2pt);
			\filldraw (hex cs:x=-2,y=3) circle (2pt);
			\filldraw (hex cs:x=0,y=-1) circle (2pt);
			\filldraw (hex cs:x=1,y=-3) circle (2pt);
			\filldraw (hex cs:x=-2,y=0) circle (2pt);
			\filldraw (hex cs:x=-3,y=2) circle (2pt);
			\filldraw (hex cs:x=-4,y=4) circle (2pt);
			\filldraw (hex cs:x=0,y=-4) circle (2pt);
			\filldraw (hex cs:x=-1,y=-2) circle (2pt);
			\filldraw (hex cs:x=-4,y=1) circle (2pt);
			\filldraw (hex cs:x=-3,y=-1) circle (2pt);
			\WeylChamber
		\end{rootSystem}
	\end{tikzpicture}
	\caption{An example of all possible $\mu\in\mathfrak{h}^*$ that could possibly satisfy $\mu\preceq\lambda$, denoted by the black dots. The grayed out area is the convex hull of all possibilities.}
	\label{fig:root_system_A_2_preceq}
\end{figure}

\subsection{Theorem of highest weight representation}\label{sec:Highest_weight_repres}
\subsubsection{Irreducible representation of a complex semisimple Lie algebra}
With that, we return to representation theory. In particular, let us begin with a definition that has been mentioned but not discussed in depth: weights. As before, let $\mathfrak{g}=\mathfrak{k}_\C$ be a complex semisimple Lie algebra, and let $\mathfrak{h}=\mathfrak{t}_\C$ be a maximal abelian subalgebra of $\mathfrak{g}$. Let $\Delta(\mathfrak{g},\mathfrak{h})=\Delta$ be the set of roots. Choose a set of positive roots $\Delta^+$ and let $\Phi$ be the set of simple roots corresponding to this choice.
\begin{definition}
	Let $(\pi,V)$ be a representation of $\mathfrak{g}$. An element $\lambda\in\mathfrak{h}^*$ is called a \emph{weight} if there exists a non-zero vector $v\in V$ such that $$\pi(H)v=\lambda(H)v$$ for all $H\in\mathfrak{h}$. In addition, we define $V_\lambda:=\{v\in V\mid\pi(H)v=\lambda(H)v\text{ for all } H\in\mathfrak{h}\}.$ We call $V_\lambda$ the \emph{weight space} if $V_\lambda\neq\{0\}$.
\end{definition}
\begin{lemma}
	Let $(\pi,V)$ be a finite dimensional representation of $\mathfrak{g}$. Then every weight $\lambda$ is integral.
\end{lemma}
\begin{proof}
	Let $\alpha\in\Delta(\mathfrak{g},\mathfrak{h})$. Then we can find $X_\alpha\in\mathfrak{g}_\alpha$ and $Y_\alpha\in\mathfrak{g}_{-\alpha}$ such that the subalgebra generated by $X_\alpha,Y_\alpha$ and $ H_\alpha$ is isomorphic to $\mathfrak{sl}(2,\C)$, see for example \cite[Thm. 7.19]{Hall} or \cite[Lm. 6.42]{Kirillov}. That is to say, $[H_\alpha,X_\alpha]=2X_\alpha$, $[H_\alpha,Y_\alpha]=-2Y_\alpha$ and $[X_\alpha,Y_\alpha]=H_\alpha$. Let us denote this subalgebra as $\mathfrak{sl}(2,\C)_\alpha$. Then $(\pi|_{\mathfrak{sl}(2,\C)_\alpha},V)$ is a representation of $\mathfrak{sl}(2,\C)_\alpha$, and we only have to show that $\lambda(H_\alpha)\in\Z$. However, the representation theory of $\mathfrak{sl}(2,\C)$ is studied extensively, and it is known that if $(\pi|_{\mathfrak{sl}(2,\C)_\alpha},V)$ is a representation of $\mathfrak{sl}(2,\C)$ (not necessarily irreducible), then all eigenvalues of $\pi(H_\alpha)$ are integers. For a proof, see for example \cite[Thm. 4.34]{Hall}. In particular, $\lambda(H_\alpha)\in\Z$.
\end{proof}

Since the weights are nothing else than eigenvalues, we immediately see that we have the decomposition of weight spaces in the same way as we would get for the eigenvalues of a linear operator, i.e., $$V=\bigoplus_{\lambda \text{ is a weight}} V_\lambda.$$ Of course, not all of the weight spaces need to be one-dimensional, as eigenvalues can occur with a multiplicity of $\geq 2$. This motivates the following definition:

\begin{definition}\label{def:weight_multiplicity}
	Let $(\pi,V)$ be a finite dimensional representation of $\mathfrak{g}$. We define the \emph{weight multiplicity} of $\lambda$ as $\dim(V_\lambda).$ We call $(\pi,V)$ \emph{weight multiplicity free} if the weight multiplicity is exactly 1 for all weights in $(\pi,V)$.
\end{definition}

\begin{theorem}\cite[Thm. 9.3]{Hall}\label{thm:weights_invariant_under_W}
	If $(\pi,V)$ is a finite dimensional representation of $\mathfrak{g}$, the weights and weight multiplicities of $\pi$ are invariant under the action of $W$.
\end{theorem}

\begin{example}\label{ex:weights_SU(N)_standard_and_dual_repr}
	Consider the example $SU(N)$ again. We continue Example \ref{ex:root_system_SU(N)}. Let us consider the representation $(\Pi_\mathrm{std},\C^N)$ of $\mathfrak{sl}(N,\C)$. Denote by $e_1,\ldots,e_N$ the standard basis of $\C^N$. We note that for $H_{\alpha_1}=\mathrm{diag}(1,-1,0\ldots,0)$ we get that $$H_{\alpha_1}e_1=e_1,\qquad H_{\alpha_j}e_1=0$$ for all $j\geq 2$. In other words, $\C e_1\subseteq (\C^N)_{\omega_1}$. In fact, $$H_{\alpha_1}e_2=- e_2,\qquad H_{\alpha_1}e_j=0$$ for all $j\geq 3$. In other words, the weight space of $\omega_1$ is exactly $\C e_1$. Similarly we see that $$H_{\alpha_1}e_2=-e_2,\qquad H_{\alpha_2}e_2=e_2,\qquad H_{\alpha_j}e_2=0,$$ which is just saying $(\C^N)_{\omega_1-\alpha-1}=\C e_2$. Continuing this for $e_3$ etc., we get
	\begin{align*}
		\C^N &= \C e_1\oplus \C e_2\oplus\cdots\oplus \C e_N\\
		&= (\C^N)_{\omega_1}\oplus (\C^N)_{\omega_1-\alpha_1}\oplus\cdots\oplus  (\C^N)_{\omega_1-\alpha_1-\cdots-\alpha_{N-1}}.
	\end{align*}
	So the weights are $\omega_1,\omega_1-\alpha_1,\ldots,\omega_1-\alpha_1-\cdots-\alpha_{N-1}$. Note that all the weight multiplicities are 1, and all weights lie in the orbit of $\omega_1$, as predicted by Theorem \ref{thm:weights_invariant_under_W}.  
	
	Similarly, one could consider $(\Pi_\mathrm{std}^*,(\C^N)^*)$. Here $((\C^N)^*)_{\omega_{N-1}}=\C e_{N}$, and going through the same steps, one sees that $\omega_{N-1},\omega_{N-1}-\alpha_{N-1},\ldots,\omega_{N-1}-\alpha_{N-1}-\cdots-\alpha_{1}$ are the weights.
\end{example}

Examining the last example, it is evident that among the weights there is only one dominant weight. For $(\Pi_\mathrm{std},\C^N)$, the dominant weight is $\omega_1$, and for $(\Pi_{\mathrm{std}}^*,(\C^N)^*)$ it is $\omega_{N-1}$. 

\begin{definition}
	Let $(\pi,V)$ be a finite dimensional representation of $\mathfrak{g}$. We say that a weight $\lambda\in\mathfrak{h}^*$ is the \emph{highest weight of $\pi$} if for every other weight $\mu$ we have $\mu\preceq \lambda$. We call any non-trivial $v\in V_{\lambda}$ a \emph{highest weight vector} of $V$.
\end{definition}

Note that in the Example \ref{ex:weights_SU(N)_standard_and_dual_repr}, $\omega_1$ and $\omega_2$ are the highest weights. However, this is not exceptional behavior, as the following theorem will tell us.

\begin{theorem}\cite[Thm. 9.4]{Hall},\cite[Thm. in Ch. 20.2]{Humphreys}\label{thm:every_irreduc_repr_has_highest_weight}
	Let $(\pi,V),(\rho,W)$ be a finite dimensional irreducible representation of $\mathfrak{g}$. Let $\Delta$ be a root system of $\mathfrak{g}$, choose a set of positive roots $\Delta^+$ with corresponding simple roots $\Phi$. Then
	\begin{enumerate}[label=(\roman*)]
		\item $(\pi,V)$ has a highest weight $\lambda$.
		\item Every other weight is of the form $\lambda-\sum_{i=1}^lk_i\alpha_i$ where $\alpha_i\in\Phi$ and $k_i\in\N_0$, i.e. every weight is $\preceq\lambda$.
		\item If $(\rho,W)$ is another finite dimensional irreducible representation with the same highest weight, then $\pi\simeq \rho.$
		\item If $\lambda$ is the highest weight of $(\pi,V)$, then $\lambda$ is dominant integral, i.e $\lambda\in P^+$.
	\end{enumerate}
\end{theorem}

In other words, for every irreducible representation $(\pi,V)$, we can find a unique weight that is both dominant integral, i.e. $\mu\in P^+$, and the highest weight as well. However, we can also go the other way, which will be the Theorem of Highest Weight. Before doing so, we first prove a small lemma.

\begin{lemma}\label{lemma:highest_weight_iff_cyclic_vector}
	Let $(\pi,V)$ be an irreducible finite dimensional representation of $\mathfrak{g}$. Then $(\pi,V)$ has highest weight $\lambda$ if and only if $V$ is generated by a vector $v\in V_\lambda$ such that $\pi(X)v=0$ for all $X\in\bigoplus_{\alpha\in\Delta^+}\mathfrak{g}_\alpha$.
\end{lemma}
\begin{proof}
	Let us first consider $\lambda$ being the highest weight of $(\pi,V).$ Take $v\in V_\lambda$. It is clear that $v$ generates $V$, as any non-zero vector generates $V$. So consider $X\in\mathfrak{g}_\alpha$ for any $\alpha\in\Delta^+$, and assume $\pi(X)v\neq 0$. Then we note that 
	\begin{align*}
		\pi(H)\pi(X)v=\pi([H,X])v+\pi(X)\pi(H)v &= \pi(\alpha(H)X)v+\lambda(H)\pi(X)v\\
		&=(\lambda+\alpha)(H)\pi(X)v.
	\end{align*}
	In other words, $\lambda+\alpha$ is a weight. But $\lambda\preceq\lambda+\alpha$, which is a contradiction.
	
	Next, let us consider that $V$ is being generated by $v\in V_\lambda$ with $\pi(X)v=0$. Before proving the statement, choose a set of positive roots $\Delta^+=\{\alpha_1,\ldots,\alpha_n\}$ and choose a basis for $\mathfrak{g}$ consisting of $Y_i\in \mathfrak{g}_{-\alpha_i}, X_i\in\mathfrak{g}_{\alpha_i}$ and $H_1,\ldots,H_r\in\mathfrak{h}$ a basis of $\mathfrak{h}$. We can extend any representation $\pi$ uniquely to a representation of the universal enveloping algebra $U\mathfrak{g}$, which we will also denote by $\pi$. By the Poincar\'e-Birkhoff-Witt the monomials $Y_1^{k_1}\cdots Y_n^{k_n}H_1^{l_1}\cdots H_r^{l_r}X_1^{m_1}\cdots X_n^{m_n}$ form a basis of $U\mathfrak{g}$.
	
	Coming back to the representation, notice that $$\pi(H)\pi(Y_i)v = (\lambda-\alpha_i)(H)\pi(Y_i)v$$ for all $i$. So $\lambda-\alpha$ is a weight if $\pi(Y)v\neq 0$. By the same argument the vector $\pi(Y_1)^{k_1}\cdots\pi(Y_n)^{k_n}v$ has weight $\lambda-k_1\alpha_1-\ldots -k_n\alpha_n$. Of course, $\pi(H_i)v=\lambda(H_i)$ and $\pi(X_i)v=0$. Since $V$ is generated by $v$, it must be spanned by elements of the form $\pi(Y_1)^{k_1}\cdots\pi(Y_n)^{k_n}v$ with weight $\lambda-k_1\alpha_1-\cdots-k_n\alpha_n$. Note that there are only finitely many choices for $k_1,\ldots,k_n$ to get the weight $\beta$, and by finite dimensionality of $V$ one cannot admit all $k_i\in\N_0$. Therefore, we have found all the weights, with $\lambda$ being the highest weight.
\end{proof}

The above lemma can be generalized, which leads to the theory of Verma modules. We will, however, skip the theory of Verma modules, and immediately go to one of the main results. For more information on Verma modules, we refer to \cite{Hall,Humphreys,Kirillov}.

\begin{theorem}\cite[Thm. 8.18]{Kirillov}\label{thm:Verma_modules_to_irreduc_repr}
	Let $\lambda\in \mathfrak{h}^*$. Then there exists a unique, up to isomorphism, irreducible representation $(\pi_\lambda,V(\lambda))$ of $\mathfrak{g}$ with highest weight $\lambda$.
\end{theorem}

\begin{theorem}\cite[Thm. in 21.2]{Humphreys}
	Let $\lambda\in\mathfrak{h}^*$ be dominant integral, i.e. $\lambda\in P^+$. Then the irreducible representation $(\pi_\lambda,V(\lambda))$ of $\mathfrak{g}$ in Theorem \ref{thm:Verma_modules_to_irreduc_repr} is finite dimensional.
\end{theorem}

So given a $\lambda\in P^+$ we have an irreducible representation $(\pi_\lambda,V_\lambda)$. We already had the other implication, as seen in Theorem \ref{thm:every_irreduc_repr_has_highest_weight}. Putting these together in a theorem, we get:

\begin{theorem}[Theorem of the Highest Weight]\label{thm:highest_weight_decomposition}
	Let $\mathfrak{g}$ be a complex semisimple Lie algebra. The map sending $\lambda\mapsto (\pi_\lambda,V(\lambda))$ is a bijection from $P^+$ to the set of equivalent classes of finite dimensional irreducible representations of $\mathfrak{g}$.
\end{theorem}
As a result, we can classify all possible finite dimensional irreducible representations of a complex semisimple Lie algebra by looking at $P^+$. 

\begin{example}\label{ex:SU(N)_representations_weights}
	We revisit our example of $SU(N)$, for which $\mathfrak{g}=\mathfrak{sl}(N,\C)$. Consider the trivial representation $(\Pi_\mathrm{triv},\C)$ of $\mathfrak{sl}(\C)$. This is clearly an irreducible representation, where $\Pi_\triv(X)v=0$ for all $v\in\C$. The highest weight is $0\in\mathfrak{h}$, which is dominant integral. Next, we saw in Example \ref{ex:weights_SU(N)_standard_and_dual_repr} that $(\Pi_\mathrm{std},\C^N)$ is an irreducible representation of $\mathfrak{sl}(N,\C)$, and $\omega_1$ is the highest weight. This shows that $$(\Pi_\mathrm{std},\C^N)\simeq (\pi_{\omega_1},V(\omega_1)).$$ In the same way, going to the dual representation $(\Pi_\mathrm{std}^*,(\C^N)^*)$, we see that $\omega_{N-1}$ is the highest weight, so $$(\Pi_\mathrm{std}^*,(\C^N)^*)\simeq (\pi_{\omega_{N-1}},V(\omega_{N-1})).$$\end{example}
Continuing Example \ref{ex:SU(N)_representations_weights}, we have other irreducible representations we wish to consider for later applications.
\begin{lemma}\label{lemma:tensor_prod_of_irreducible_descends_to_irr_symm_alg}
	Consider the tensor product representation $(\bigotimes^k\Pi_\mathrm{std},T^k\C^N)$ of the standard representation of $\mathfrak{sl}(N,\C)$. Then the canonical representation $(\pi,S^k\C^N)$, where $\pi$ is the resulting map of factoring $\bigotimes^k\Pi_\mathrm{std}$ through to the symmetric power $S^k\C^N$, is an irreducible representation for all $k\in\N$. In particular, $$(\pi,S^k\C^N)\simeq (\pi_{k\omega_1},V(k\omega_1)),$$ where $\omega_1$ is the fundamental weight as in Example \ref{ex:SU(N)_representations_weights}.
\end{lemma}
\begin{proof}
	Before describing the representation on $S^k\C^N$, let us first consider $(\Pi_\mathrm{std},\C^N)$ more. We consider the matrix $E_{ij}\in\mathfrak{g}$ with $i\neq j$, where $(E_{ij})_{\mu,\nu}=\delta_{i,\mu}\delta_{j,\nu}$. We see that $$\Pi_{\mathrm{std}}(E_{ij})e_m=E_{ij}e_m=\delta_{j,m}e_i.$$ Now consider the tensor product $T^k\C^N$. By definition, the tensor product representation is given by $$\pi(X)(x_1\otimes \cdots\otimes x_k):=\Pi_{\mathrm{std}}(X)x_1\otimes x_2\otimes\cdots x_n+\cdots+x_1\otimes \cdots\otimes \Pi_{\mathrm{std}}(X)x_n$$ for any $x_1,\ldots,x_n\in\C^N$ and $X\in\mathfrak{sl}(N,\C)$. This factors through to the symmetric algebra $S^k\C^N$. In particular $$\pi(E_{ij})e_1^{m_1}\cdots e_n^{m_n}=m_je_1^{m_1}\cdots e_{j-1}^{m_{j-1}}e_{j}^{m_j-1}e_{j+1}^{m_{j+1}}\cdots e_{i-1}^{m_{i-1}}e_i^{m_i+1}e_{i+1}^{m_{i+1}}\cdots e_n^{m_n},$$ where $m_j\in\N_0$ such that $\sum_{j=1}^n m_j=k$. We claim that the monomial $e_1^k$ generates $S^k\C^N$. Since all monomials generate $S^k\C^N$, we just have to show that all monomials are generated from $e_1^k$. Note $$\pi(E_{i1})^{m_1}e_1^k=\frac{k!}{m_1!}e_1^{k-m_1}e_i^{m_1}$$ and $$\pi(E_{j1})^{m_2}\pi(E_{i1})^{m_1}e_1^k=\frac{k!}{(m_1+m_2)!}e_1^{k-m_1-m_2}e_i^{m_1}e_j^{m_2}$$ and so forth. So all monomials can be found by letting $\mathfrak{g}$ act on $e_1^k$, which is another way of stating that $e_1^k$ generates $S^k\C^N$. Note that $$\pi(H_\alpha)e_1^k=ke_1^{k-1}(\pi(H_\alpha)e_1)=k\omega_1(H_\alpha)e_1^k.$$ In addition, $E_{ij}\in\mathfrak{g}_{\alpha_i-\alpha_j}$ satisfies $\pi(E_{ij})e_1^{k}=0$ for all $i>j$. So $\pi\left(\bigoplus_{\alpha\in\Delta^+}\mathfrak{g}_\alpha\right)e_1^k=0$. In other words, $k\omega_1$ is the highest weight with highest weight vector $e_1^k$. By the proof of Lemma \ref{lemma:highest_weight_iff_cyclic_vector}, it shows that $S^k\C^N$ is irreducible with highest weight $k\omega_1$.
\end{proof}
Similarly, we have the following lemma.
\begin{lemma}\label{lemma:tensor_prod_of_dual_descends_to_irr_symm_alg}
	Let $(\Pi_{\mathrm{std}}^*,(\C^N)^*)$ be the dual of the standard representation of $\mathfrak{sl}(N,\C)$. Then the representation on $S^k(\C)^*$ is also irreducible and $$(\pi,S^k(\C^N)^*)\simeq(\pi_{k\omega_{N-1}},V(k\omega_{N-1})).$$
\end{lemma}
\begin{proof}
	Since $\Pi_\mathrm{std}^*(E_{ij})=-E_{ji}$, the proof of Lemma \ref{lemma:tensor_prod_of_irreducible_descends_to_irr_symm_alg} shows that $(\pi,S^k(\C^N)^*)$ is irreducible. Going through the same steps with $(e_N)^*$ instead of $e_1$, where $(e_N)^*$ is defined as $(e_N)^*(e_j)=\delta_{N,j}$, we see that $V(k\omega_{N-1})\simeq S^k(\C^N)^*.$
\end{proof}
Taking Lemma \ref{lemma:tensor_prod_of_irreducible_descends_to_irr_symm_alg} into account, we will denote the representation acting irreducibly on $S^k\C^N$ by $\pi_{k\omega_1}$. In the same way, we will denote $\pi_{k\omega_{N-1}}$ as the representation acting irreducibly on $S^k(\C^N)^*$.

As an aside, we note that if $\lambda\in P^{+}$, then $w_0\lambda$ lies in the negative Weyl chamber. Here $w_0$ is the longest Weyl group element. In particular, we see that $-w_0\lambda$ lies again in $P^+$. Since the dual representation $((\pi_{\lambda})^*,V(\lambda)^*)$ has as lowest weight $-\lambda$, it shows that $-w_0\lambda$ is the highest weight of $((\pi_\lambda)^*,V(\lambda)^*)$. In other words $$((\pi_\lambda)^*,V(\lambda)^*)\simeq (\pi_{-w_0\lambda},V(-w_0\lambda))$$

The correspondence between irreducible representations and dominant integral weights is just the beginning of the analysis of finding irreducible subrepresentations. We will discuss some results on finding irreducible subrepresentations in the tensor product of two irreducible representations in Section \ref{sec:tensor_prod} using highest weights; however, some results can be gathered already.

\begin{proposition}\cite[Prop. in 21.3]{Humphreys}\label{prop:weights_are_W_invariant_and_saturated}
	Let $\lambda\in P^+$. Any $\mu\in\mathfrak{h}^*$ is a weight of $(\pi_\lambda,V(\lambda))$ if and only if $\nu\preceq \lambda$ for all $\nu\in W(\mu)$.
\end{proposition}
More precisely, all the weights lie in the convex hull of the orbit of the highest weight under $W$, and if the weight $\mu$ appears, then also $\mu-k\alpha$ appears, where $k\in\N_0$ with $k\leq \mu(H_\alpha)$, for all $\alpha\in\Delta$ such that $(\mu,\alpha)>0$ \cite{Hall}.

Other impactful results are the Weyl character formula (see \cite{KnappRepresentation,Kirillov,Hall}), Freudenthal's formula (see \cite{Humphreys}) or Kostant's Multiplicity formula (see \cite{Hall}). We will, however, just mention a corollary of Weyl's character formula, which is the Weyl dimension formula.

\begin{theorem}[The Weyl dimension formula]\cite[Thm. 10.18]{Hall}\label{thm:Weyl_dim_formula}
	Let $(\pi_\lambda,V(\lambda))$ be a finite dimensional irreducible representation of $\mathfrak{g}$ with highest weight $\lambda$. Then  
	\begin{equation}\label{eq:Weyl_dimension_formula}
		\dim(V(\lambda))=\frac{\prod_{\alpha\in\Delta^+}(\alpha,\lambda+\rho)}{\prod_{\alpha\in\Delta^+}(\alpha,\rho)} = \frac{\prod_{\alpha\in\Delta^+}(\lambda+\rho)(H_\alpha)}{\prod_{\alpha\in\Delta^+}\rho(H_\alpha)}
	\end{equation} where $\rho=\frac{1}{2}\sum_{\alpha\in\Delta^+}\alpha$.  
\end{theorem}

\subsubsection{Irreducible representations of a compact group}\label{sec:highest_weight_repr_compact_grp}

After the excursion into the world of complex semisimple Lie algebras, we return to the continent of compact connected Lie groups. We recall that we are mostly interested in the $K=SU(N)$ case. Thus, by Lemma \ref{lemma:representation_G_and_Lie_alg_isomorphic_simply_connected}, we also have an immediate characterization of all irreducible representations of $K$:

\begin{theorem}[Theorem of the Highest Weight for simply connected compact groups]\label{thm:highest_weight_semisimple_groups}
	Let $K$ be a simply connected compact Lie group, and let $T$ be a maximal torus in $K$. Let $\Delta=\Delta(\mathfrak{k}_\C,\mathfrak{t}_\C)$ be the corresponding root system with a choice of positive roots $\Delta^+$. Then the set $P^{+}$ is in bijection to the set of inequivalent finite dimensional irreducible representations $(\pi,V)$ of $K$, where the correspondence is that $\lambda\in P^{+}$ is the highest weight of $\pi$.
\end{theorem}
\begin{proof}
	Let us first consider a finite dimensional irreducible representation $(\pi,V)$ of $K$. As $(\pi,V)$ is a representation of $K$, the derivative $(T_e\pi,V)$ is an irreducible representation of $\mathfrak{k}$, and hence an irreducible representation of $\mathfrak{k}_\C$, which is semisimple because $\mathfrak{k}$ is. By Theorem \ref{thm:highest_weight_decomposition} there corresponds a highest weight to this representation, let us call this $\lambda\in (\mathfrak{t}_\C)^*$. Note that $\lambda$ is dominant integral.
	
	On the other hand, let $\lambda\in\mathfrak{t}_\C^*$ be a dominant integral weight. So by Theorem \ref{thm:highest_weight_decomposition} we have an irreducible representation $(\pi_\lambda,V(\lambda))$ of $\mathfrak{k}_\C$. Restricting $\pi_\lambda$ to $\mathfrak{k}$ gives an irreducible representation of $\mathfrak{k}$, and since $K$ is simply connected, this defines an irreducible representation on $K$.
	
	To see that two representations with highest weight $\lambda$ are equivalent, let $(\pi,V)$ and $(\pi',W)$ be irreducible representations of $K$ with highest weight $\lambda\in\mathfrak{t}_\C^*$. However, if the representations $(T_e\pi,V)$ and $(T_e\pi',W)$ have the same highest weight, then by Theorem \ref{thm:highest_weight_decomposition} they are isomorphic. Since $K$ is connected, this means that also $(\pi,V)$ and $(\pi',W)$ are isomorphic. 
\end{proof}

\begin{example}
	Let us go back to the case of $SU(N)$, and specifically Example \ref{ex:SU(N)_representations_weights}. Note that the set $P^+=\bigoplus_{j=1}^{N-1}\N_0\omega_j$ gives all inequivalent finite dimensional representations of $SU(N)$. We saw that $(\pi_{k\omega_1},S^{k}\C^N)$ and $(\pi_{k\omega_{N-1}},S^k(\C^N)^*)$ are irreducible representations of $\mathfrak{su}(N)$, hence by Theorem \ref{thm:highest_weight_semisimple_groups} these are also irreducible representations of $SU(N)$.
\end{example}

\subsection{Complexifying groups}\label{sec:complexification}

Now that we have classified the irreducible representations of a simply connected compact Lie group, one could be interested in how this impacts other non-compact groups as well. We have already seen that the representation theories of $\mathfrak{k}$ and $\mathfrak{k}_\C$ are the same, so one could consider whether there exists a group that is ``the complexification'' of the group $K$, whatever that means. There are also other arguments to consider complexifying a group, for all compact Lie groups are defined over $\R$, as we have seen in Lemma \ref{lemma:compact_grp_complex_is_finite}, and $\R$ is clearly not algebraically closed, which algebraic group theory often requires. In short, we are interested in a complexification. This section will be dedicated to the notion of complexifying a compact Lie group.

\begin{definition}\label{def:complexification_of_Lie_grp}
	Let $K$ be a Lie group. We define an \emph{analytic complexification} of $K$, or just the \emph{complexification} of $K$, as a complex analytic group $G$ together with a Lie group homomorphism $i:K\rightarrow G$ such that whenever $\phi:K\rightarrow H$ is a Lie group homomorphism into a complex analytic Lie group $H$, there exists a unique analytic Lie group homomorphism $\Phi:G\rightarrow H$ such that $$\phi=\Phi\circ i.$$ We denote $K_\C:=G$.
\end{definition}
We note that there also exist other ways of complexifying a group, some more algebraically inclined or more Lie theoretic, and not all of them are necessarily equivalent. See for example \cite{Bump,KnappBeyond}. We will however use Definition \ref{def:complexification_of_Lie_grp} as it is the one that is most commonly used in the literature.

\begin{example}\label{ex:SU(n)_complexification_is_SL(n,C)}
	\begin{itemize}
		\item Consider the circle group $U(1)=\{z\in\C\mid|z|=1\}$. Then $U(1)_\C=\C^\times$.
		\item Let $n\in\N$, and consider the compact Lie group $SU(N)$. We have already seen that $$\mathfrak{su}(N)_\C=\mathfrak{sl}(N,\C)=\{ A\in M(N,\C)\mid\Tr(A)=0\}$$ which is the Lie algebra of $SL(N,\C)$. Since $SL(N,\C)$ is simply connected, any homomorphism $\psi:\mathfrak{sl}(N,\C)\rightarrow\mathfrak{h}$ comes from a unique $\Psi:SL(N,\C)\rightarrow H$. In other words, if we have a homomorphism $\phi:SU(N)\rightarrow H$, then $T_e\phi:\mathfrak{su}(N)\rightarrow \mathfrak{h}$ is a homomorphism that can be extended to the map $\psi:\mathfrak{sl}(N,\C)\rightarrow \mathfrak{h}$ with $\psi|_{\mathfrak{su}(N)}=T_e\phi$, which again comes from a homomorphism $\Psi:SL(N,\C)\rightarrow H$, and $\Psi|_{SU(N)}=\phi$. So $SL(N,\C)=SU(N)_\C$, thus extending $\phi$, making $SL(N,\C)$ the complexification of $SU(N)$.
	\end{itemize}
\end{example}
\noindent Note that this definition is just the universal property that defines the Lie group $G$ uniquely up to isomorphism. This definition is particularly nice, for now any representation $(\pi,V_\pi)$ of $K$ immediately lift to a representation $(\tilde{\pi},V_\pi)$ of $K_\C$, and $\tilde{\pi}|_{i(K)}$ will give the representation on $K$ back. However, with this definition, it is not clear whether $K_\C$ exists. However, the following result proves this for the cases we are interested in.
\begin{theorem}\cite[Thm. 24.1]{Bump}\label{thm:complexification}
	Let $K$ be a compact connected Lie group. Then $K$ has an analytic complexification $K\rightarrow K_\C$. In addition $\mathrm{Lie}(K_\C)=\mathfrak{k}_\C$, where $\mathfrak{k}_\C=\mathfrak{k}\otimes \C$ is the complexification of $\mathfrak{k}$, and the fundamental group of $K_\C$ is isomorphic to the fundamental group of $K$.
	In addition, any faithful complex representation of $K$ can be extended to a faithful analytic complex representation of $K_\C$
\end{theorem}
Interestingly, there are multiple ways to construct the analytical complexification, such as an algebraic way. To do this, one considers $K$ as an algebraic group (which can always be done, see for example \cite{Procesi}). The process of finding the complexification through algebraic means goes by the name of the Tannaka-Krein duality. We will not discuss it, but for more details regarding the Tannaka-Krein duality and how compact Lie groups are related to algebraic groups, we refer to \cite{Broecker_tom_Dieck,Procesi}. 

Next, one would be interested in the representation theory of $K_\C$. It should be noted that $K_\C$ is not compact, and therefore not all irreducible representations are necessarily finite dimensional. For example, for $K=SU(2)$, we have $K_\C=SL(2,\C)$, and there exist irreducible representations of $SL(2,\C)$ that are not finite dimensional, such as the principal series. See for example \cite[Prop. 2.6]{KnappRepresentation}. However, when we restrict ourselves to finite dimensional ones, we have a one-to-one correspondence with the irreducible representations of $K$.

\begin{proposition}\label{prop:irred_K_one-to-one_with_irred_K_C}
	Let $K$ be a connected compact Lie group and $K_\C$ its complexification. If $(\pi,V)$ is an irreducible representation of $K$, then $(\pi_\C,V)$ is a finite dimensional irreducible holomorphic representation of $K_\C$, where $\pi_\C$ is the lift of $\pi$ to $K_\C$. In addition, if $(\pi,V)$ is a finite dimensional irreducible holomorphic representation of $K_\C$, then $(\pi|_{K},V)$ is an irreducible representation of $K$.
\end{proposition}
\begin{proof}
	Let $(\pi,V)$ be an irreducible representation of $K$. Consider $(\pi_\C,V)$ as representation of $K_\C$. If $\pi_\C(K_\C)W\subseteq W$ for some non-trivial subspace $W\subseteq V$, then for sure $\pi(K)W=\pi_\C(K)W\subseteq W$. However, $(\pi,V)$ is irreducible, which is a contradiction. So $(\pi_\C,V)$ is also an irreducible representation of $K_\C$.
	Now consider a finite dimensional irreducible holomorphic representation $(\pi,V)$ of $K_\C$. Consider $(\pi|_{K},V)$ as representation of $K$. Assume to the contrary that there exists a non-trivial subspace $W$ such that $\pi(K)W\subseteq W$. We remark that $$\pi(\exp(X))=\exp(T_e\pi(X))$$ for any $X\in\mathfrak{k}$. In particular, since $W$ is closed, this means that $T_e\pi(X)W\subseteq W$, so $(T_e\pi,V)$ is a reducible representation of $\mathfrak{k}$. The map $T_e\pi$ can be extended to a complex linear map, hence defines a representation $(T_e\pi,V)$ of $\mathfrak{k}_\C$, which is still reducible for $W$ is a vector space over $\C$. However, $T_e\pi$ is already complex linear, and the fact that $W$ is closed gives that $(\pi,V)$ must be reducible, which is a contradiction. So $(\pi|_{K},V)$ is an irreducible representation of $K$. 
\end{proof}

\begin{corollary}\label{cor:highest_weight_decomposition_complexification}
	Let $K=SU(N)$, and $K_\C=SL(N,\C)$ be the complexification of $K$. Then the finite dimensional irreducible representations of $K_\C$ are in a one-to-one correspondence with the dominant integral weights of $\mathfrak{k}$, which in turn are in a one-to-one correspondence with the irreducible representations of $K$.
\end{corollary}
\begin{proof}
	Use Proposition \ref{prop:irred_K_one-to-one_with_irred_K_C} together with Theorem \ref{thm:highest_weight_semisimple_groups}.
\end{proof}

\subsection{The power of being minuscule}
After having discussed the Theorem of Highest Weight and some theory about the complexification of a group, we require the notion of minuscule weights. We mostly follow \cite{Bourbaki_ch7_8} here. Let $\mathfrak{g}$ be a complex semisimple Lie algebra with Cartan subalgebra $\mathfrak{h}$. Let $\Delta(\mathfrak{g},\mathfrak{h})$ be the root system, and $\Phi$ the simple roots of $\mathfrak{g}$.
\begin{definition}
	Let $\lambda\in P^+$ be non-zero. We call $\lambda$ \emph{minuscule} if it is minimal with respect to the partial order $\preceq$ on $P^+$. Minimal here means that there exists no $\alpha\in P^+$ such that $\alpha\prec \lambda$.
\end{definition}

\begin{example}
	\begin{enumerate}[label=(\roman*)]
		\item Consider $SU(3)$ with the root system as before. The weight lattice is given in Figure \ref{fig:weight_space_A_2}. Just by looking at the picture, we see that the minuscule weights are given by $\omega_1$ and $\omega_2$. 
		\begin{figure}
			\begin{subfigure}[t]{0.48\textwidth}
				\centering
				\begin{tikzpicture}[scale=0.6]
					\foreach \angle in {30,90,...,330} {
						\draw[->, thick, black] (0,0) -- ({2.6*cos(\angle)}, {2.6*sin(\angle)});
					}
					\begin{rootSystem}[weight length=1.5cm,weight radius=1.5pt]{A}
						\node [above] at \Root {1}{0} {\(\alpha_1\)};
						\node [right] at \Root {0}{1} {\(\alpha_2\)};
						\node [right] at \Root {1}{-1} {\(\)};
						\fundamentalweights
						\node [right] at \weight {1}{0} {\(\omega_2\)};
						\node [right] at \weight {0}{1} {\(\omega_1\)};
						\WeylChamber
					\end{rootSystem}
				\end{tikzpicture}
				\caption{The root system for $A_2$. The weight lattice is given by the corners of every triangle, and the minimal non-zero dominant integral weights are given by $\omega_1$ and $\omega_2$.}
				\label{fig:weight_space_A_2}
			\end{subfigure}
			\hspace{0.2cm}
			\begin{subfigure}[t]{0.48\textwidth}
				\centering
				\begin{tikzpicture}[scale=0.6]
					\draw[step=1cm,gray,very thin] (-3.9,-3.9) grid (3.9,3.9);	
					\fill[gray!60, opacity=0.5] (0,0) rectangle (3.9,3.9);\foreach \angle in {1,2} {
						\draw[->, thick, black] (0,0) -- ({2*cos(90*(\angle-1))}, {2*sin(90*(\angle-1))}) node[anchor=south west] {$\alpha_\angle$};
					}
					\foreach \angle in {3,4} {
						\draw[->, thick, black] (0,0) -- ({2*cos(90*(\angle-1))}, {2*sin(90*(\angle-1))});}
					
					\fill (1,0) circle (2pt) node[anchor=south east] {$\omega_1$};   
					\fill (0,1) circle (2pt) node[anchor=south east] {$\omega_2$};
					\fill (1,1) circle (2pt) node[anchor=south west] {$\omega_1+\omega_2$};
				\end{tikzpicture}
				\caption{The root system for $A_1\times A_1$ inside its weight lattice, where the weights are given by the corners of every square. The fundamental weights are given by $\omega_1$ and $\omega_2$. Note that the minuscule weights are given by $\omega_1,\omega_2$ and $\omega_1+\omega_2$.}
				\label{fig:weight_space_A_1xA_1}
			\end{subfigure}
		\end{figure}
		\item We now consider the case $A_1\times A_1$. See Figure \ref{fig:weight_space_A_1xA_1} for the weight lattice. Here the minuscule weights are clearly $\omega_1,\omega_2$ and $\omega_1+\omega_2$.
	\end{enumerate}
\end{example}

Some authors use different definitions of a minuscule weight, see for example \cite{Bourbaki_ch4_6,Stembridge,Mathieu}. However, most are equivalent which we will show now.

\begin{lemma}\label{lemma:repres._equiv_def_minuscule_weight}
	Let $\lambda$ be a non-zero dominant integral weight. Then the following are equivalent:
	\begin{enumerate}[label=(\roman*)]
		\item $\lambda$ is minuscule,
		\item the Weyl group $W$ acts transitively on the set of weights of $V(\lambda)$,
		\item $\lambda$ is the unique dominant integral weight for $V(\lambda)$.
	\end{enumerate}
\end{lemma}

\begin{proof}
	$\textrm{(i)}\iff \textrm{(ii)}$ is clear, as the weights of $V(\lambda)$ are exactly the set $\cup_{\mu\in P^+, \mu\preceq \lambda}W(\mu)$ by Proposition \ref{prop:weights_are_W_invariant_and_saturated}. If there is only one such $\mu\preceq \lambda$, namely $\mu=\lambda$, the weights are exactly $W(\lambda)$.
	
	$\textrm{(ii)}\implies \textrm{(iii)}.$ Let $\mu$ be another dominant integral weight of $V(\lambda)$. By assumption, there exists $s\in W$ such that $\mu=s(\lambda)$. However, the Weyl group orbit $W(\lambda)$ intersects the set of dominant weights only at one point, which is at $\lambda$. Therefore $s=e$ and thus $\lambda=\mu$.
	
	$\textrm{(iii)}\implies \textrm{(i)}$. Let $\lambda$ be the unique dominant integral weight for $V(\lambda)$. This means that if $\Phi=\{\alpha_i\}_i$, then $\lambda-\alpha_i$ is not a dominant integral weight anymore, for all $i$. Hence, $\lambda$ is minimal with respect to $\preceq$, so minuscule.
\end{proof}

Lemma \ref{lemma:repres._equiv_def_minuscule_weight} shows that there are not that many options for minuscule weights. Indeed, there are only finitely many options, and when considering simple Lie algebras, they are subsets of the fundamental weights. To show that, we require one lemma.

\begin{lemma}\label{lemma:repres._equiv_def_minuscule_weight_2}
	Let $\lambda$ be a dominant integral weight. Then $\lambda$ is minuscule weight if and only if $\lambda(H_\alpha)\in\{-1,0,1\}$ for all roots $\alpha\in\Delta$.
\end{lemma}
\begin{proof}
	Assume $\lambda$ is minuscule. By definition, it follows that $\lambda-\alpha_i$ is not dominant for every simple root $\alpha_i\in\Phi$. In particular this means that there exists a $j_i$ such that $$(\lambda-\alpha_i)(H_{\alpha_{j_i}})<0.$$ But recall that $\alpha_i(H_{\alpha_j})\leq 0$ if $i\neq j$, thus previous inequality can only be satisfied if $j_i=i$. Now $\alpha_i(H_{\alpha_i})=2$ by definition, which means that $$(\lambda-\alpha_i)(H_{\alpha_i})=\lambda(H_{\alpha_i})-2<0.$$ By definition of dominant we have $\lambda(H_{\alpha_i})\geq 0$, hence $\lambda(H_{\alpha_i})=0$ or $1$ for all $i$. Note that $H_{-\alpha_i}=-H_{\alpha_i}$, so $\lambda(H_\alpha)=-1,0$ or $1$ for all roots $\alpha$. 
	Now assume $\lambda(H_\alpha)\in\{-1,0,1\}$. We note that $$(\lambda-\alpha_i)(H_{\alpha_i})=\lambda(H_{\alpha_i})-2<0$$ for all $i$. In other words, $\lambda-\alpha_i$ is not dominant, for all $i$, so $\lambda$ is minimal with respect to $\preceq$, hence minuscule.
\end{proof}

\begin{corollary}
	Let $\mathfrak{g}$ be a simple Lie algebra, and let $\lambda\in P^+$ be minuscule. Then $\lambda$ is a fundamental weight, see Definition \ref{def:fundamental_weight}.
\end{corollary}

\begin{proof}
	Let $\lambda$ be minuscule. Since it is dominant integral, we can expand it in the basis of fundamental weights, i.e., $\lambda=n_1\omega_1+\cdots + n_k\omega_k$ where $\omega_j$ are the fundamental weights, and $n_j\in \N_0$. Note $\lambda(H_{\alpha_j})\in\{0,1\}$ for all simple roots by Lemma \ref{lemma:repres._equiv_def_minuscule_weight_2}, but also $\lambda(H_{\alpha_j})=n_j$. So $n_j=0$ or 1. To get the result, we recall that the coroots $\{H_\alpha\mid \alpha\in\Delta\}$ form a root system $R^\vee$, with simple roots $\{H_\alpha\mid \alpha\in\Phi\}$ \cite{Bourbaki_ch4_6}. Then there is also a highest root $H_\mu$ in $R^\vee$, which can be written as $$H_\mu=m_1H_{\alpha_1}+m_2H_{\alpha_2}+\cdots m_kH_{\alpha_{k}}$$ with $m_j\geq 0$. But recall $\lambda(H_\mu)\in\{-1,0,1\}$, so $$n_1m_1+\cdots+n_km_k\in\{-1,0,1\}$$ while $n_j=0$ or 1. This can only happen if there is an $i$ such that $n_i=1$ and $n_j=0$ for all $j\neq i$.
\end{proof}

\begin{example}\label{ex:fundamental_weights_A_N_minuscule}
	For the root system $A_n$, all fundamental weights are minuscule; see for example \cite[p. 132]{Bourbaki_ch4_6}. This need not happen, for in $B_2$, only the fundamental weight corresponding to the short root is a minuscule weight, while in $G_2$, there are no minuscule weights.
\end{example}

\subsection{Decomposing tensor products}\label{sec:tensor_prod}

The last few results of representation theory that we need to discuss Mathieu's work are in relation to decomposing tensor products into irreducibles. The lemmas in this and upcoming section were discussed in \cite{Mathieu}, however the proofs are novel. 

This section is based on \cite{Hall, Kumar, Wallach}. In this section, we will consider tensor products of irreducible representations. To be more specific, let $K$ be a compact connected Lie group, and let $G=K_\C$ be its complexification; see Section \ref{sec:complexification}. Let $(\pi,V)$ and $(\rho,W)$ be irreducible representations of $K$, which lift to irreducible representations of $G$ as well by Proposition \ref{prop:irred_K_one-to-one_with_irred_K_C}. Since any finite dimensional representation of $K$ is completely reducible, see Corollary \ref{cor:compact_group_completely_reducible}, we have $$V\otimes W\simeq \bigoplus_{j=1}^n V_j$$ as representations of $K$, where $(\pi_i,V_i)$ are irreducible representations of $K$ with $\pi_i(k)=\pi(k)|_{V_i}$ for all $k\in K$. In this section, we will be interested in which irreducible representations $(\pi_i,V_i)$ of $K$ will occur in the above decomposition and with what multiplicity (recall Definition \ref{def:multiplicity}). Since the correspondence between irreducible representations of $K$ and finite dimensional irreducible representations of $G$ is one-to-one, we will focus mostly on the representations of $G$ because the Lie algebra is over $\C$. However, note that the difference here is minimal by Proposition \ref{prop:irred_K_one-to-one_with_irred_K_C}, so sometimes we will switch back to $K$ to use certain properties that are only available in the compact case (such as having a left and right Haar measure, characters, etc.). The theory of decomposing tensor products is rich. For a more in depth discussion, we refer to \cite{Kumar} for a more general decomposition.

To start the discussion, we need to set the stage. Let $K$ be a simply connected compact Lie group and $G=K_\C$. Let $\mathfrak{k}=\mathrm{Lie}(K)$ and $\mathfrak{g}:=\mathfrak{k}_\C=\mathrm{Lie}(G)$. We follow the strategy of Section \ref{sec:Root_system_compact_group}. We choose a maximal torus $T\subseteq K$ with Lie algebra $\mathfrak{t}$, and let $\mathfrak{h}:=\mathfrak{t}_\C$ be the Cartan subalgebra corresponding to $\mathfrak{t}$, which dictates a root system $\Delta(\mathfrak{g},\mathfrak{h})$. Choose a set of positive roots $\Delta^+(\mathfrak{g},\mathfrak{h})$, and let $\Phi$ be the set of simple roots.

Let us denote the set of equivalence classes of finite dimensional irreducible representations of $G$ as $\widehat{G}$. By Corollary \ref{cor:highest_weight_decomposition_complexification}, we can identify $\widehat{G}$ with $\widehat{K}$ which again can be identified with the set of dominant integral weights of $\mathfrak{h}$, which we shall denote by $P^{+}$. For a weight $\lambda\in P^{+}$ we shall denote the highest weight representation of $K$ corresponding to $\lambda$ by $(\pi_\lambda,V(\lambda))$.

We note that $-w_0\lambda$ is the highest weight of $V(\lambda)^*$, the dual representation of $V(\lambda)$. Here, $w_0$ is the longest Weyl group element acting on $\lambda$ in the usual way. We shall write $\lambda^*:=-w_0\lambda$ for the highest weight of $V(\lambda)^*$. Finally, we shall denote by $\triv$ the trivial representation of $K$, which corresponds to the highest weight $0\in P^{+}$.


\begin{remark}
	Before we proceed, we would like to take a moment to point out a notation we will use. 
	By Theorem \ref{thm:highest_weight_semisimple_groups} and Corollary \ref{cor:highest_weight_decomposition_complexification}, there is a one-to-one correspondence between the finite dimensional irreducible representations of $K$, of $G$, and of $\mathfrak{g}$ with integral highest weight. Because of this, we will use the same letters for the representation of $K$, $G$ and $\mathfrak{g}$. That is to say, when $(\pi_\lambda,V(\lambda))$ is a representation of $K$ with $\lambda\in P^{+}$, we also use the same notation for the representation of $G$ or $\mathfrak{g}$, and vice versa. We will always note which group or algebra we are using for the representation but note that the same letters will be used, given the one-to-one correspondence. This notation will be kept throughout this paper.
\end{remark}
\begin{remark}
	Since we are on the topic of notations, we will discuss some tensor products. Giving a name for every representation often dilutes the proofs. Therefore, when it is clear what the representation should be, we will omit the letter representing the representation. So for example if $(\pi,V)$ is a representation of $G$, and $v\in V$, then we write $gv$ instead of $\pi(g)v$, with $g\in G$, to simplify our notations.
\end{remark} 

With the tools discussed at hand, we consider two lemmas, used to find the multiplicities of irreducibles occurring in the representation $(\pi_\lambda\otimes\pi_\mu,V(\lambda)\otimes V(\mu))$. Parts of the proofs are based on \cite{Kumar, Stembridge}.

\begin{lemma}\label{lemma:Mathieu_3.1}
	Let $\lambda,\mu\in P^{+}$. Then
	\begin{enumerate}[label=\roman*)]
		\item The irreducible representation $(\pi_{\lambda+\mu},V(\lambda+\mu))$ occurs with multiplicity one in $V(\lambda)\otimes V(\mu)$.
		\item Assume $\lambda-\mu^*$ is dominant. Then $(\pi_{\lambda-\mu^*},V(\lambda-\mu^*))$ occurs with multiplicity one in $V(\lambda)\otimes V(\mu)$.
		\item The trivial representation $(\pi_\triv,\C)$ appears in $V(\lambda)\otimes V(\mu)$ if and only if $\mu=\lambda^*$, and the multiplicity is one.
	\end{enumerate}
\end{lemma}
\begin{proof}
	For part $\textrm{(i)}$, we consider the representations as representations of $\mathfrak{g}$. Recall that $\mathfrak{g}$ acts on $V(\lambda)\otimes V(\mu)$ as $X(v\otimes w)=\pi_\lambda(X)v\otimes w+ v\otimes \pi_\mu(X)w$ where $X\in \mathfrak{g}$. Considering the highest weight vectors $v^+\in V(\lambda)$ and $w^+\in V(\mu)$, we get
	$$H(v^+\otimes w^+)=(\lambda(H)+\mu(H))(v^+\otimes w^+)\qquad \forall H\in\mathfrak{h},$$ and $X(v^+\otimes w^+)=0$ for all $X\in\bigoplus_{\alpha\in\Delta^+}\mathfrak{g}_\alpha$, which shows that $\lambda+\mu$ is the highest weight in $V(\lambda)\otimes V(\mu)$. In addition, $v^+\otimes w^+$ generates a $\mathfrak{g}$-invariant subspace isomorphic to $V(\lambda+\mu)$. Since $v^+\otimes w^+$ is the only vector with weight $\lambda+\mu$, the multiplicity is 1, proving \textrm{(i)}.
	
	\textrm{(ii)}: Here we see the representations as representations of $K$. Especially, we make us of Lemma \ref{lemma:Repr_theory_review_character_results}. Let $\chi_{\lambda-\mu^*}$ be the character of $V(\lambda-\mu^*)$, and $\chi_{\lambda\otimes \mu}$ be the character of $V(\lambda)\otimes V(\mu)$. Then
	\begin{align*}
		\langle \chi_{\lambda\otimes \mu},\chi_{\lambda-\mu^*}\rangle &= \int_K \chi_{\lambda\otimes \mu}(k)\overline{\chi_{\lambda-\mu^*}(k)}dk\\
		&=\int_K \chi_{\lambda}(k)\chi_{\mu}(k)\overline{\chi_{\lambda-\mu^*}(k)}dk\\
		&=\int_K \chi_{\lambda}(k)\overline{\chi_{\mu^*}(k)}\overline{\chi_{\lambda-\mu^*}(k)}dk\\
		&=\int_K \chi_{\lambda}(k)\overline{\chi_{\mu^*\otimes(\lambda-\mu^*)}(k)}dk\\
		&=\overline{\langle \chi_{\mu\otimes (\lambda^*-\mu)},\chi_{\lambda^*}\rangle} = 1
	\end{align*}
	where in the last equation we used part $(\textrm{i})$ since $\lambda^*=\mu+(\lambda^*-\mu)$.
	
	\textrm{(iii):} This is proven by character calculation as well:
	\begin{align*}
		\langle \chi_{\lambda\otimes \mu},\chi_\triv\rangle &=\int_K \chi_{\lambda\otimes\mu}(k)\overline{\chi_{\triv}(k)}dk = \int_K \chi_{\lambda\otimes\mu}(k)dk\\
		&=\int_K \chi_{\lambda}(k)\chi_\mu(k)dk = \int_K \chi_{\lambda}(k)\overline{\chi_{\mu^*}(k)}dk = \delta_{\lambda,\mu^*}
	\end{align*}
	where $\delta_{\lambda,\mu^*}=1$ if $\lambda=\mu^*$ and zero otherwise.
\end{proof}
For the next result, we recall the definition of a weight multiplicity of a weight $\lambda$, given in Definition \ref{def:weight_multiplicity}, which is just the dimension of $V_\lambda$. We called $(\pi,V)$ weight multiplicity free if the weight multiplicity is exactly one.
\begin{lemma} \cite[Lemma 3.2]{Mathieu}\label{lemma:Mathieu_3.2}
	Let $\lambda,\mu,\nu\in P^{+}$.
	\begin{enumerate}[label=\roman*)]
		\item Assume $V(\mu)$ is weight multiplicity free. Then $V(\lambda)\otimes V(\mu)$ is multiplicity free. In addition, if $V(\nu)$ appears in $V(\lambda)\otimes V(\mu)$, then $\nu-\lambda$ is a weight of $V(\mu)$.
		\item If $\mu$ is minuscule, then $V(\nu)$ appears in $V(\lambda)\otimes V(\mu)$ if and only if $\nu-\lambda$ is a weight of $V(\mu)$.
	\end{enumerate}
\end{lemma}
\begin{proof}
	The proof is based on \cite{Kumar}. For part $(\textrm{i})$, note that since $V(\nu)$ is finite dimensional, we have by Lemma \ref{lemma:highest_weight_iff_cyclic_vector} that there exists a $v_\nu^+\in V$ such that $\pi(H)v_\nu^+=\nu(H)v_\nu^+$ and $\pi(X)v_\nu^+=0$ for all $H\in\mathfrak{h}$ and $X\in\bigoplus_{\alpha\in\Delta^+}\mathfrak{g}_\alpha$. In particular, any linear map $T$ commuting with the action of $\mathfrak{g}$ is completely defined by what $T(v_\nu^+)$ is. Using Schur's Lemma, we thus see that
	\begin{align*}
		\langle \chi_{\lambda\otimes \mu},\chi_\nu\rangle&=\dim_\C \mathrm{Hom}_\mathfrak{g}(V(\nu),V(\lambda)\otimes V(\mu))\\
		&=\dim_\C \mathrm{Hom}_\mathfrak{b}(\C v_\nu^+,V(\lambda)\otimes V(\mu))\\
		&=\dim_\C \mathrm{Hom}_\mathfrak{b}(\C v_\nu^+\otimes V(\lambda)^*, V(\mu)),
	\end{align*}
	where $\mathfrak{b}=\mathfrak{h}\oplus \bigoplus_{\alpha\in\Delta^+}\mathfrak{g}_\alpha$. If we assume $V(\nu)$ actually appears in $V(\lambda)\otimes V(\mu)$ we have $\langle \chi_{\lambda\otimes \mu},\chi_\nu\rangle\geq 1$. So, the right-hand side must be at least one. In other words, there exists a non-zero linear map $T\in \mathrm{Hom}_\mathfrak{b}(\C v_\nu^+\otimes V(\lambda)^*,V(\mu))$. Note that $\C v_\nu^+\otimes V(\lambda)^*$ is generated by $v_\nu^+\otimes v_{-\lambda}^-$, where $v_{-\lambda}^-\in V(\lambda)^*$ is the lowest weight vector generating $V(\lambda)^*$ by $Hv_{-\lambda}^-=-\lambda(H)v_{-\lambda}^-$ for all $H\in\mathfrak{h}$ and $Xv_{-\lambda}^-=0$ for all $X\in\bigoplus_{\alpha\in\Delta^+}\mathfrak{g}_{-\alpha}$. This shows that $$HT(v_\nu^+\otimes v_{-\lambda}^-)=(\nu(H)-\lambda(H))T(v_\nu^+\otimes v_{-\lambda}^-).$$
	Since $T$ is not zero, $T(v_\nu^+\otimes v_{-\lambda}^-)$ cannot be zero. In other words, we see that $\nu-\lambda$ is a weight in $V(\mu)$. In fact, any $T\in \mathrm{Hom}_\mathfrak{b}(\C v_\nu^+\otimes V(\lambda^*),V(\mu))$ is completely determined by $T(v_\nu^+\otimes v_{-\lambda}^-)$, so the dimension can be at most 1. So $\langle \chi_{\lambda\otimes \mu},\chi_\nu\rangle\leq 1$ proving \textrm{(i)}.
	
	For a proof of \textrm{(ii)}, we first claim that $V(\mu)$ is weight multiplicity free. Since $W$ acts transitively by Lemma \ref{lemma:repres._equiv_def_minuscule_weight}, any weight $\mu$ that appears in $V(\lambda)$ satisfies $$\dim(V(\lambda)_\mu)=\dim(V(\lambda)_{w(\mu)})=\dim(V(\lambda)_\lambda)=1.$$ So $V(\mu)$ is weight multiplicity free. Using part $\textrm{(i)}$ we find $\implies$. So we only need to show $\impliedby$. We refer to \cite[Corollary (3.4)]{Kumar} to show that it is enough to show that if $\mu$ is minuscule, then $\lambda+w\mu$ is nearly dominant, which is to say $(\lambda+w\mu)(H_{\alpha_i})\geq -1$ for all simple roots $\alpha_i$ and all $w\in W$. Note that it is enough to prove that $w\mu(H_{\alpha_i})\geq -1$, since $\lambda(H_{\alpha_i})\geq 0$ by definition. We now remark that $$w\mu(H_{\alpha_i}) = \frac{2(w\mu,\alpha_i)}{(\alpha_i,\alpha_i)}=\frac{2(\mu,w^{-1}\alpha_i)}{(w^{-1}\alpha_i,w^{-1}\alpha_i)}=\mu(H_{w^{-1}\alpha_i}).$$ Using Lemma \ref{lemma:repres._equiv_def_minuscule_weight_2} proves the lemma. 
\end{proof}

\section{Application to $SU(N)$ and $SL(N,\C)$}\label{sec:prerequisites_Mathieu}

Having covered a significant quantity of prerequisites, we turn to the application to $SU(N)$. We already applied most of the theory to $SU(N)$, but the decomposition of tensor products into certain subrepresentations of $SU(N)$ is essential in Mathieu's proof. For this decomposition, we need certain equivariant maps (see Definition \ref{def:intertwiner}). These maps will be used in the proof of Mathieu, which is discussed in Section \ref{chapter:Mathieu_implies_Jacobian}. The proofs in this section are novel as well.

This section continues with the one example we have been looking at, namely $SU(N)$, with complexification $SL(N,\C)$ (see Example \ref{ex:SU(n)_complexification_is_SL(n,C)}). We continue from Example \ref{ex:SU(N)_representations_weights}. We already saw that $V(k\omega_1)\simeq S^k\C^N$ in Lemma \ref{lemma:tensor_prod_of_irreducible_descends_to_irr_symm_alg} and that $V(k\omega_{N-1})\simeq S^k(\C^N)^*$ by Lemma \ref{lemma:tensor_prod_of_dual_descends_to_irr_symm_alg}. We also know that $\omega_1$ and $\omega_{N-1}$ are minuscule weights, see Example \ref{ex:fundamental_weights_A_N_minuscule}. With an eye on the previous section, we get the following lemma.

\begin{lemma}\label{lemma:Mathieu_tensor_product_=_direct_sum}
	Let $N\in\N$. Then $$S^k\C^N\otimes (\C^N)^*\simeq V(k\omega_1+\omega_{N-1})\oplus V((k-1)\omega_1)$$ as representation of $SU(N)$, of $SL(N,\C)$, of $\mathfrak{su}(N)$ and of $\mathfrak{sl}(N,\C)$.
\end{lemma}
\begin{proof}
	By Corollary \ref{cor:highest_weight_decomposition_complexification}, the finite dimensional representations of $SU(N)$, $SL(N,\C)$ and $\mathfrak{su}(N)$ are the same, so we can choose which picture to use. We will prove the lemma by seeing all representations as representations of $\mathfrak{sl}(N,\C)$. Using Lemma \ref{lemma:tensor_prod_of_irreducible_descends_to_irr_symm_alg} and \ref{lemma:tensor_prod_of_dual_descends_to_irr_symm_alg} we get $$V:=S^k\C^N\otimes (\C^N)^*\simeq V(k\omega_1)\otimes V(\omega_{N-1}).$$ Making use of Lemma \ref{lemma:Mathieu_3.1}, we notice that $V(k\omega_1+\omega_{N-1})$ occurs in $V$ with multiplicity one. In addition, by Example \ref{ex:SU(N)_representations_weights}, $\omega_1^*=\omega_{N-1}$, hence $$k\omega_1-\omega_{N-1}^*=k\omega_1-\omega_1=(k-1)\omega_1,$$ which is dominant. By Lemma \ref{lemma:Mathieu_3.1} we can conclude that also $V((k-1)\omega_{1})$ occurs in $V$ with multiplicity one. So we see that the space $$W:=V(k\omega_1+\omega_{N-1})\oplus V((k-1)\omega_1)$$ occurs in $V$.
	Assume to the contrary that there exists another $\nu\in P^{+}$ such that $V(\nu)$ occurs in $V$ with multiplicity $m\geq 1$. Now note that $\omega_{N-1}$ is minuscule, so by Lemma \ref{lemma:Mathieu_3.2}, $V(\nu)$ appears in $V$ if and only if $\nu-k\omega_1$ is a weight in $V(\omega_{N-1})$. The weights of $(\pi_{\omega_{N-1}},V(\omega_{N-1}))\simeq (\Pi_{\mathrm{std}}^*,(\C^N)^*)$ have already been calculated, namely in Example \ref{ex:weights_SU(N)_standard_and_dual_repr}, and are given by: $$\omega_{N-1},\,\omega_{N-1}-\alpha_{N-1},\,\ldots,\,\omega_{N-1}-\sum_{j=1}^{N-1}\alpha_{N-j}.$$ So the options for $\nu$ are: $$\nu=k\omega_1+\omega_{N-1},\quad\nu=k\omega_1+\omega_{N-1}-\alpha_{N-1},\quad\ldots,\quad  \nu=k\omega_1+\omega_{N-1}-\sum_{j=1}^{N-1}\alpha_{N-j}.$$  We recall from Equation (\ref{eq:omega_1_explicit}) and (\ref{eq:omega_{n-1}_explicit}) that
	\begin{align}
		\omega_1&=\frac{1}{N}((N-1)\alpha_1+(N-2)\alpha_2+\cdots +\alpha_{N-1}),\\
		\omega_{N-1}&=\frac{1}{N}(\alpha_1+2\alpha_2+\cdots+(N-1)\alpha_{N-1}).
	\end{align}
	and thus we see that only $\nu=k\omega_1+\omega_{N-1}$ and $\nu=k\omega_{1}+\omega_{N-1}-\sum_{j=1}^{N-1}\alpha_{N-j} = k\omega_1-\omega_1=(k-1)\omega_1$ are actually dominant integral weights. The other options give negative integers when one considers $\nu(H_{\alpha_{j}})$, where $j$ is the subscript in the last simple root in the possible options for $\nu$. Since we already know that $V(k\omega_1+\omega_{N-1})$ and $V((k-1)\omega_1)$ appear with multiplicity one, that means that no other $V(\nu)$ can occur in $V$. 
\end{proof}

\begin{proof}[Alternative proof]
	For an alternative proof, we again note that $W\subseteq V$. We argue $W=V$ through a dimensional argument. We have that $\dim(S^mV)=\binom{\dim(V)+m-1}{m}$ for any finite dimensional vector space $V$. So $$\dim(S^k\C^N\otimes (\C^N)^*)=N\frac{ (N+k-1)!}{k!(N-1)!}.$$ On the other hand, we have $W$ which has dimension
	\begin{align}
		\dim(W)&=\dim(V(k\omega_1+\omega_{N-1}))+\dim(V((k-1)\omega_1))\nonumber\\
		&= \dim(V(k\omega_1+\omega_{N-1}))+\dim(S^{k-1}\C^N)\nonumber\\
		&= \dim(V(k\omega_1+\omega_{N-1}))+\binom{N+k-2}{k-1}\label{eq:alternative_proof}.
	\end{align}
	To find $\dim(V(k\omega_1+\omega_{N-1}))$, we make use of Weyl's dimension formula, given in Theorem \ref{thm:Weyl_dim_formula}. For that, we recall that $\Delta^+=\{\epsilon_{i}-\epsilon_j\mid i\geq j\}.$ We can write $\epsilon_i-\epsilon_j = \alpha_i+\alpha_{i+1}+\cdots+\alpha_{j}$, see for example \cite{Bourbaki_ch4_6}. With that, we have 
	\begin{align*}
		\Delta^+=\{&\alpha_1,\alpha_2,\ldots,\alpha_{N-1},\\
		&\alpha_1+\alpha_2,\alpha_1+\alpha_2+\alpha_3,\ldots,\alpha_1+\cdots+\alpha_{N-2},\alpha_1+\cdots+\alpha_{N-1},\\
		&\alpha_2+\alpha_3,\ldots,\alpha_2+\cdots+\alpha_{N-1},\ldots,\alpha_{N-2}+\alpha_{N-1}\}.
	\end{align*}
	It is known that for any simple root $\alpha_j$, we have $\frac{(\rho,\alpha_j)}{(\alpha_j,\alpha_j)}=1$, see for example \cite[Lemma 7.36]{Kirillov} or \cite[Prop. 2.69]{KnappBeyond}. Multiplying Equation (\ref{eq:Weyl_dimension_formula}) with 1, where $1=\frac{\prod_{\alpha\in\Delta^+}(\alpha,\alpha)}{\prod_{\alpha\in\Delta^+}(\alpha,\alpha)}$, the denominator can be written as (evaluating the product at $\alpha\in \Delta^+$ in the same order as in the above)
	\begin{align*}
		\prod_{\alpha\in\Delta^+}\frac{(\alpha,\rho)}{(\alpha,\alpha)}=&1\cdot1\cdot\cdots\cdot 1\cdot2\cdot3\cdot\cdots (N-2)\cdot(N-1)\cdot 2\cdots (N-2)\cdot\cdots 2\\
		=&(N-1)!(N-2)!\cdots 2!.
	\end{align*}
	Similarly, we see that
	\begin{align*}
		\prod_{\alpha\in\Delta^+}\frac{(\alpha,k\omega_1+\omega_{N-1}+\rho)}{(\alpha,\alpha)}=& \prod_{\alpha\in\Delta^+}\left[k\frac{(\alpha,\omega_1)}{(\alpha,\alpha)}+\frac{(\alpha,\omega_{N-1})}{(\alpha,\alpha)}+\frac{(\alpha,\rho)}{(\alpha,\alpha)}\right]\\
		=&(k+0+1)\cdot1\cdot\cdots \cdot (0+1+1)\cdot \\
		&(k+2)\cdot (k+3)\cdot\cdots (k+N-2)\cdot (k+1+N-1)\cdot\\
		&\cdot\cdots\cdot (1+N-2)\cdot\cdots 3.\\
		=&\frac{(k+N)!}{k!(k+N-1)}\cdot\frac{(N-1)!(N-2)!\cdots 2!}{(N-2)!} .
	\end{align*}
	Hence 
	\begin{align*}
		\dim(V(k\omega_1+\omega_{N-1}))&=\prod_{\alpha\in\Delta^+}\frac{(\alpha,k\omega_1+\omega_{N-1}+\rho)}{(\alpha,\rho)}\\
		&=\frac{(k+N)!}{k!(N-2)!(k+N-1)}=(k+N)(N-1)\frac{(k+N-2)!}{k!(N-1)!}.
	\end{align*}
	Putting that back in, into Equation $(\ref{eq:alternative_proof})$ we get:
	\begin{align*}
		\dim(W)&= (k+N)(N-1)\frac{(k+N-2)!}{k!(N-1)!}+\frac{(k+N-2)!}{(k-1)!(N-1)!}\\
		&= \frac{[(k+N)(N-1)+k](k+N-2)!}{k!(N-1)!}\\
		&= \frac{[N(N+k-1)](k+N-2)!}{k!(N-1)!}=\dim(V),
	\end{align*}
	which proves the statement.
\end{proof}

We now turn to the representation $(\pi_{k\omega_1}\otimes \pi_{l\omega_{N-1}},S^k\C^N\otimes S^l(\C^N)^*)$ of $SU(N)$, or equivalently $SL(N,\C)$, and see how it decomposes into irreducibles. Before describing the representation in more detail, we first construct the space itself. Let $x_1,x_2,\ldots,x_N$ be an orthonormal basis of $\C^N$. We identify $S\C^N$ as the set of polynomials on $(\C^N)^*$, and we can identify $S(\C^N)^*$ with the algebra of symmetric differential operators on $S\C^N$. In other words, we identify $S\C^N=\C[x_1,\ldots,x_N]$ and $S(\C^N)^*=\C\left[\partial_1,\ldots,\partial_N\right]$ where $\partial_i(x_j)=\delta_{ij}$. This way, we can identify $S^k\C^N$ as the algebra of homogeneous polynomials of degree $k$, and $S^l(\C^N)^*$ can be identified with the algebra of homogeneous linear partial differential operators of degree $l$.

We recall that the action of $\mathfrak{sl}(N,\C)$ on $S^k\C^N$ is given in Lemma \ref{lemma:tensor_prod_of_irreducible_descends_to_irr_symm_alg}, which is to say that
\begin{align*}
	\pi_{k\omega_1}(E_{ij})x_1^{m_1}\cdots x_n^{m_n}=m_jx_1^{m_1}\cdots x_{j-1}^{m_{j-1}}x_{j}^{m_j-1}x_{j+1}^{m_{j+1}}\cdots x_{i-1}^{m_{i-1}}x_i^{m_i+1}x_{i+1}^{m_{i+1}}\cdots x_n^{m_n},
\end{align*}
where $m_j\in\N_0$ with $\sum_{j=1}^nm_j=k$, and
$$\pi_{k\omega_1}(E_{ii}-E_{jj})x_1^{m_1}\cdots x_n^{m_n}=(m_i-m_j)x_1^{m_1}\cdots x_n^{m_n}.$$ Here the matrix $(E_{ij})_{\mu,\nu}=\delta_{i,\mu}\delta_{j,\nu}$. As an action of $SL(N,\C)$, we have $$\left(\otimes^k\pi_{\mathrm{std}}\right)(g)(y_1\otimes\cdots\otimes y_{k})=(\pi_{\mathrm{std}}(g)y_1)\otimes\cdots\otimes (\pi_{\mathrm{std}}(g)y_k)$$ for any $y_i\in \C^N$ and thus it factors through to an action on $S^k\C^n$ as $$\pi_{k\omega_1}(g)(x_1^{m_1}\cdots x_k^{m_k})=(\pi_{\mathrm{std}}(g)x_1)^{m_1}\cdots(\pi_{\mathrm{std}}(g)x_k)^{m_k}.$$ Similar results hold for the dual representation $(\pi_{k\omega_{N-1}},S^k(\C^N)^*)$, where one should take note that $\pi_{k\omega_{N-1}}(E_{ij})=(-1)^k\pi_{k\omega_1}(E_{ji})$.

The rest of this section is dedicated to building three $SL(N,\C)$-equivariant maps to get some insight into which irreducible representations occur in the tensor representation $S^k\C^N\otimes S^l(\C^N)^*$. These maps will be needed in Section \ref{chapter:Mathieu_implies_Jacobian}.

\begin{definition}
	Let $k,l\in\N$. Define the operator
	\begin{align*}
		\mathrm{div}:S^k\C^N\otimes S^l(\C^N)^* &\rightarrow S^{k-l}\C^N\\
		\sum_{i} P_i\otimes Q_i&\mapsto \sum_i Q_i(\partial_1,\ldots,\partial_n)(P_i)
	\end{align*}
	if $k\geq l$. If $k<l$ then $S^{k-l}\C^N=\{0\}$ by definition, so $\mathrm{div}=0$. Here we denoted $Q_i(\partial_1,\ldots,\partial_n)\in\C[\partial_1,\ldots,\partial_n]$ such that $Q_i(\partial_1\ldots,\partial_n)(P_i)$ is the polynomial of differential operators acting on the polynomial $P_i$. Note in the case of $k=l$, we get the Fischer inner product (for more information, see, for example, \cite{Faraut_Koranyi_1990,Arazy_1992}).
\end{definition}
\begin{lemma}\label{lemma:Mathieu_div_G_intertwiner}
	The map $\mathrm{div}$ is a $SL(N,\C)$-equivariant map for all $k,l$.
\end{lemma}
\begin{proof}
	First, we note that if $k<l$, then the map $\mathrm{div}=0$, which is trivially an intertwiner. So, let us consider $k\geq l$. We start by noting that $\mathrm{div}$ can be defined on the tensor algebra in the following way
	\begin{align*}
		\mathrm{div}:T^k \C^N\otimes T^l(\C^N)^*&\rightarrow T^{k-l}\C^N\\
		x_{i_1}\otimes\cdots \otimes x_{i_k}\otimes \partial_{j_1}\otimes \cdots\otimes\partial_{j_l}&\mapsto \\
		\frac{1}{k!l!}\sum_{\sigma\in S_k}\sum_{\rho\in S_l}\big(\partial_{\rho(j_1)}&(x_{\sigma(i_1)})\cdots\partial_{\rho(j_l)}(x_{\sigma(i_l)})\big)\,x_{\sigma(i_{l+1})}\otimes\cdots\otimes x_{\sigma(i_{k})}
	\end{align*}
	and linearly extend it. This definition maps, by construction, into $S^{k-l}\C^N$, and descends to the map $\mathrm{div}:S^k\C^N\otimes S^l(\C^N)^*\rightarrow S^{k-l}\C^N$ we defined before. But now note that $\mathrm{div}$ is $SL(N,\C)$-equivariant, for 
	\begin{align*}
		\mathrm{div}(g(&x_1\otimes \cdots\otimes x_k\otimes \partial_1\otimes \cdots\otimes \partial_l))= \mathrm{div}(gx_1\otimes \cdots\otimes gx_k\otimes g\partial_1\otimes \cdots\otimes g\partial_l)\\
		&=\frac{1}{k!l!}\sum_{\sigma\in S_k}\sum_{\rho\in S_l}\left(g\partial_{\rho(1)}(gx_{\sigma(1)})\cdots g\partial_{\rho(l)}(gx_{\sigma(l)})\right)\,gx_{\sigma(l+1)}\otimes\cdots\otimes gx_{\sigma(k)}.
	\end{align*}
	We recall that the representation of $SL(N,\C)$ on $\partial_j$ is given by the dual representation, so $$(g\partial_i)(gx_j)=\partial_i(g^{-1}gx_j)=\partial_i(x_j).$$ This gives that
	\begin{align*}
		\mathrm{div}(g(x_1&\otimes \cdots\otimes x_k\otimes \partial_1\otimes \cdots\otimes \partial_l))\\
		&=\frac{1}{k!l!}\sum_{\sigma\in S_k}\sum_{\rho\in S_l}\left(\partial_{\rho(1)}(x_{\sigma(1)})\cdots \partial_{\rho(l)}(x_{\sigma(l)})\right)\,gx_{\sigma(l+1)}\otimes\cdots\otimes gx_{\sigma(k)}\\
		&=g\circ \mathrm{div}(x_1\otimes \cdots\otimes x_k\otimes \partial_1\otimes \cdots\otimes \partial_l),
	\end{align*}
	which means $g\circ \mathrm{div}=\mathrm{div}\circ g$. So the map $\mathrm{div}:T^k\C^N\otimes T^l(\C^N)^*\rightarrow T^{k-l}\C^N$ is a $SL(N,\C)$-intertwiner. Restricting it to the map $\mathrm{div}:S^k\C^N\otimes S^l(\C^N)^*\rightarrow S^{k-l}\C^N$ is then also a $SL(N,\C)$-intertwiner.
\end{proof}

\begin{corollary}\label{cor:Mathieu_div_is_an_isomorphism_onto_subspace}
	The multiplicity of the representation $V((k-l)\omega_1)$ occurring in $S^k\C^N\otimes S^l(\C^N)^*$ is exactly one, and $\mathrm{div}$ is an isomorphism on the $V((k-l)\omega_1)$-component and zero otherwise. 
\end{corollary}
\begin{proof}
	We recall that $SL(N,\C)=SU(N)_\C$, and $SU(N)$ is compact, so restricting the representation $(\pi_{k\omega_1}\otimes \pi_{l\omega_{N-1}},S^k\C^N\otimes S^l(\C^N)^*)$ to $SU(N)$ gives a representation that is completely reducible by Lemma \ref{cor:compact_group_completely_reducible}. In other words, $$S^k\C^N\otimes S^l(\C^N)^*\simeq\oplus_{\lambda\in P^{+}}d_\lambda V(\lambda)$$ as representations of $SU(N)$, where $d_\lambda$ is the multiplicity of $V(\lambda)$ in $S^k\C^N\otimes S^l(\C^N)^*$. Note that by Corollary \ref{cor:highest_weight_decomposition_complexification} the above decomposition also holds as representations of $SL(N,\C)$. To prove the first part, note that by Lemma \ref{lemma:Mathieu_3.1}\textrm{.(ii)}, the representation $V((k-l)\omega_1)$ occurs with multiplicity 1 in $S^k\C^N\otimes S^l(\C^N)^*$ since $$S^k\C^N\otimes S^l(\C^N)^*\simeq V(k\omega_1)\otimes V(l\omega_{N-1}).$$ To prove the second part, by Lemma \ref{lemma:Mathieu_div_G_intertwiner}, the map $\mathrm{div}:S^k\C^N\otimes S^l(\C^N)^*\rightarrow S^{k-l}\C^N\simeq V((k-l)\omega_1)$ is a non-zero $SL(N,\C)$-intertwiner, so by Schur's Lemma the map must be an isomorphism as map $\mathrm{div}|_{V((k-l)\omega_1)}$, and $\mathrm{div}|_{V(\lambda)}=0$ for any other $\lambda\in P^{+}$. 
\end{proof}

\begin{definition}
	Define the map $E:S^k\C^N\rightarrow (\C^N)^*\otimes S^{k+1}\C^N$ given by $$E(f)=\sum_{i=1}^N\partial_i\otimes(f\cdot x_i).$$
\end{definition}
\begin{remark}
	Consider the space $C^\infty(\R^n)$, and define the maps:
	\begin{align*}
		\Delta(f)&:=\sum_{i=1}^n\frac{\partial^2}{\partial x_i^2}f,\\
		E(f)&:=\sum_{i=1}^n x_i\frac{\partial}{\partial x_i}f.
	\end{align*}
	Then the maps $\{-\mathrm{id}-2E, \Delta,\frac{1}{2}\sum_i x_i^2\}$ form a $\mathfrak{sl}(2)$ triple, where the commutator is given by $[A,B]=A\circ B-B\circ A$. The map $E$ is called the Euler map. In other words, our definition of $E$ can be seen as a generalization of the Euler map. For some information on the Euler map, see for example \cite[pg. 182]{Faraut}.
\end{remark}
\begin{lemma}\label{lemma:Mathieu_Euler_map_G_intertwiner}
	The map $E:S^k\C^N\rightarrow (\C^N)^*\otimes S^{k+1}\C^N$ is $SL(N,\C)$-equivariant.
\end{lemma}
\begin{proof}
	We will apply the same strategy as the proof of Lemma \ref{lemma:Mathieu_div_G_intertwiner}.	We define the map $E_\mathrm{tensor}:T^k\C^N\rightarrow (\C^N)^*\otimes T^{k+1}\C^N$ by 
	\begin{align*}
		E_{\mathrm{tensor}}(x_{i_1}\otimes \cdots\otimes x_{i_k}) := \frac{1}{k+1}\sum_{j=1}^N&\big( \partial_j\otimes x_j\otimes x_{i_1}\otimes \cdots\otimes x_{i_k} + \\
		&\partial_j\otimes x_{i_1}\otimes x_j \otimes \cdots \otimes x_{i_k}+ \cdots+\\
		&\partial_j\otimes x_{i_1}\otimes \cdots x_{i_k}\otimes x_j\big),
	\end{align*}
	which, by construction of the symmetric algebra, descends to the map $E:S^k\C^N\rightarrow (\C^N)^*\otimes S^{k+1}\C^N$. Now let $X\in\mathfrak{sl}(N,\C)$. We need to show $E_{\mathrm{tensor}}\circ X = X\circ E_{\mathrm{tensor}}$. We claim that it is enough to show that 
	\begin{align}
		\sum_{j=1}^N&X(\partial_j)\otimes x_j\otimes x_1\otimes \cdots \otimes x_k + \partial_j\otimes X(x_j) \otimes x_1\otimes \cdots \otimes x_k = 0,\label{eq:Mathieu_E_intertwiner}\\
		\sum_{j=1}^N&X(\partial_j)\otimes x_1\otimes x_j\otimes \cdots \otimes x_k + \partial_j \otimes x_1\otimes X(x_j)\otimes \cdots \otimes x_k = 0,\\
		\vdots\nonumber\\
		\sum_{j=1}^N&X(\partial_j)\otimes x_1\otimes \cdots \otimes x_k \otimes x_j + \partial_j \otimes x_1\otimes \cdots \otimes x_k \otimes X(x_j)= 0\label{eq:Mathieu_E_intertwiner2}
	\end{align}	
	to show $E_{\mathrm{tensor}}$ is a $\mathfrak{sl}(N,\C)$-intertwiner. To prove this, note that if one writes out the exact action of $X$, we see that
	\begin{align*}
		(X\circ E_{\mathrm{tensor}}-&E_{\mathrm{tensor}}\circ X)(x_{i_1}\otimes \cdots\otimes x_{i_k})= \\
		\sum_{j=1}^N& X(\partial_j)\otimes x_j\otimes x_1\otimes \cdots \otimes x_k + \partial_j\otimes X(x_j) \otimes x_1\otimes \cdots \otimes x_k+\\
		&X(\partial_j)\otimes x_1\otimes x_j\otimes \cdots \otimes x_k + \partial_j \otimes x_1\otimes X(x_j)\otimes \cdots \otimes x_k +\cdots +\\
		&X(\partial_j)\otimes x_1\otimes \cdots \otimes x_k \otimes x_j + \partial_j \otimes x_1\otimes \cdots \otimes x_k \otimes X(x_j),
	\end{align*} which is exactly Equations (\ref{eq:Mathieu_E_intertwiner})-(\ref{eq:Mathieu_E_intertwiner2}).
	Each of these equations is proven in the same way, so we focus on Equation (\ref{eq:Mathieu_E_intertwiner}). Note that $\mathfrak{sl}(N,\C)$ is spanned by the matrices $E_{ij}$ for $i\neq j$ and $E_{ii}-E_{(i+1)(i+1)}$ for $i=1,\ldots, n-1$. So it is enough to prove Equation (\ref{eq:Mathieu_E_intertwiner}) with the elements $E_{ij}$.
	
	Recall that the representation of $\mathfrak{sl}(N,\C)$ is the standard representation, in other words $E_{ij}x_k=\delta_{j,k}x_i$. In the same way, $\mathfrak{sl}(N,\C)$ acts on $(\C^N)^*$ as the dual representation, i.e. $$E_{ij}\partial_k(x_l)=\partial_k\circ (-E_{ij})(x_l) = \partial_k(-\delta_{j,l}x_i)=-\delta_{j,l}\delta_{i,k} = -\delta_{i,k}\partial_j(x_l).$$ In other words, $E_{ij}\partial_k=-\delta_{i,k}\partial_j$. Then the left hand side of Equation (\ref{eq:Mathieu_E_intertwiner}) becomes 
	\begin{align*}
		\sum_{j=1}^N&E_{ab}(\partial_j)\otimes x_j\otimes x_{i_1}\otimes \cdots\otimes x_{i_k} + \partial_i\otimes E_{ab}(x_i) \otimes x_{i_1}\otimes \cdots\otimes x_{i_k} \\
		=&-\partial_b\otimes x_a\otimes x_{i_1}\otimes \cdots\otimes x_{i_k} + \partial_b\otimes x_a \otimes x_{i_1}\otimes \cdots\otimes x_{i_k} = 0. 
	\end{align*}
	which shows that $E_{\mathrm{tensor}}$, and thus $E$, is a $\mathfrak{sl}(N,\C)$-intertwiner. Now note that
	\begin{align*}
		\exp(X)\circ E=\sum_{n=1}^\infty \frac{1}{n!}X^n\circ E=\sum_{n=1}^\infty \frac{1}{n!}E\circ X^n= E\circ \exp(X)
	\end{align*}
	for all $X\in\mathfrak{sl}(N,\C)$ by linearity of $E$. By connectedness of $SL(N,\C)$ we thus have $g\circ E = E\circ g$ for all $g\in SL(N,\C)$, so $E$ is a $SL(N,\C)$-intertwiner.
\end{proof}

After having discussed two equivariant maps, we combine them to obtain a third map $\Psi$. The map $\Psi$ is the key player in proving the Jacobian Conjecture, assuming Mathieu's Conjecture. In particular, Lemma \ref{lemma:Mathieu_Psi_isomorphism} is the most important result of this section.

\begin{definition}\label{def:Mathieu_Psi}
	Define the map $\Psi:S^k\C^N\otimes S^l(\C^N)^*\rightarrow  S^{k-l+1}\C^N\otimes (\C^N)^*$ by 
	\begin{align*}
		\Psi :=& \sigma\circ(\mathrm{id}_{(\C^N)^*}\otimes \mathrm{div})\circ (E\otimes \mathrm{id}_{S^l(\C^N)^*}),\\
		\Psi&:S^k\C^N\otimes S^l(\C^N)^*\xrightarrow{E\otimes \mathrm{id}} (\C^N)^*\otimes S^{k+1}\C^N\otimes S^l(\C^N)^*\\
		&\xrightarrow{\mathrm{id}\otimes \mathrm{div}} (\C^N)^*\otimes S^{k-l+1}\C^N\xrightarrow{\sigma} S^{k-l+1}\C^N\otimes (\C^N)^*,
	\end{align*}
	where the map $\sigma:(\C^N)^*\otimes S^{k+1-l}\C^N\rightarrow S^{k+1-l}\C^N\otimes (\C^N)^*$ is defined by $$\sigma\left(\sum_i\partial_i\otimes P_i\right)=\sum_i P_i\otimes \partial_i.$$ More concretely, we define the map $\Psi$ by $$\Psi\left(\sum_{j}P_j\otimes Q_j\right)= \sum_j\sum_{i=1}^N Q_j(\partial)(P_j\cdot x_i)\otimes \partial_i,$$ where by $Q_j(\partial)(P_j\cdot x_i)$ we mean the polynomial of differential operators $Q$ with variables $\partial_1,\ldots,\partial_N$ applied to the polynomial $P_j\cdot x_i$.
\end{definition}
\begin{lemma}\label{lemma:Mathieu_Psi_isomorphism}
	Let $k,l\in\N$ be such that $k\geq l-1$. The map $\Psi:S^k\C^N\otimes S^l(\C^N)^*\rightarrow S^{k+1-l}\C^N\otimes (\C^N)^*$ is a $SL(N,\C)$-equivariant surjective map and the restrictions $\Psi|_{V((k-l)\omega_1)}$ and $\Psi|_{V((k-l+1)\omega_1+\omega_{N-1})}$ are isomorphisms upon their images.
\end{lemma}
\begin{proof}
	It is clear that $\Psi$ is a $SL(N,\C)$-intertwiner since it is a composition of three $SL(N,\C)$-intertwining maps. By Lemma \ref{lemma:Mathieu_3.2}, the image of $\Psi$ lies in $$S^{k+1-l}\C^N\otimes (\C^N)^*\simeq V((k+1-l)\omega_1+\omega_{N-1})\oplus V((k-l)\omega_1).$$ Schur's Lemma tells us then that there exists either an invariant subspace $W\subseteq S^k\C^N\otimes S^l(\C^N)^*$ such that $\Psi|_{W}$ is an isomorphism, meaning that either $V((k+1-l)\omega_1+\omega_{N-1})$ or $V((k-l)\omega_1)$, or both $V((k+1-l)\omega_1+\omega_{N-1})$ and $V((k-l)\omega_1)$ occurs in $S^k\C^N\otimes S^l(\C^N)^*$, or $\Psi|_W$ is the zero map. It should be clear that $\Psi\not\equiv 0$. We claim that there exist irreducible subspaces $W_1,W_2\subseteq S^k\C^N\otimes S^l(\C^N)^*$ such that $W_1\simeq V((k+1-l)\omega_1+\omega_{N-1})$ and $W_2\simeq V((k-l)\omega_1)$, and such that $\Psi|_{W_i}$ is an isomorphism for $i=1,2$, which is another way of saying that $\Psi$ is a surjective isomorphism. Since $\Psi$ is an intertwiner, it is enough to show that $\Psi$ maps both into the $V((k-l)\omega_1)$- and into the $V((k-l+1)\omega_1 + \omega_{N-1})$-component of $S^{k+1-l}\C^N\otimes (\C^N)^*$. Let us choose $x_1^k\otimes \partial_1^{l-1}\partial_2$. Then
	\begin{align*}
		\Psi(x_1^k\otimes \partial_1^{l-1}\partial_2) &= \sum_{i=1}^N\partial_1^{l-1}\partial_2(x_1^kx_i)\otimes \partial_i\\
		&= \partial_1^{l-1}(x_1^k)\otimes \partial_2\neq 0
	\end{align*} since $k\geq l-1$.
	Now we use Corollary \ref{cor:Mathieu_div_is_an_isomorphism_onto_subspace}, which states that $\mathrm{div}$ is an isomorphism of the $V((k-l)\omega_1)$-component of $S^{k-l+1}\C^N\otimes(\C^N)^*$ into $V((k-l)\omega_1)$. We find that  $$\mathrm{div}(\Psi(x_1^k\otimes \partial_1^{l-1}\partial_2)) = \mathrm{div}(\partial_1^{l-1}(x_1^k)\otimes \partial_2)=  0.$$ Therefore, we see that $\Psi(x_1^k\otimes \partial_1^{l-1}\partial_2)$ lies in the $V((k-l+1)\omega_1+\omega_{N-1})$-component of $S^{k-l+1}\C^N\otimes (\C^N)^*$. Since any vector in an irreducible finite dimensional representation generates the entire representation space, we can fully generate $W_1$ from $\Psi(x_1^k\otimes \partial_1^{l-1}\partial_2)$. And since $S^{k-l+1}\C^N\otimes (\C^N)^*\simeq V((k-l)\omega_1)\oplus V((k-l+1)\omega_1+\omega_{N-1})$ by Lemma \ref{lemma:Mathieu_tensor_product_=_direct_sum}, we can also generate $W_2$ by taking any other vector, such as $\Psi(x_1^k\otimes \partial_1^l)$, since $\mathrm{div}(\Psi(x_1^k\otimes \partial_1^l))\neq 0$. In other words, both restriction maps are non-zero and, therefore, isomorphisms, proving the lemma.
\end{proof}

\section{Mathieu's Conjecture implies the Jacobian Conjecture}\label{chapter:Mathieu_implies_Jacobian}
	
In this section, we will apply the representation theory tools to show that Mathieu's Conjecture implies the Jacobian Conjecture. As mentioned in Section \ref{chapter:Intro}, there is a tremendous amount of work published on the Jacobian Conjecture and major simplifications have been made. We will begin with a discussion on some of these simplifications.
	
\subsection{Some facts about the Jacobian Conjecture}\label{sec:facts_Jacobian_Conjecture}
Before actually going into details, let us state the most general version of the Jacobian Conjecture
\begin{conjecture}[The Jacobian Conjecture]\label{conj:Jacobian_conj}
	Let $n\in\N$ and let $k$ be an algebraically closed field of characteristic zero. Let $f:k^n\rightarrow k^n$ be a polynomial map, that is to say, that every $f_i:k^n\rightarrow k$ is a polynomial map, in such a way that the Jacobian $J(f)=1$. Then $f$ is invertible with polynomial inverse.
\end{conjecture}
	
The difference with Conjecture \ref{conj:Intro_Jacobian_conj} is that we allow any algebraically closed field of characteristic zero. However, there is no loss of generality to consider Conjecture \ref{conj:Intro_Jacobian_conj}, i.e., taking $k=\C$, for there is a common trick that is used in algebraic geometry, which is sometimes called the Lefschetz Principle:
\begin{lemma}[`The Lefschetz Principle' at the start of Section I in \cite{Jacobian_overzicht_artikel}]\label{lemma:Mathieu_Jacobian_enough_to_show_over_C}
	If the Jacobian Conjecture is true for the algebraically closed field $\C$, then the Jacobian Conjecture is true for any algebraically closed field $k$ of characteristic zero.
\end{lemma}
As mentioned, the Jacobian Conjecture has a rich history. For a more detailed outline of both the history and a discussion on possible ways to tackle the Jacobian Conjecture, we refer to sources like \cite{Jacobian_overzicht_artikel,vdEssen,vdEssen_en_de_rest}. Some reductions have been found, which we will briefly discuss.
	
\begin{conjecture}\cite[Introduction]{Jacobian_overzicht_artikel}\label{conj:d_restricted_Jacobian_conj}
	Let $n\in \N$, and $f:\C^n\rightarrow \C^n$ be a polynomial map with Jacobian $J(f)=1$. Assume that $f_i=x_i-h_i$ where $h_i:\C^n\rightarrow \C$ is a homogeneous polynomial of degree $d$, where $d\in\N$. Then $f$ is invertible with polynomial inverse.
\end{conjecture}
\begin{theorem}\cite[Thm. II.2.1 and Cor. II.2.2]{Jacobian_overzicht_artikel}\label{lemma:Mathieu_Generalized_Conjecture_implies_Jacobian}
	Conjecture \ref{conj:d_restricted_Jacobian_conj} with $d=3$ implies the Jacobian Conjecture.
\end{theorem}
	
Combining Conjecture \ref{conj:d_restricted_Jacobian_conj}, Lemma \ref{lemma:Mathieu_Jacobian_enough_to_show_over_C} and Theorem \ref{lemma:Mathieu_Generalized_Conjecture_implies_Jacobian} gives the version of the Jacobian Conjecture we will be focus on in this section:
	
\begin{conjecture}\label{conj:Jacobian_conj_that_actually_matters}
	Let $n\in \N$, and $f:\C^n\rightarrow \C^n$ be a polynomial map with Jacobian $J(f)=1$. In addition, assume that $f_i=x_i-h_i$ where $h_i:\C^n\rightarrow \C$ is a homogeneous polynomial of degree $3$. Then $f$ is invertible with polynomial inverse.
\end{conjecture}
	
One of the key tools that will be required for the rest of this section is actually inverting the polynomial $f$. Let $f$ satisfy the hypotheses of Conjecture \ref{conj:Jacobian_conj_that_actually_matters}. Then $f(0)=0$ and $J(f)\neq 0$, so by the Holomorphic Inverse Function Theorem, there exists a local inverse $F$ around $0$ that can be written as a power series which converges around $0$. Thus, formally one can invert $f$ by the formal power series $F\in\C[[X_1,\ldots,X_n]]$.
	
Proving the existence and finding coefficients of implicitly defined formal power series is a well-established branch of combinatorics that dates back to the work of people like Jacobi \cite{Jacobi}. Nowadays, this principle is known as the multidimensional Lagrange inversion formula, and many equivalent versions of the inversion formula exist, see for example \cite{Abhyankar,Jacobian_overzicht_artikel,Good,Gessel_multi,Gessel}. 
	
In particular, we can find some properties of $f$ and the formal inverse $F$ that will be required for later proofs. To do that, we prove the following theorem, as given in \cite{Abhyankar,Jacobian_overzicht_artikel}. We follow the notation in \cite{Jacobian_overzicht_artikel}. Consider any $\mathbb{Q}$-algebra $A$, and let $a=(a_1,\ldots,a_n)\in A^n$. For $\alpha=(\alpha_1,\ldots,\alpha_n)\in \N_0^n$ we denote $$a^\alpha:=a_1^{\alpha_1}a_2^{\alpha_2}\cdots a_n^{\alpha_n},\qquad \alpha!:=\alpha_1!\alpha_2!\cdots \alpha_n!.$$ The Binomial Theorem can then be rewritten as $$\frac{1}{\alpha!}(a+b)^\alpha=\sum_{\beta+\gamma=\alpha}\frac{1}{\beta!\gamma!}a^\beta b^{\gamma}.$$
	
\begin{theorem}\cite[Thm. II.2.1]{Jacobian_overzicht_artikel}\label{thm:Jacobian_general_prerequisites_for_sum_equaling_1_and_h}
	Let $f=(f_1,\ldots,f_n):\C^n\rightarrow \C^n$ be a polynomial map such that $f_i=x_i-h_i$ with $h_i$ a polynomial map of degree $\geq 2$ and such that $f$ is invertible in the ring of formal power series, which is to say that there exists a formal inverse $F=(F_1,\ldots,F_n):\C^n\rightarrow \C^n$ of $f$ where $F_i$ is a formal power series. For any $U\in \C[[X_1,\ldots,X_n]]$ define the following:
	\begin{equation}
		\{U,f\}:=\sum_{(\alpha_1,\ldots,\alpha_n)\in \N_0^n } \frac{1}{\alpha!}\partial^{\alpha}(U(f)\cdot j(f)\cdot (X-f)^{\alpha}).
	\end{equation} 
	Then $\{U,f\}=U$. Here $\partial^{\alpha}=\partial_1^{\alpha_1}\cdots \partial_n^{\alpha_n}$ where $\partial_j=\frac{\partial}{\partial X_j}$, and $j(f)=\det(J(f))$ where $J(f)$ is the Jacobian of $f$.
\end{theorem}
\begin{proof}
	The proof consists of four steps, where in each step we reduce the problem to a simpler case.
	
	\noindent\textit{Step 1:} We first assume that $f=G\circ H$ where both $G$ and $H$ satisfy the theorem, which is to say that $\{U,G\}=U$ and $\{U,H\}=H$ for all $U\in \C[[X_1,\ldots,X_n]]$. Note that by the chain rule
	\begin{align*}
		j(G(H))=\det(J(G(H)))=\det(J(G)(H)\cdot J(H))=j(G)(H)\cdot j(H),
	\end{align*}
	where the notation $J(G)(H)$ means the Jacobian of $G$ evaluated at $H$. Using the Binomial Theorem, we find
	\begin{align*}
		\{U,&F\}=\sum_{\alpha\in \N_0^n } \frac{1}{\alpha!}\partial^{\alpha}(U(G(H))\cdot j(G(H))\cdot (X-H+H-G(H))^{\alpha})\\
		&=\sum_{\alpha\in \N_0^n } \partial^{\alpha}\left[U(G(H))\cdot j(G(H))\cdot \sum_{\substack{\beta,\gamma\in\N_0^n,\\
		\beta+\gamma=\alpha}}\frac{(X-H)^{\gamma}}{\gamma!}\frac{(H-G(H))^{\beta}}{\beta!}\right]\\
		&=\sum_{\beta\in \N_0^n } \partial^{\beta}\left[\sum_{\gamma\in\N_0^n}\partial^\gamma [U(G(H))\cdot j(G)(H)\cdot j(H)\cdot \frac{(H-G(H))^{\beta}}{\beta!}\frac{(X-H)^{\gamma}}{\gamma!}]\right]\\
		&=\sum_{\beta\in \N_0^n } \partial^{\beta}\,\left(\left\{U(G)j(G)\frac{(X-G)^{\beta}}{\beta!},H\right\}\right)\\
		&=\sum_{\beta\in \N_0^n } \partial^{\beta}\left(U(G)j(G)\frac{(X-G)^{\beta}}{\beta!}\right)\\
		&= \{U, G\} = U.
	\end{align*}
	So whenever $f$ is a composition of polynomials for which the theorem is true, it is true for $f$ as well.
		
	\noindent\textit{Step 2:} Define $G_n=f_n$ and $G_j$ such that $$G_j(X_1,\ldots,X_{j},f_{j+1},\ldots,f_{n})=f_j.$$ Since $f_j$ is invertible, $G_j$ is uniquely defined in $\C[[X_1,\ldots,X_n]].$ Then note that the maps $$H_j:=(X_1,\ldots,X_{j-1},G_j,X_{j+1},\ldots,X_n)$$ satisfy $$H_i\circ H_{i+1}\circ\cdots\circ H_n = (X_1,\ldots,X_{i-1},f_i,\cdots,f_n).$$ So $H_1\circ\cdots\circ H_n=f$. Using Step 1, it is enough to prove the theorem for $H_i$.
		
	\noindent\textit{Step 3:} By Step 2, it is enough to prove the theorem for $f=(f_1,X_2,\cdots,X_n)$. Then immediately we get $j(f)=\partial_1 f_1$ and $$(X-f)^{\alpha}=(X_1-f_1)^{\alpha_1}(X_2-X_2)^{\alpha_2}\cdots (X_n-X_n)^{\alpha_n}.$$ In other words, if $\alpha_j>0$ for $j=2,\ldots,n$, then $(X-f)^\alpha=0$. This reduces $\{U,f\}$ to
	\begin{align*}
		\{U,f\}=\sum_{\alpha_1=0}^\infty\frac{1}{\alpha_1!} \partial_1^{\alpha_1}\left[U(f) (\partial_1 f_1)(X-f_1)^{\alpha_1}\right].
	\end{align*}
	Now if we write $U=\sum_{m=0}^\infty U_m(X_2,\ldots,X_n)X_1^m$ with each $U_m\in \C[[X_2,\ldots,X_n]]$, we get that $$\{U,f\}=\sum_{\alpha_1=0}^\infty \frac{1}{\alpha_1!}\partial_1^{\alpha_1}\left[\sum_{m=0}^\infty U_m f_1^m (\partial_1 f_1) (X-f_1)^{\alpha_1}\right].$$ Swapping the two sums, and noting that $U_m$ is not dependent on $X_1$, we get 
	\begin{align*}
		\{U,f\}=\sum_{m=0}^\infty U_m \left[\sum_{\alpha_1=0}^\infty \frac{1}{\alpha_1!}\partial_1^{\alpha_1}\left(f_1^m (\partial_1 f_1)(X-f_1)^{\alpha_1}\right)\right] = \sum_{m=0}^\infty U_m \{X_1^m,f_1\}.
	\end{align*}
	So we can focus on the case $\{X_1^m,f_1\}$.
		
	\noindent\textit{Step 4:} We compute $\{X_1^m,f_1\}$. For simplicity, let us write $f$ instead of $f_1$. Now, we use the final assumption, that is, $f=X_1-h$, where $h$ is a homogeneous polynomial of degree $\geq 2$. This gives $\partial_1 f= 1- \partial_1 h$, and thus
	\begin{align}
		\{X_1^m,f\}&=\sum_{p=0}^\infty \frac{1}{p!}\partial_1^p ((X_1-h)^m(1-\partial_1 h)h^{p})\nonumber\\
		&=\sum_{p=0}^\infty \frac{1}{p!}\partial_1^{p} \left(\sum_{j=0}^m\binom{m}{j}(-1)^jX_1^{m-j}(1-\partial_1 h)h^{j+p}\right),\label{eq:Mathieu_equations_for_h_and_1}
	\end{align}
	where we used the Binomial Theorem again. Note that $$\frac{1}{p!}\partial^{p}(FG)=\sum_{\substack{r,s\in\N_0\\
			r+s=p}}\frac{1}{r!s!}\partial^{r}F\cdot \partial^{s} G,$$ so applying that to Equation (\ref{eq:Mathieu_equations_for_h_and_1}) gives
	\begin{align}\label{eq:{U,F}_bijna_klaar}
		\{X_1^m,f\}&=\sum_{p=0}^\infty \sum_{j=0}^m\sum_{r=0}^{\min(p,m-j)}(-1)^j\binom{m}{j}\binom{m-j}{r}\frac{X_1^{m-j-r}}{(p-r)!}\partial_1^{p-r}[(1-\partial_1 h)h^{j+p}].
	\end{align}
	We relabel the indices over which we sum. Define $t=r+j$. We claim that 
	\begin{align}\label{eq:Mathieu_enorme_som_swap_sums}
		\{X_1^m,f\}=\sum_{s=0}^\infty \sum_{t=0}^m\sum_{j=0}^t (-1)^j\binom{m}{j}\binom{m-j}{t-j}\frac{X_1^{m-t}}{s!}\partial_1^{s}[(1-\partial_1 h)h^{t+s}].
	\end{align}
	Note that substituting $t=r+j$ and $s=p-r$ in Equation (\ref{eq:{U,F}_bijna_klaar}) give the same summand as in Equation (\ref{eq:Mathieu_enorme_som_swap_sums}). So we can solely focus on the index ranges. Note that if we choose $s$ as our index, we can choose $r\in\{0,1,\ldots, m-j\}$ for $\partial_1^s X^{m-j}=0$ for any $r>m-j$. This also fixes $p=s+r$. Abusing notation by removing the summand for a second, we see that Equation (\ref{eq:{U,F}_bijna_klaar}) becomes $$\sum_{p=0}^\infty\sum_{j=0}^m\sum_{r=0}^{\min(p,m-j)}=\sum_{s=0}^\infty \sum_{j=0}^m\sum_{r=0}^{m-j}.$$ Next, note that $$\sum_{j=0}^m\sum_{r=0}^{m-j}= \sum_{j=0}^m\sum_{t=j}^m = \sum_{t=0}^m\sum_{j=0}^t$$ where we used in the first equation the equality $t=j+r$, and in the second we swapped the sums. Since the summands were already equal, we have Equation (\ref{eq:Mathieu_enorme_som_swap_sums}).
	
	Now that Equation (\ref{eq:Mathieu_enorme_som_swap_sums}) has been derived, we continue from there. Note that $\binom{m}{j}\binom{m-j}{t-j}=\binom{m}{t}\binom{t}{j}$. We thus get
	\begin{align*}
		\{X_1^m,f\}=\sum_{s=0}^\infty \sum_{t=0}^m\left[\sum_{j=0}^t (-1)^j\binom{t}{j}\right]\binom{m}{t}\frac{X_1^{m-t}}{s!}\partial_1^{s}[(1-\partial_1 h)h^{t+s}].
	\end{align*}
	But note that $\sum_{j=0}^t\binom{t}{j}(-1)^j=(1-1)^t=\delta_{t,0}$, so filling that in gives
	\begin{align*}
		\{X_1^m,f\}&=\sum_{s=0}^\infty \frac{X^{m}}{s!}\partial_1^{s}[(1-\partial_1 h)h^{s}]\\
		&=\sum_{s=0}^\infty \frac{X_1^{m}}{s!}\partial_1^s\bigg[h^{s}-\frac{1}{s+1}\partial_1(h^{s+1})\bigg]\\
		&=\sum_{s=0}^\infty X_1^{m}\bigg[\frac{1}{s!}\partial_1^{s}(h^{s})-\frac{1}{(s+1)!}\partial_1^{s+1} (h^{s+1})\bigg].
	\end{align*}
	The last sum is a telescoping sum, which just gives $X_1^m$. In short we have $\{X_1^m,f\}=X_1^m$. This proves the theorem.
\end{proof}
\begin{corollary} \cite[Cor. III.2.2 + III.2.3]{Jacobian_overzicht_artikel}\label{cor:Mathieu_inverse_of_polynomials_satisfy_some_relations}
	Let $n\in\N$ and let $f=(f_1,\ldots,f_n):\C^n\rightarrow \C^n$ be as in Conjecture \ref{conj:d_restricted_Jacobian_conj}. Denote $F=(F_1,\ldots,F_n)$ to be the formal inverse of $f$, and let $f_i=x_i-h_i$. Then, the following equalities are true
	\begin{align*}
		\sum_{(\alpha_1,\ldots,\alpha_n)\in\N_0^n}\frac{1}{\alpha!}\partial^{\alpha}h^\alpha = 1
	\end{align*} 
	and 
	\begin{align*}
		\sum_{(\alpha_1,\ldots,\alpha_n)\in\N_0^n}\frac{1}{\alpha!}\partial^{\alpha}(h^\alpha\cdot x_i) = F_i.
	\end{align*} 
\end{corollary}
\begin{proof}
	Apply Theorem \ref{thm:Jacobian_general_prerequisites_for_sum_equaling_1_and_h} with $U=1$ and $U=F$ to get the respective equalities.
\end{proof}

\subsection{The Mathieu Conjecture and the Jacobian Conjecture}
\subsubsection{Mathieu's proof}\label{sec:Mahtieu's_proof}
	
With some knowledge of the Jacobian Conjecture, we turn to Mathieu's main result: the Mathieu Conjecture implies the Jacobian Conjecture. The original proof of Mathieu relied heavily on algebraic (group) arguments, which sometimes could be circumvented. We note that this section is meant as a expository version of Mathieu's proof. Some proofs differ from Mathieu's, some details have been added where necessary, and the structure has been reorganized. For Mathieu's original paper, we refer to \cite{Mathieu}. The rest of the section is organized in the following way: in this section we will first prove the Jacobian Conjecture of degree $d$, given in Conjecture \ref{conj:Jacobian_conj_that_actually_matters}, assuming Conjecture \ref{conj:Mathieu_implied_conjecture}. In Section \ref{sec:simplification_Mathieu}, we will show Mathieu's Conjecture implies Conjecture \ref{conj:Mathieu_implied_conjecture}.
	
We recall that, as mentioned in Section \ref{sec:complexification}, the complexification $K_\C$ is an algebraic group. This opens a new window, which we will make use of sometimes. To do that, we need a few definitions. We follow here \cite{Goodman,Procesi}.
	
\begin{definition}
	Let $G=K_\C$ where $K\subseteq M(n,\C)$ is a connected compact Lie group. A function $f:G\rightarrow \C$ is called \emph{algebraic} if $f$ is the restriction to $G$ of a polynomial in $\C[x_{11},x_{12},\ldots,x_{nn},\det(x)^{-1}]$ where $x_{ij}$ are the matrix entry functions on $M(n,\C)$, and $\det(x)$ is the determinant function. Let us denote the set of algebraic functions by $\C[G]$.
\end{definition}
\begin{definition}
	Let $G=K_\C$ where $K$ is a connected compact Lie group. We say a representation $(\pi,V_\pi)$ of $G$ is \emph{algebraic} if $\dim(V_\pi)<\infty$ and if the matrix coefficients $m_{v,w}^\pi$ are algebraic for all $v,w\in V_\pi$.
	
	A representation $(\pi,V_\pi)$ of $G$ is called \emph{rational} if $V_\pi$ is the union of finite-dimensional subrepresentations which are algebraic.
\end{definition}

\begin{definition}
	Let $A$ be a commutative algebra over $\C$. We call $A$ a \emph{$G$-algebra} if $A$ admits a rational representation $(\pi,A)$ of $G$ such that $\pi:G\rightarrow \mathrm{Aut}(A)$ is a homomorphism.
\end{definition}
\begin{example}\label{ex:G-algebra}
	Let $K=SU(N)$ and $G=SL(N,\C)$. Then the representations $(\pi_{k\omega_1},S^k\C^N)$ and $(\pi_{k\omega_{N-1}},S^{k}(\C^N)^*)$ are algebraic representations of $SL(N,\C)$ for all $k\geq 1$. For future convenience, let $d\in\N$. Then it follows that $$(\pi_{dk\omega_1}\otimes \pi_{dk\omega_{N-1}},S^{dk}\C^N\otimes S^k(\C^N)^*)$$ is an algebraic representation of $SL(N,\C)$ for all $k$, and that also the algebra $$A:=\bigoplus_{k\in\N}S^{dk}\C^N\otimes S^k(\C^N)^*$$ is a rational representation of $SL(N,\C)$. In particular, $A$ is a $SL(N,\C)$-algebra.
\end{example}
\begin{theorem}\cite[Prop. 1 in 7.3.1]{Procesi} \label{thm:reductive_means_completely_reducible}
	Let $G=K_\C$, where $K$ is a compact connected Lie group. Then\, \emph{any} rational representation of $G$ is completely reducible.
\end{theorem}
\begin{theorem}\cite[Prop. in 8.7.1 + Thm. 3 in 8.7.2]{Procesi}\label{thm:C[G]_are_just_finite-type_functions}
	Let $G=K_\C$, where $K$ is a compact connected Lie group. Then the restriction map onto $K$, as a map going from $\C[G]$ to the continuous functions on $K$, is an isomorphism on the space of finite-type functions. In addition, $K$ lies Zariski dense in $G$, which means that any algebraic function can be found by analytically extending a finite-type function.
\end{theorem}

As we already hinted towards in Theorem \ref{thm:big_theorem_Mathieu_implies_Jacobian}, we only need to discuss the situation where $K=SU(N)$ and $G=SL(N,\C)$. So for the rest of the paper, let $K=SU(N)$ and $G=SL(N,\C)$.

Let $(\pi_\lambda,V(\lambda))$ be a finite dimensional irreducible representation of $K$ with highest weight $\lambda$. By Corollary \ref{cor:highest_weight_decomposition_complexification}, we can extend this uniquely to a finite dimensional irreducible representation of $G$. Then using Theorem \ref{thm:C[G]_are_just_finite-type_functions} shows that this representation is algebraic. In other words, any algebraic representation can be decomposed into irreducible representations, each labeled by the highest weight.
	
Let us move towards the result of Mathieu. Taking Theorem \ref{thm:reductive_means_completely_reducible} into account, we begin with some notation.
	
\begin{definition}
	If $A$ is a $G$-algebra, we have $A=\oplus_{\pi\in \widehat{G}} A_\pi$ where $A_\pi$ is the $\pi$-isotypical component of $A$. Let $f\in A$. We denote by $(f)_\pi$ the projection upon the $\pi$-isotypical component of $A$. Since the finite dimensional irreducible representation $\pi$ correspond to $\lambda\in P^{+}$, we also use the notation $(f)_\lambda$ instead of $(f)_\pi$.
\end{definition}

\begin{conjecture}\label{conj:Mathieu_implied_conjecture}
	Let $A=\bigoplus_{n\in\N}S^{dn}\C^N\otimes S^k(\C^N)^*$ be a graded $SL(N,\C)$-algebra. Let $\tau,\mu\in P^{+}$, and let $f\in A$. If $(f^n)_{n\tau}=0$ for all $n\in\N$, then $(f^n)_{n\tau+\mu}=0$ for $n$ large enough.
\end{conjecture}

\begin{theorem}\cite[Thm. 5.3 + 5.4]{Mathieu}\label{thm:Mathieu_implied_Mathieu_implies_Jacobian}
	Assume Conjecture \ref{conj:Mathieu_implied_conjecture}. Then Conjecture \ref{conj:d_restricted_Jacobian_conj} is true for all $d\in\N$. In particular, Conjecture \ref{conj:Jacobian_conj_that_actually_matters} is true (which implies the Jacobian Conjecture).
\end{theorem}

\begin{remark}
	 Note that we can set $d=3$ by Theorem \ref{lemma:Mathieu_Generalized_Conjecture_implies_Jacobian}, but will show it works for any $d\in\N$.
\end{remark}

\begin{proof}[Proof of Theorem \ref{thm:Mathieu_implied_Mathieu_implies_Jacobian}]
	Let $f=(f_1,\ldots,f_n):\C^n\rightarrow\C^n$ be a polynomial map with Jacobian $J(f)=1$, and let $f_j(x_1,\ldots,x_n)=x_j-h_j(x_1,\ldots,x_n)$ where $h_j:\C^n\rightarrow\C$ is a homogeneous polynomial of degree $d$. In addition, let $F:=(F_1,\ldots,F_n)$ be the formal inverse of $f$. As we discussed in Section \ref{sec:facts_Jacobian_Conjecture}, we have that $$F_i=\sum_{k=0}^\infty F_i^{(k)}$$ is a formal power series with each $F_i^{(k)}$ being a homogeneous polynomial of degree $k$. If we can show that there exists an $N\in\N$ such that for all $k\geq N$, we have $F_i^{(k)}=0$ for all $i$, then the inverse is a polynomial, and we are done. This will be the goal of the proof. By Corollary \ref{cor:Mathieu_inverse_of_polynomials_satisfy_some_relations}, $$F_i=\sum_{\alpha\in\N_0^n} \frac{1}{\alpha!}\partial^{\alpha}(h^\alpha \cdot x_i),$$ where we recall $\partial^{\alpha}=\prod_{i=1}^n \frac{\partial^{\alpha_i}}{\partial x_i^{\alpha_i}}$, $\alpha!=\alpha_1!\cdots\alpha_n!$, and $h^\alpha = h_1^{\alpha_1} \cdots h_n^{\alpha_n}$. Counting the degree of each term, one immediately gets that 
	\begin{align}\label{eq:proof_implied_Mathieu_implies_Jacobian_eq_1}
			F_i^{(k(d-1)+1)}=\sum_{\alpha\in \N_0^n\text{ such that }\sum_{j=1}^n\alpha_j=k}\frac{1}{\alpha!}\partial^{\alpha}(h^\alpha\cdot x_i),
	\end{align}
	and $F_i^{(k)}=0$ for all other $k\in\N$. Thus, we only need to show that $F_i^{(k(d-1)+1)}=0$ for sufficiently large $k$ to prove the theorem. 
		
	To analyze Equation (\ref{eq:proof_implied_Mathieu_implies_Jacobian_eq_1}) more, we focus on $\frac{1}{\alpha!}\partial^{\alpha}(h^\alpha\cdot x_i)$ specifically. To connect back to Section \ref{sec:prerequisites_Mathieu}, we define the element $$Q:=\sum_{i=1}^nh_i\otimes \partial_i \in S^d\C^n\otimes (\C^n)^*.$$ It is then clear that $Q^k \in S^{dk}\C^n\otimes S^k(\C^n)^*$ for all $k\in\N$. This means that we are working in the graded algebra $\bigoplus_{k\in\N}S^{dk}\C^n\otimes S^k(\C^n)^*$. Let $k$ be fixed for now, i.e., we work in $S^{dk}\C^n\otimes S^k(\C^n)^*$ to be able to use the tools developed in Section \ref{sec:prerequisites_Mathieu}. Note that 
	\begin{align*}
		\Psi\left(Q^k\right) &= k!\sum_{i=1}^n\,\sum_{\substack{\alpha\in\N_0^n,\\ \sum_{j=1}^n\alpha_j=k}} \frac{1}{\alpha!}\partial^{\alpha}(h^{\alpha}x_i)\otimes \partial_i\\
		&= k!\sum_{i=1}^n F_i^{(k(d-1)+1)}\otimes \partial_i,
	\end{align*} where $\Psi$ is as in Definition \ref{def:Mathieu_Psi}. We argued it was enough to show $F_i^{(k(d-1)+1)}=0$ for sufficiently large $k$, so showing that $\Psi(Q^k)=0$ for sufficiently large $k$ is enough to prove the theorem. This will be the goal for the rest of the proof. 
	
	We shift gears for the moment and consider another $G$-equivariant map to see if we can gather more information on $Q^k$. Recall from Corollary \ref{cor:Mathieu_div_is_an_isomorphism_onto_subspace} that the map $\mathrm{div}:S^{k-l+1}\C^n\otimes (\C^n)^*\rightarrow  S^{k-l} {\C}^n$ is a projection upon the $V((k-l)\omega_1)\simeq S^{k-l}\C^n$ component of $S^{k-l+1}\C^n\otimes (\C^n)^*$. Applying $\mathrm{div}$ to $Q^k$ gives 
	\begin{align}\label{eq:Mathieu_evaluatie_div(Q^k)}
		\mathrm{div}(Q^k) = k!\sum_{\substack{\alpha \in \N_0^n,\\ \sum_j\alpha_j=k}}\frac{1}{\alpha!}\partial^{\alpha}h^\alpha.
	\end{align}
	We claim that $\mathrm{div}(Q^k)=0$. We know from Corollary \ref{cor:Mathieu_inverse_of_polynomials_satisfy_some_relations} that 
	\begin{align}\label{eq:div_is_zero}
		\sum_{\alpha\in\N_0^n}\frac{1}{\alpha!}\partial^{\alpha}h^\alpha=1,
	\end{align} where $1$ is constant function. The degree of $\partial^{\alpha}h^{\alpha}$ is $(d-1)\sum_j\alpha_j$. Since $d\geq 2$ and the function $1$ has degree 0, Equation \ref{eq:div_is_zero} can only be true if $$\sum_{\sum_{j}\alpha_j=k}\frac{1}{\alpha!}\partial^{\alpha}h^\alpha=0$$ for all $k\geq1$. Hence $$\mathrm{div}(Q^k)=0\qquad \forall k\geq 1.$$ By Corollary \ref{cor:Mathieu_div_is_an_isomorphism_onto_subspace}, the irreducible representation $(\pi_{(k-l)\omega_1},V((k-l)\omega_1))$ occurs in $S^k\C^n\otimes S^l(\C^n)^*$ with multiplicity 1. In other words we have $$(Q^k)_{k(d-1)\omega_1}=0$$ for all $k\geq 1$.
	
	Now we are in the position to use Conjecture \ref{conj:Mathieu_implied_conjecture}, for we recall that we are working in the commutative graded $SL(n,\C)$-algebra $A=\bigoplus_{k=0}^\infty S^{dk}\C^n\otimes S^k(\C^n)^*$. The conjecture gives us $$(Q^k)_{k(d-1)\omega_1+\mu}=0$$ for sufficiently large $k$, where we choose $\mu$ to be $\mu=\omega_1+\omega_{n-1}$. So $$(Q^k)_{(k(d-1)+1)\omega_1+\omega_{n-1}}=0$$ for all $k$ large enough. But by Lemma \ref{lemma:Mathieu_Psi_isomorphism}, $\Psi$ is the zero map on all other isotypical components of $S^{dk}\C^n\otimes S^k(\C^n)^*$. In other words, $\Psi(Q^k)=0$ for all $k$ large enough, hence proving the theorem.
\end{proof}
\begin{remark}
	We would like to remark that Conjecture \ref{conj:Mathieu_implied_conjecture} is slightly altered with respect to the paper of Mathieu. We only allow $\mu\in P^{+}$, while Mathieu allows any integral weight $\mu\in\mathfrak{t}_\C^*$. This stricter claim is chosen because if one allows any integral weight $\mu$, some parts of the lemma are not well-defined. For example, we can choose $\tau$ to lie on the boundary of the positive Weyl chamber and choose $\mu$ in such a way that $n\tau+\mu$ is not dominant for all $n\in\N$. An example for this would be the Lie algebra $\mathfrak{sl}(3,\C)$. We choose $\tau=\omega_1$ and $\mu=\alpha_1$. This way, $(f^n)_{n\tau+\mu}$ is not defined. This correction does not change the overall paper, for it only uses $\mu=\omega_1+\omega_{n-1}$, which is dominant anyway.
\end{remark}

\subsubsection{Simplifying Conjecture \ref{conj:Mathieu_implied_conjecture}}\label{sec:simplification_Mathieu}
The rest of this section is dedicated showing that the Mathieu Conjecture implies Conjecture \ref{conj:Mathieu_implied_conjecture}. We first recall the Mathieu Conjecture for any connected compact Lie group $K$. As was discussed above, it is enough to consider $G=SL(N,\C)$ and $K=SU(N)$.
	
\begin{conjecture}[The Mathieu Conjecture] \cite[Main Conjecture]{Mathieu}\label{conj:Mathieu}
	Let $K$ be a compact connected Lie group, and let $f$ and $h$ be complex-valued finite-type functions. Assume $$\int_K f^n(k)dk=0$$ for all $n\in \N$. Then $$\int_K f^n(k)h(k)dk=0$$ for all $n$ large enough. Here $dk$ is the Haar measure on $K$.
\end{conjecture}

\begin{lemma}
	Let $G=K_\C$ be the complexification of $K$, where $K$ is as in Conjecture \ref{conj:Mathieu}. Define $L:\C[G]\rightarrow \C[G]$, $f\mapsto f_{\triv}$. Then $$L(f)=\int_Kf(k)dk.$$ Here $dk$ is the normalized Haar measure on $K$. 
\end{lemma}
\begin{proof}
	Let $f\in\C[G]$. By Theorem \ref{thm:C[G]_are_just_finite-type_functions}, restricting to $K$ makes $f$ a finite-type function on $K$, i.e., a sum of matrix coefficients. Using Lemma \ref{lemma:Schur_orthogonality} together with the fact that the constant matrix coefficient $m_{1,1}^{\triv}(k):=1$ for all $k\in K$, we get the result.
\end{proof}
	
\begin{lemma}\label{lemma:Mathieu_1.3}
	Let $(\pi,V_\pi)\in\widehat{G}$ and $f\in \C[G]$ such that $(f^n)_\triv=0$ for all $n\in\N$. Assume the Mathieu Conjecture is true (see Conjecture \ref{conj:Mathieu}). Then $(f^n)_\pi=0$ for all $n$ large enough. 
\end{lemma}
\begin{proof}
	Let $\pi\in\widehat{G}$. We assume Conjecture \ref{conj:Mathieu}, which gives that for any finite type function $h$ the integral $$\int_K f^n(k)h(k)dk=0$$ for all $n\geq N(h)$ with $N(h)\in\N$. Consider the dual of $\pi$, i.e., $(\pi^*,V_{\pi^*})$ and take $h=m_{u,v}^{\pi^*}$ where $u,v\in V_{\pi^*}$ are arbitrary. Note that $m_{u,v}^{\pi^*}(g)=\overline{m_{v,u}^{\pi}(g)}$. So $$0=\int_K f^n(k)h(k)dk=\int_{K}f^n(k)m_{u,v}^{\pi^*}(k)dk=\int_{K}f^n(k)\overline{m_{v,u}^\pi(k)}dk$$ for all $n\geq N(h)$. Since $\pi$ is a finite dimensional representation, the space of matrix coefficients of $\pi$ is finite dimensional, so we can choose a basis, say $\{m_1,\ldots,m_p\}$. Apply the above arguments for each $m_j$, which gives a natural number $N(m_j)$ for which $\int f^n m_j \,dk=0$ for all $n\geq N(m_j)$. Thus, we get a finite collection of integers $N(h)$. So there is a maximum, say $N$. Hence, by construction, we have $$\int_K f^n(k)h(k)dk=0$$ for all $h\in \C[G]_{\pi^*}$ for all $n\geq N$, and thus by Schur orthogonality relations we must have $(f^n)_\pi=0$ for $n\geq N$.
\end{proof}
	
For the remainder of this section, we require some knowledge of commutative algebra and filtered colimits. We do not assume the reader is too familiar with this, so we included a brief discussion on the necessary material in Section \ref{sec:cheat_sheet}.
	
\begin{lemma}\label{lemma:Mathieu_more_general_conjecture_G_algebra_implied_by_Mathieu_conj}
	Let $K=SU(N)$ and $ G=SL(N,\C)$. Consider $\tau\in P^{+}$ and let $A=\bigoplus_{k\in\N}S^k\C^N\otimes S^k(\C^N)^k\otimes V(k\tau^*)$. Consider $\pi\in\widehat{G}$ and $f\in A$ such that $(f^n)_\triv=0$ for all $n\in\N$. Assume Mathieu's Conjecture to be true for $K$, then $(f^n)_\pi=0$ for all $n$ large enough. 
\end{lemma}
\begin{proof}
	Assume the Mathieu Conjecture is true for $K$, and let $A=\bigoplus_{k\in\N}S^k\C^N\otimes S^k(\C^N)^k\otimes V(k\tau^*)$. The idea is to reduce the problem to the case $A=k[G]$ for some field $k$, embed $k$ into $\C$ and then apply Lemma \ref{lemma:Mathieu_1.3}. We will do so in three steps. 
	
	\textit{Step 1:} We remark that $A$ is reduced. In addition, we can reduce the lemma to a noetherian subring of $A$. To show this, note that $A$ can be written as the filtered colimit of noetherian algebras. To be more specific, consider $I$ to be the set of finite unions of irreducible subrings of $A$; see Section \ref{sec:filtered_colim}. Then $(I,\subseteq)$ is a filtered poset, and thus $\mathrm{colim}_{i\in I}\C[i]=A$. In particular, $\C[i]$ is noetherian, reduced and $G$-invariant. 
		
	By Theorem \ref{thm:reductive_means_completely_reducible}, both $A$ and $\C[i]$ are completely reducible as representations of $G$. So they decompose as a direct sum of irreducible representations, i.e. $$A=\bigoplus_{\tau\in\widehat{G}}A_\tau,\qquad\qquad \C[i]=\bigoplus_{\rho\in\widehat{G}}(\C[i])_\rho.$$ 
	Note that these decompositions are compatible, as in $(\C[i])_\tau\subseteq A_\tau$ by Schur's Lemma. These direct sums are in the category of $\C$ vector spaces. Since the direct sum can be seen as a colimit, they commute with other filtered colimits. Using Schur's Lemma, we see that these colimits commute in the category of $G$-algebras as well, see Lemma \ref{lemma:categorical_nonsense_why_swapping_is_okay}, thus giving
	\begin{align*}
		(\mathrm{colim}_{i\in I} \C[i])_\tau &= \big(\mathrm{colim}_{i\in I} \bigoplus_{\rho\in\widehat{G}} (\C[i])_\rho\big)_\tau \\
		&\simeq \big(\bigoplus_{\rho\in \widehat{G}}[\mathrm{colim}_{i\in I}\C[i]_\rho]\big)_\tau\\
		&= \mathrm{colim}_{i\in I}(\C[i]_\tau).
	\end{align*}
	This is true for all $\tau\in \widehat{G}$, so if we can prove the lemma for the noetherian algebra $\C[i]$, $i\in I$, we have it for $A=\mathrm{colim}_{i\in I}\C[i]$. In other words, in the lemma we can replace $A$ by $\C[i\in I]$, where $I$ is a finite union of irreducible subrings, i.e., a noetherian algebra.
	
	\textit{Step 2:} By Theorem \ref{thm:reduced_ring_embedded_in_product_of_fields}, any reduced noetherian ring $R$ can be embedded as a subring of $\Phi:=\prod_{j=1}^m k_j$ where $k_j$ is a field. Let us call this map $\xi:R\rightarrow \Phi$.
	
	Let us concentrate on the action of $G$ on $A$. The action of $G$, denoted by $\pi:G\rightarrow \mathrm{Aut}(A)$, induces a map $\Delta:G\times A\rightarrow A$ by $\Delta(g,a)=\pi(g)a$. Then the pullback is a map $$\Delta^*:\C[A]\rightarrow \C[G\times A]\simeq \C[G]\otimes_\C \C[A],$$ where $\Delta^*(f)=f\circ \Delta$. Note that we can equip the spaces $\C[A]$ and $\C[G\times A]$ with two actions of $G$ as well. Namely let $G$ act on $\C[A]$ by $\rho(g)f(a)=f(\pi(g^{-1})a)$, and let $G$ act on $\C[G\times A]$ by $\sigma(g)f(h,a)=f(g^{-1}h,a)$. With these actions, we claim that $\Delta^*$ is an injective $G$-equivariant map.
	
	We first prove that $\Delta^*$ is injective. Let $\Delta^*(f)=\Delta^*(h)$. Then $f\circ \Delta = h\circ \Delta$. Since $\Delta$ is surjective, for every $a\in A$ there exist $(g,b)\in G\times A$ such that $\Delta(g,b)=a$. This means $$f(a)=(f\circ \Delta)(g,b)=(h\circ \Delta)(g,b)=h(a).$$ This is true for all $a\in A$, so $f=h$. So $\Delta^*$ is injective. To show it is $G$-equivariant, we see
	\begin{align*}
		[(\sigma(g)\circ\Delta^*)(f)](h,a)&=[\sigma(g)(f\circ \Delta)](h,a)=(f\circ\Delta)(g^{-1}h,a)\\
		&=f(\pi(g^{-1}h)a)= f(\pi(g^{-1})\pi(h)a)\\
		&=(\rho(g)\circ f)(\pi(h),a) = (\rho(g)\circ f)(\Delta(h,a))\\
		&=\Delta^*(\rho(g)\circ f)(h,a) = [(\Delta^*\circ \rho(g))(f)](h,a).
	\end{align*}
	So $\sigma(g)\circ\Delta^*=\Delta^*\circ \rho(g)$. 
	Note that $\C[A]\simeq A$ and $\C[G\times A]\simeq \C[G]\otimes\C[A]$. Denoting $A_\triv$ as the algebra $A$ with the trivial $G$-action, these isomorphisms translate to that $\Delta^*$ is an injective $G$-equivariant map as
	\begin{align}\label{eq:Mathieu_injective_G-equivariant_map}
		\Delta^*:A\rightarrow \C[G]\otimes_\C A_\triv.
	\end{align}
	Together with the injective map $\xi:A\rightarrow \Phi$ we have an injective $G$-equivariant map $$A\xrightarrow{\Delta^*} \C[G]\otimes_\C A_\triv\xrightarrow{\mathrm{id}\otimes \xi} \C[G]\otimes_\C \Phi.$$
	
	Summarizing, we have an injective $G$-equivariant map $A\rightarrow \C[G]\otimes_\C \Phi$, which allows us to embed $A$ into $\C[G]\otimes_\C\Phi$ as a subalgebra. However, since $\Phi=\prod_{j=1}^m k_j$, we have 
	$$\C[G]\otimes_\C \Phi = \C[G]\otimes_\C \prod_{j=1}^m k_j \simeq \prod_{j=1}^m \C[G]\otimes_\C k_j \simeq \prod_{j=1}^m k_j[G]$$ because $\C[G]$ is finitely generated. So $A$ can be embedded in $\prod_{j=1}^m k_j[G]$ as an algebra. In fact, the lemma we have to prove is stated in such a way that if we have shown it for $\prod_{j=1}^m k_j[G]$, then it is automatically true for any subalgebra $A\subseteq \prod_{j=1}^m k_j[G]$. So we can just assume $A = \prod_{j=1}^m k_j[G]$. By construction, if $f=(f_1,f_2,\ldots,f_m)\in \prod_{j=1}^m k_j[G]$ and $\tau\in\widehat{G}$, then $$(f)_\tau = (f_1,f_2,\ldots, f_m)_\tau = ((f_1)_\tau,\ldots, (f_m)_\tau).$$ In other words, if we can prove the lemma for the $k_j[G]$ for some field $k_j$, then it can be shown for $A= \prod_{j=1}^m k_j[G]$. 
		
	\textit{Step 3:} So it is enough to show the lemma for $A=k[G]$ where $k$ is some field. Let $f\in k[G]$. Then $f$ is a polynomial with finitely many coefficients. Let $k'$ be the subfield generated by the coefficients of $f$. This subfield can be embedded in $\C$ by choosing a transcendental basis. In other words, we can assume $f\in \C[G]$. However, we are now in the position to apply Lemma \ref{lemma:Mathieu_1.3}, which proves the lemma.
\end{proof}
	
\begin{theorem}\label{thm:Mathieu_Mathieu_implies_Mathieu_inplied}
	Let $K=SU(N)$, $G=SL(N,\C)$ and let $A=\bigoplus_{k\in\N}S^k\C^N\otimes S^k(\C^N)^k$. Let $\tau,\mu\in P^{+}$, and let $f\in A$. Assume the Mathieu Conjecture to be true for $K$. Then Conjecture \ref{conj:Mathieu_implied_conjecture} is true, i.e. if $(f^k)_{k\tau}=0$ for all $k\in\N$, then $(f^k)_{k\tau+\mu}=0$ for $k$ large enough.
\end{theorem}	
	
\begin{proof}
	Let $f\in A$ such that $(f^k)_{k\tau}=0$ for all $k\in \N$. Note $\mu+k\tau$ is a dominant integral weight for all $k\in\N$. Let $\xi\in V(\tau^*)$ be a non-zero highest weight vector of $\tau^*$. Identifying $\xi^k$ as an element in $S^kV(\tau^*)$, i.e., $\xi^k\in S^kV(\tau^*)$, we see that $\xi^k$ is a highest weight vector with weight $k\tau^*$, so $\xi^k\in V(k\tau^*)$. 	
	
	Now let us consider $L:= f\otimes \xi$. Then $L^k = f^k\otimes \xi^k\in A\otimes V(k\tau^*)$. By Lemma \ref{lemma:Mathieu_3.1}\textrm{.(iii)}, we get $$(L^k)_\triv = \bigoplus_{\pi\in\widehat{G}}(f^k)_\pi\otimes (\xi^k)_{\pi^*}= (f^k)_{k\tau}\otimes \xi^k$$ where we used that $\bigoplus_{\pi\in\widehat{G}}(\xi^k)_{\pi^*}=(\xi^k)_{k\tau^*}=\xi^k$. Since we assumed that $(f^k)_{k\tau}=0$ for all $k\in\N$, we find that $(L^k)_\triv=0$ for all $k\in\N$. 
	
	We note that $L\in \sum_{j=0}^\infty A_j\otimes V(j\tau^*)$, where $A_j=S^j\C^N\otimes S^j(\C^N)^*$, so we are in the position to use Lemma \ref{lemma:Mathieu_more_general_conjecture_G_algebra_implied_by_Mathieu_conj}, which means that $(L^k)_{\mu}=0$ for sufficiently large $k$. In particular, that means
	\begin{align}\label{eq:Mathieu_Mathieu_implies_Conj}
		0=(L^k)_{\mu}=(f^k\otimes \xi^k)_{\mu} &= \left[\bigoplus_{\nu(\mu)\in P^{+}} (f^k)_{\nu(\mu)}\right]\otimes (\xi^k)_{k\tau^*}
	\end{align}
	for some unknown $\nu(\mu)\in P^{+}$. Before we continue, we note that there are finitely many $\nu\in P^{+}$ such that $V(\mu)$ occurs in $V(\nu)\otimes V(k\tau^*)$. To prove this, we interpret the representations in question, such as $V(\mu)$, as representations of $K$. This allows us to consider $\langle \chi_{\mu},\chi_{\nu\otimes k\tau^*}\rangle$. But using the same trick as before, we see 
	\begin{align*}
		\langle \chi_{\mu},\chi_{\nu\otimes k\tau^*}\rangle &=\int_{K}\chi_{\mu}(g)\overline{\chi_{\nu}(g)\chi_{k\tau^*}(g)}\,dg\\
		&=\int_{K}\chi_{\mu}(g)\chi_{k\tau}(g)\overline{\chi_{\nu}(g)}\,dg=\overline{\langle \chi_{\nu},\chi_{\mu\otimes k\tau}\rangle}.
	\end{align*}
	Since $\mu$ and $k\tau$ are fixed and both are finite dimensional, the inner product is finite, and thus there are only finitely many $\nu\in P^{+}$ that satisfy the above. In particular note that by Lemma \ref{lemma:Mathieu_3.1} $\nu=\mu+k\tau$ is one of the weights satisfying the above.
	
	In short, there are only finitely many tensor products that have non-zero $\mu$-component in Equation (\ref{eq:Mathieu_Mathieu_implies_Conj}). This means that we can swap the sum and the tensor product, giving
	\begin{align*}
		0=(L^k)_\mu&= \left[\bigoplus_{\lambda\in P^+} (f^k)_{\lambda}\right]\otimes (\xi^k)_{k\tau^*}=\bigoplus_{\lambda\in P^+} \left[(f^k)_{\lambda}\otimes (\xi^k)_{k\tau^*}\right].
	\end{align*}
	But this must mean that $(f^k)_{\lambda}\otimes (\xi^k)_{k\tau^*} = 0$, since the sum is direct. We know that $\lambda=\mu+k\tau$ is one of the summands in the direct sum, so $(f^k)_{\mu+k\tau}\otimes (\xi^k)_{k\tau^*}=0$. Since $\xi^k\neq 0$ for all $k\in\N$, we must find that $(f^k)_{\mu+k\tau}=0$, which was what we had to show.
\end{proof}
Finally, we get to Theorem \ref{thm:big_theorem_Mathieu_implies_Jacobian}.
\begin{corollary}
	Let $N\in\N$ and consider $K=SU(N)$. Assume the Mathieu Conjecture is true for $SU(N)$. Then the Jacobian Conjecture for $\C^N$ is true.
\end{corollary}
\begin{proof}
	Combine Theorem \ref{thm:Mathieu_Mathieu_implies_Mathieu_inplied} and Theorem \ref{thm:Mathieu_implied_Mathieu_implies_Jacobian}.
\end{proof}

\newpage
\begin{appendices}
\section{Some results on commutative algebra and category theory}\label{sec:cheat_sheet}
This appendix is written for readers not familiar with some of the techniques that are used in this section. This section will mostly be a collection of results and definitions. It is based on Atiyah and MacDonald \cite{Atiyah}, Matsumura \cite{Matsumura}, and the Stack Project \cite[\href{https://stacks.math.columbia.edu/tag/00AO}{Tag 00AO}]{stacks-project}, and we recommend these sources for more information. In this section, all rings and algebras are considered to be unital and commutative.
	
\subsection{Some elementary commutative algebra results}
We start with some definitions and early results.
\begin{definition}
	Let $R$ be any ring. We call $R$ an \emph{integral domain} if there are no zero-divisors, i.e., if $xy=0$ for any $x,y\in R$, then either $x=0$ or $y=0$. We call $R$ a \emph{field} if for every $x\in R\setminus\{0\}$ there exists an $x^{-1}\in R$ such that $xx^{-1}=1$. 
\end{definition}
\begin{definition}
	Let $R$ be any ring. We call $R$ \emph{noetherian} if every ideal in $R$ is finitely-generated, which is to say that for every ideal $I$ there exist $x_1,\ldots,x_n\in R$ such that $I=(x_1,\ldots,x_n)$.
\end{definition}
One of the notions we will need is the nilradical:
\begin{definition}
	Let $R$ be a ring, then we define the \emph{nilradical}, denoted $\mathfrak{nilrad}(R)$, by $$\mathfrak{nilrad}(R):=\{x\in R\mid x^n=0 \text{ for some } n\in\N\}.$$ 
\end{definition}
\begin{definition}
	Let $R$ be any ring. We call $R$ \emph{reduced} if $R$ contains no non-zero nilpotent elements.
\end{definition}
\begin{lemma}
	Let $R$ be any ring. Then $R$ is reduced if and only if $\mathfrak{nilrad}(R)=0$. In particular $R/\mathfrak{nilrad}(R)$ is reduced.
\end{lemma}
Another essential concept is that of prime ideals.
\begin{definition}
	A proper ideal $I$ in $R\neq 0$ is called a \emph{prime ideal} if $R/I$ is an integral domain, that is to say, if $xy\in I$ then $x\in I$ or $y\in I$. We will adopt the notation of $\mathfrak{p}$ for a prime ideal in $R$.
\end{definition}
\begin{lemma}\label{lemma:nilrad_is_intersection_prime_ideals}
	Let $R$ be a ring, then the nilradical is exactly $$\mathfrak{nilrad}(R)=\bigcap_{\mathfrak{p} \text{ is prime}}\mathfrak{p}$$
\end{lemma} 
\begin{proof}
	Let $x\in\mathfrak{nilrad}(R)$, then $x^n=0$ for some $n$. We note that $0\in \mathfrak{p}$ for any ideal. So if $\mathfrak{p}$ is a prime ideal, then $x\cdot x^{n-1}=0\in \mathfrak{p}$, so $x$ or $x^{n-1}\in \mathfrak{p}$. If $x\in\mathfrak{p}$ we are done, and if not, repeat the process on $x\cdot x^{n-2}=x^{n-1}\in \mathfrak{p}$. By induction, $x\in \mathfrak{p}$. The other inclusion requires Zorn's Lemma, for which we refer to \cite[Prop. 1.8]{Atiyah}.
\end{proof}
		
\begin{lemma}\label{lemma:nilrad_nilpotent}
	Let $R$ be a noetherian ring. Then $\mathfrak{nilrad}(R)$ is nilpotent, that is to say there exists $n\in \N$ such that $\mathfrak{nilrad}(R)^n=0$.
\end{lemma}
\begin{proof}
	Since $R$ is noetherian, there exist $x_1,\ldots,x_m\in R$ such that $\mathfrak{nilrad}(R)=(x_1,\ldots,x_m)$. But since $x_j\in\mathfrak{nilrad}(R)$, there exist $k_j$ such that $x_j^{k_j}=0$. Now observe that if we define $n:=\sum_{j=1}^mk_j$ then $\mathfrak{nilrad}(R)^n$ is generated by elements of the form $x_1^{l_1}\cdots x_m^{l_m}$ with $\sum_{j=1}^ml_j=n$. Thus, there is at least one $l_j\geq k_j$, meaning $x_1^{l_1}\cdots x_m^{l_m}=0$.
\end{proof}
		
\subsection{Embedding of rings and primary decomposition}
Next, we go to some of the more specific results, with this section dedicated to embeddings of (noetherian) rings into other rings.

\begin{lemma}\cite[\href{https://stacks.math.columbia.edu/tag/00FR}{Tag 00FR}]{stacks-project}
	Let $R$ be a noetherian ring. Then $R$ has finitely many minimal prime ideals $\mathfrak{p}_1,\ldots,\mathfrak{p}_n$, where minimal means with respect to inclusion.
\end{lemma}
\begin{corollary}
	Let $R$ be a noetherian reduced ring. Let $\mathfrak{p}_1,\ldots,\mathfrak{p}_n$ be the minimal prime ideals. Then the map $$\phi:R\rightarrow \prod_{i=1}^nR/\mathfrak{p}_i,\qquad r\mapsto (r+\mathfrak{p}_1,\ldots,r+\mathfrak{p}_n)$$ is injective, meaning $R$ can be embedded in a product of integral domains.
\end{corollary}
\begin{proof}
	Note $x\in \ker(\phi)$. Then $x\in \mathfrak{p}_i$ for all $i$. Thus, $x\in \bigcap_{i=1}^n\mathfrak{p}_i$, which is exactly $\mathfrak{nilrad}(R)$ by Lemma \ref{lemma:nilrad_is_intersection_prime_ideals}. Since $R$ is reduced, $x=0$. This shows injectivity. Since $\mathfrak{p}_i$ are prime ideals, $R/\mathfrak{p}_i$ is an integral domain.
\end{proof}

However, we might not always be in reduced rings. We still wish to have some embedding, which requires a primary decomposition.

\begin{definition}
	Let $R$ be a ring and $I$ be an ideal. We call $I$ \emph{primary} if $xy\in I$ gives $x\in I$ or $y^n\in I$ for some $n\in\N$. 
\end{definition}
\begin{lemma}
	Let $R$ be a ring. Then the ideal $I$ is primary if and only if $R/I\neq 0$ and has only nilpotent zero divisors.
\end{lemma}
Clearly, any prime ideal $\mathfrak{p}$ is primary. On the other hand, not every primary ideal is prime. However, similar to how one decomposes numbers into a product of prime numbers, we can decompose ideals into primary ideals. For that, we make the following definition
\begin{definition}
	Let $R$ be a ring, and let $I$ be an ideal. We say that $I$ has \emph{a primary decomposition} if there exist primary ideals $\mathfrak{q}_i$ with $i=1,\ldots,n$ such that $$I=\bigcap_{i=1}^n\mathfrak{q}_i.$$
\end{definition}
Of course, it is not necessary that such a decomposition exists, nor that it is unique. In fact, in most cases, there is no uniqueness result. In noetherian rings, however, we have an existence result.

\begin{theorem}\cite[Thm. 7.13]{Atiyah}\label{thm:Noetherian_ring_primary_decomp}
	Let $R$ be a noetherian ring, and let $I$ be an ideal. Then $I$ has a primary decomposition.
\end{theorem}

\begin{lemma}\label{lemma:Noetherian_rings_embedded_in_primary_algebras}
	Let $R$ be a noetherian ring. Then there exist finitely many primary ideals $\mathfrak{q}_i$ such that $$\phi:R\rightarrow \prod_{i=1}^nR/\mathfrak{q}_i$$ is injective. 
\end{lemma}
\begin{proof}
	Since $R$ is noetherian, by Theorem \ref{thm:Noetherian_ring_primary_decomp} the zero ideal $(0)$ has a primary decomposition, i.e. there exist primary ideals $\mathfrak{q}_i$, $i=1,\ldots, n$ such that $$0=\bigcap_{i=1}^n\mathfrak{q}_i.$$ Now consider the map $\phi:R\rightarrow \prod_{i=1}^nR/\mathfrak{q}_i$. Then $\ker(\phi)=\bigcap_{i=1}^n\mathfrak{q}_i=0$, so this map is injective.
\end{proof}

\subsection{Some results about quotients}
Next we introduce the theory of quotient rings.
\begin{definition}
	Let $R$ be an integral domain. We call a field $F$ \emph{the fraction field of $R$} if it is the smallest field in which $R$ can be embedded. The construction is as follows:
	
	Let $\sim$ be an equivalence relation on $R\times R\setminus\{(0,0)\}$, given by $(n,d)\sim (n',d')$ if and only if $nd'=n'd$. Any equivalence class can be denoted as $[(n,d)]=:\frac{n}{d}$. Then $F$ is given by $F=(R\times R\setminus\{(0,0)\})/\sim$ with the addition and multiplication $$\frac{a}{b}+\frac{c}{d}=\frac{ad+bc}{ad},\qquad \frac{a}{b}\frac{c}{d}=\frac{ac}{bd}.$$
\end{definition}
\begin{definition}
	The above can be done more generally. Let $R$ be some ring, and $S\subseteq R$ such that $1\in S$, and $S$ is closed under multiplication. Define the equivalence relation $\equiv$ on $R\times S$ by $$(a,s)\equiv (b,t):\Leftrightarrow \text{ there exists }u\in R\text{ such that }(at-bs)u=0.$$ This is an equivalence relation, and we denote $S^{-1}R$ by the set of equivalence classes. If we define the equivalence class as $\frac{a}{s}:=[(a,s)]$ then the set $S^{-1}R$ becomes a ring with the addition $+$ and multiplication $\cdot$ where $$\frac{a}{s}+\frac{b}{t}=\frac{at+sb}{st},\qquad \frac{a}{s}\frac{b}{t}=\frac{ab}{st}.$$ We call $S^{-1}R$ the \emph{ring of fractions} of $R$ with respect to $S$.
\end{definition}
\begin{lemma}\cite[Prop. 3.1]{Atiyah}
	The ring of fractions of $R$ with respect to $S$ has a universal property, as in that any ring homomorphism $f:R\rightarrow B$ such that $f(s)$ is a unit in $B$ for every $s\in S$, descends to a unique ring homomorphism $g:S^{-1}R\rightarrow B$. 
\end{lemma}

\begin{lemma}\cite[Prop. 7.3]{Atiyah}\label{lemma:Noetherian_ring_fraction_field}
	Let $R$ be a noetherian ring and $S$ any multiplicatively closed subset of $R$. Then $S^{-1}R$ is noetherian.
\end{lemma}
We will be using quotient rings multiple times in different circumstances, so we will develop the theory in more detail.
\begin{definition}
	Let $R$ be a ring. We call $R$ a \emph{local ring} if there exists exactly one maximal ideal $\mathfrak{m}$. The field $R/\mathfrak{m}$ is called the \emph{residue field}.
\end{definition}
\begin{definition}
	Let $R$ be a ring and let $\mathfrak{p}$ be some prime ideal in $R$. Define $S=R\setminus \mathfrak{p}$. We denote $S^{-1}R$ as $R_\mathfrak{p}$. The ring $R_\mathfrak{p}$ is often called the \emph{localization of $R$ at $\mathfrak{p}$}.
\end{definition}
\begin{lemma}\label{lemma:localization_is_a_local_field}
	Let $R$ be a ring, and $\mathfrak{p}$ a prime ideal. Then $R_\mathfrak{p}$ is a local ring with $\mathfrak{m}_\mathfrak{p}:=\mathfrak{p} R_\mathfrak{p}$ being the maximal ideal.
\end{lemma}
\begin{proof}
	It is clear that $\mathfrak{m}_\mathfrak{p}$ is an ideal, for $\frac{a}{b}\in\mathfrak{m}_\mathfrak{p}$ if and only if $a\in\mathfrak{p}$ and $b\in R\setminus\mathfrak{p}$. To show that it is maximal, let $\frac{a}{b}\notin \mathfrak{m}_\mathfrak{p}$. Then $a\notin\mathfrak{p}$, so $a\in R\setminus\mathfrak{p}$ which means that $a$ is a unit. So $\mathfrak{m}_\mathfrak{p}$ is maximal. It is clearly unique for if $\mathfrak{a}\neq \mathfrak{m}_\mathfrak{p}$ is another maximal ideal, then $\mathfrak{a}$ contains an element not in $\mathfrak{m}_\mathfrak{p}$, which is a unit by the same argument as before. 
\end{proof}

\begin{lemma}
	Let $R$ be a noetherian ring. Then $R_\mathfrak{p}$ is noetherian.
\end{lemma}
\begin{proof}
	Apply Lemma \ref{lemma:Noetherian_ring_fraction_field}.
\end{proof}
		
\begin{theorem}\cite[page 29]{Reid},\cite[\href{https://stacks.math.columbia.edu/tag/02LV}{Tag 02LV}]{stacks-project}\label{thm:reduced_ring_embedded_in_product_of_fields}
	Let $R$ be a reduced ring. Then
	\begin{enumerate}[label=(\roman*)]
		\item Let $\mathfrak{p}$ be a minimal prime ideal. Then $R_\mathfrak{p}$ is a field.
		\item $R$ is a subring of products of rings, and $$R\rightarrow \prod_{\mathfrak{p} \text{ minimal prime}}R_\mathfrak{p}$$ is an embedding, i.e. $R$ can be embedded in a product of fields.
	\end{enumerate} 
\end{theorem}

\subsection{Colimits}\label{sec:filtered_colim}
We end this section with a small discussion about (filtered) colimits. This section is based on \cite{Riehl,stacks-project}. 

\begin{definition}
	Let $F:I\rightarrow C$ be a functor. A \emph{cone under $F$ with nadir $c$} is a collection of morphisms $\phi_i:F(i)\rightarrow c$, indexed by $i\in I$, such that for all morphisms $f:j\rightarrow k$ with $j,k\in I$, the following diagram commutes:
\[\begin{tikzcd}
		F(j) \arrow[rr, "F(f)"] \arrow[dr, "\phi_j"'] & & \arrow[dl, "\phi_k"] F(k) \\
		& c &
	\end{tikzcd}\]
\end{definition}
\begin{definition}
	A \emph{colimit of $F:I\rightarrow C$} is an object in $C$, denoted by $\mathrm{colim}_{i\in I} F$, together with a cone under $F$ with nadir $\mathrm{colim}_{i\in I} F$ satisfying the following: consider any other cone under $F$ with nadir $d$ with morphisms $\psi_j:F(j)\rightarrow d,$ with $j\in I$, then there exists a unique morphism $u:\mathrm{colim}_I F\rightarrow d$ such that the following diagram commutes:
	\[\begin{tikzcd}
		F(j) \arrow[rr, "\phi_j"] \arrow[dr, "\psi_j"'] & & \arrow[dl, "u"] \mathrm{colim}_IF \\
		& d &
	\end{tikzcd}\]
\end{definition}
\begin{example}\label{ex:colim_direct_sum}
	Let $I$ be some discrete set, and recall that $\mathrm{Vect}_\C$ is the category of $\C$-vector spaces over. Consider $F:I\rightarrow\mathrm{Vect}_\C$ where $F(i)=V_i$. Then the colimit is given by the direct sum, i.e. $$\mathrm{colim}_{i\in I} F= \bigoplus_{i\in I}V_i.$$ The morphisms $\phi_j$ here are exactly the injective maps $V_i\rightarrow \bigoplus_{i\in I}V_i$.
\end{example}
\begin{definition}
	Let $(I,\leq)$ be a poset. If for all $i,j\in I$ there exists $k\in I$ such that $i\leq k$ and $j\leq k$, we call $(I,\leq)$ a \emph{filtered poset}. When considering a colimit of $F:I\rightarrow C$ where $I$ is filtered, we call the colimit a \emph{filtered colimit}.
\end{definition}
\begin{example}\label{ex:filtered_colim}
	\begin{itemize}
		\item Let $A$ be a set, $I$ the set of finite subsets of $A$. Then $(I,\subseteq)$ is a filtered poset. Considering the constant functor $F:I\rightarrow \mathrm{Set}$, given by $i\mapsto i$ and morphisms $f\mapsto \mathbf{1}_i$, we get that $$\mathrm{colim}_{i\in I}\,i = A.$$
		\item Let $R$ be a ring, and $I$ the set of finite subsets of $R$. Let $i\in I$, and consider $\C[i]$. Note that $\C[i]$ is noetherian since $i$ is finite, then one sees by the previous example that
		\begin{align*}
			\mathrm{colim}_{i\in I}\mathbb{C}[i]=R.
		\end{align*} Note that if $R$ is reduced, then $\C[i]$ is also reduced. 
		\item Let $G=K_\C$ be as in Section \ref{sec:Mahtieu's_proof}, and let $R$ be any $G$-algebra. We note then that, since $G$ is reductive, that $$R\simeq \bigoplus_{\pi \in\widehat{G}}R_\pi$$ where $R_\pi$ is the isotypical component of $\pi\in\widehat{G}$. Note that $R_\pi$ can be written as a (possible infinite) sum of copies of the irreducible representation $(\pi,V_\pi)$, i.e. $R_\pi\simeq\bigoplus_{i\in I_\pi}V_\pi$
		where $I_\pi$ is some index set. In other words, $$R\simeq\bigoplus_{\pi\in\widehat{G}}\,\bigoplus_{i\in I_\pi}V_\pi.$$ Take $J$ as the set of finite unions of irreducible subrepresentations of $R$. Then again $(J,\subseteq)$ is a filtered poset, and using the same argument in previous example, we see that $$\mathrm{colim}_{j\in J}\C[j]=R.$$ In particular, note that each $\C[j]$ is now a $G$-invariant noetherian subring of $R$. 
	\end{itemize}
\end{example}

\begin{lemma}\label{lemma:categorical_nonsense_why_swapping_is_okay}
	Let $R$ be a $G$-algebra, and let $I$ the set of finite unions of irreducible subrepresentations of $R$. Then $$\mathrm{colim}_{i\in I}\bigoplus_{\tau\in \widehat{G}}\C[i]_\tau=\bigoplus_{\tau\in \widehat{G}}\mathrm{colim}_{i\in I}\C[i]_\tau.$$
\end{lemma}

To prove Lemma \ref{lemma:categorical_nonsense_why_swapping_is_okay}, we need the following theorem:

\begin{theorem}\cite[Thm. 3.8.1]{Riehl}\label{thm:colims_commute}
	Consider a functor $F:I\times J\rightarrow C$. If $\mathrm{colim}_{i\in I}\mathrm{colim}_{j\in J} F(i,j)$ and $\mathrm{colim}_{j\in J}\mathrm{colim}_{i\in I} F(i,j)$ exist in $C$, they are isomorphic and define the limit $\mathrm{colim}_{I\times J}F(i,j)$.
\end{theorem}

\begin{proof}[Proof of Lemma \ref{lemma:categorical_nonsense_why_swapping_is_okay}]
	Note that we are in the situation of Example \ref{ex:filtered_colim}, and that $\bigoplus_{\tau\in\widehat{G}}$ is just a colimit  by Example \ref{ex:colim_direct_sum}. We would like to apply Theorem \ref{thm:colims_commute} to be done.
	
	To do this in the category of $G$-algebras, define the functor $F:I\times \widehat{G}\rightarrow G\mathrm{-Alg}$ by $F(i,\tau)=\C[i]_\tau$, with morphisms given by $F(i,\tau)\rightarrow F(j,\tau)$ for any $\tau\in\widehat{G}$. This functor is well-defined, because Schur's Lemma states that only non-zero morphisms $f:(i,\tau)\rightarrow(j,\rho)$ can occur whenever $\tau\simeq \rho$. Here $(i,\tau),(j,\rho)\in I\times \widehat{G}$. 
	
	Also note that by Schur's Lemma, there are no non-zero morphisms $\mathrm{colim}_{i\in I}\C[i]_\tau\rightarrow \mathrm{colim}_{i\in I}\C[i]_{\rho}$ whenever $\tau\not\simeq\rho$. This in particular defines $\bigoplus_{\tau\in\widehat{G}}\mathrm{colim}_{i\in I}\C[i]_\tau$. On the other hand, $\mathrm{colim}_{i\in I}\bigoplus_{\tau\in\widehat{G}}\C[i]_\tau$ also exists, for there are no non-zero morphisms $\C[i]_\tau\rightarrow \C[i]_\rho$ whenever $\tau\not\simeq\rho$. Applying Theorem \ref{thm:colims_commute} gives the result.
\end{proof}
\end{appendices}

\newpage
\printbibliography

@book{KnappBeyond,
	author		= {Anthony W. Knapp},
	title		= {Lie Groups Beyond an Introduction, 2nd ed.},
	year 		= {2002},
	publisher	= {Birkh\"auser Boston},
	isbn		= {0-8176-4259-5}
}

@book{KnappRepresentation,
	author		= {Anthony W. Knapp},
	title		= {Representation Theory of Semisimple Groups: An Overview Based on Examples},
	year 		= {1986},
	publisher	= {Princeton University Press},
	isbn		= {0-691-08401-7}
}

@book{Kirillov,
	author    = {Kirillov Jr., Alexander}, 
	title     = {{An introduction to Lie groups and Lie algebras}},
	publisher = {Cambridge University Press},
	year      = 2008,
	isbn      = {978-0-521-88969-8}
}

@inproceedings{Mathieu,
	title={{Some Conjectures About Invariant Theory and their Applications}},
	author={Mathieu, Olivier},
	editor={Alex, J. and Cauchon, G.},
	booktitle={{Alg\`ebre non commutative, groupes quantiques et invariants}},
	address={Reims},
	volume={2},
	pages={263--279},
	year={1997},
	publisher={Soci\'et\'e Math\'ematique de France}
}

@book{Bump,
	author    = {Bump, Daniel}, 
	title     = {{Lie Groups, Second Edition}},
	publisher = {Springer New York},
	year      = 2013,
	isbn	  = {978-1-4614-8023-5}
}

@article{Duistermaat,
	title={{Constant terms in powers of a Laurent polynomial}},
	author={Duistermaat,  Johannes Jisse and Van Der Kallen, Wilberd},
	journal={Indagationes Mathematicae},
	volume={9},
	number={2},
	pages={221--231},
	year={1998},
	publisher={Elsevier}
}

@book{Procesi,
	author    = {Procesi, Claudio}, 
	title     = {{Lie Groups, An Approach through Invariants and Representations}},
	publisher = {Springer New York},
	year      = 2007,
	isbn	  = {978-0-387-28929-8}
}

@book{Humphreys,
	author    = {Humphreys, James E.}, 
	title     = {{Introduction to Lie Algebras and Representation Theory}},
	publisher = {Springer-Verlag New York},
	year      = 1972,
	isbn	  = {0-387-90053-5}
}

@book{Hall,
	author    = {Hall, Brian C.}, 
	title     = {{Lie Groups, Lie Algebras, and Representations}},
	publisher = {Springer Cham},
	year      = 2015,
	isbn	  = {978-3-319-13466-6}
}

@article{Keller,
	title={{Ganze Cremona-Transformationen}},
	author={Keller, Ott Heinrich},
	journal={Monatshefte f\"ur Mathematik und Physik},
	pages={299–306},
	volume = {47},
	year={1939}
}

@article{Kumar,
	title={{Tensor Product Decomposition}}, 
	author={Kumar, Shrawan},
	journal = {Proceedings of the International Congress of Mathematicians 2010 (ICM 2010)},
	pages = {1226-1261},
	year = {2010},
	doi = {10.1142/9789814324359_0094}
}

@article{Jacobian_overzicht_artikel,
	title={{The Jacobian conjecture: reduction of degree and formal expansion of the inverse}},
	author={Bass, Hyman and Connell, Edwin H. and Wright, David},
	journal = {Bulletin of the American Mathematical Society},
	volume = {7},
	number = {2},
	year={1982}
}

@book{Wallach,
	author    = {Wallach, Nolan R.}, 
	title     = {{Real Reductive Groups I}},
	publisher = {Academic Press Inc.},
	year      = 1988,
	isbn	  = {0-12-732960-9}
}

@book{vdEssen,
	author    = {van den Essen, Arno}, 
	title     = {{Polynomial Automorphisms and the Jacobian Conjecture}},
	publisher = {Birkh\"auser Basel},
	year      = 2000,
	isbn	  = {3-7643-6350-9}
}

@book{vdEssen_en_de_rest,
	author    = {van den Essen, Arno and Kuroda, Shigeru and Crachiola, Antony J.}, 
	title     = {{Polynomial Automorphisms and the Jacobian Conjecture, New Results from the Beginning of the 21st Century}},
	publisher = {Birkh\"auser},
	year      = 2021,
	isbn	  = {978-3-030-60535-3}
}

@book{Reid,
	author    = {Reid, Miles}, 
	title     = {{Undergraduate Commutative Algebra}},
	publisher = {Cambridge University Press},
	year      = 1995,
	isbn	  = {0-521-45889-7}
}

@book{Broecker_tom_Dieck,
	author    = {Br\"ocker,Theodor and tom Dieck, Tammo}, 
	title     = {{Representations of Compact Lie Groups}},
	publisher = {Springer Berlin, Heidelberg},
	year      = 2003,
	isbn	  = {978-3-540-13678-1}
}

@article{MohJacobian,
	title={{On the Jacobian conjecture and the configurations of roots}},
	author={Moh, T. T.},
	journal={Journal f{\"u}r die reine und angewandte Mathematik},
	volume={340},
	pages={140-212},
	year={1983}
}

@article{WangJacobian,
	title={{A Jacobian Criterion for Separability}},
	author={Wang, Stuart Sui-Sheng},
	journal = {Journal of Algebra},
	volume = {65},
	number = {2},
	pages = {453-494},
	year = {1980},
	issn = {0021-8693},
	doi = {https://doi.org/10.1016/0021-8693(80)90233-1}
}

@book{Abhyankar,
	title={{Lectures on expansion techniques in algebraic geometry}},
	author={Abhyankar, Shreeram Shankar and Singh, Balwant},
	year={1977},
	publisher={Tata Institute of Fundamental Research, Bombay}
}

@article{Nakai-Baba,
	title={{A generalization of Magnus' theorem}},
	author={Nakai, Yoshikazu and Baba, Kiyoshi},
	journal = {Osaka Journal of Mathematics},
	volume = {14},
	pages = {403-409},
	year = {1977}
}

@article{WrightJacobian,
	title={{The amalgamated free product structure of $GL_2(k[X_1,\ldots,X_n])$ and the weak Jacobian Theorem in two variables}},
	author={Wright, David},
	journal = {Journal of Pure and Applied Algebra},
	volume = {12},
	pages = {235-251},
	year = {1978}
}

@article{BondtVanDenEssen,
	title={{A reduction of the Jacobian conjecture to the symmetric case}},
	author={De Bondt, Michiel and Van den Essen, Arno},
	journal={Proceedings of the American Mathematical Society},
	volume={133},
	number={8},
	pages={2201--2205},
	year={2005}
}

@article{Smale,
	title={{Mathematical Problems for the Next Century}},
	author={Smale, Steve},
	journal={The Mathematical Intelligencer},
	volume={20},
	pages={7-15},
	year={1998},
	doi={https://doi.org/10.1007/BF03025291}
}

@article{Wright_On_the_Jacobian_Conj,
	title={{On the Jacobian Conjecture}},
	author={Wright, David},
	journal={Illinois Journal of Mathematics},
	volume={25},
	number={3},
	pages={423-440},
	year={1981}	
}

@article{Nori,
	title={{The Integral of Powers of a Function}},
	author={Nori, Madhav V.},
	journal={Contemporary Mathematics},
	volume={312},
	pages={163-175},
	year={2002}	
}

@article{Zhao_Mathieu_Subspaces,
	title = {{Generalizations of the image conjecture and the Mathieu conjecture}},
	journal = {Journal of Pure and Applied Algebra},
	volume = {214},
	number = {7},
	pages = {1200-1216},
	year = {2010},
	issn = {0022-4049},
	doi = {https://doi.org/10.1016/j.jpaa.2009.10.007},
	author = {Zhao, Wenhua}
}

@book{Bourbaki_ch4_6,
	author    = {Bourbaki, Nicolas}, 
	title     = {{Groupes et alg\`ebres de Lie, Chapitres IV, V et VI}},
	publisher = {Springer Berlin, Heidelberg},
	year      = 2006,
	isbn	  = {978-3-540-34490-2}
}

@book{Bourbaki_ch7_8,
	author    = {Bourbaki, Nicolas}, 
	title     = {{Groupes et alg\`ebres de Lie, Chapitres VII et VIII}},
	publisher = {Springer Berlin, Heidelberg},
	year      = 2006,
	isbn	  = {978-3-540-33939-7}
}

@article{Stembridge,
	title = {{The Partial Order of Dominant Weights}},
	author = {Stembridge, John R.},
	journal = {Advances in Mathematics},
	volume = {136},
	number = {2},
	pages = {340-364},
	year = {1998}
}

@book{Duistermaat_Kolk,
	author    = {Duistermaat, Johannes Jisse and Kolk, Johan A. C.}, 
	title     = {{Lie Groups}},
	publisher = {Springer Berlin, Heidelberg},
	year      = 1999,
	isbn 	  = {978-3-540-15293-4}
}

@article{Faraut_Koranyi_1990,
	title={{Function spaces and reproducing kernels on bounded symmetric domains}},
	volume={88}, 
	DOI={10.1016/0022-1236(90)90119-6}, 
	number={1}, 
	journal={Journal of Functional Analysis}, 
	author={Faraut, J. and Koranyi, A.}, 
	year={1990},
	pages={64–89} 
}

@article{Arazy_1992,
	title = {{Realization of the invariant inner products on the highest quotients of the composition series}},
	author = {Arazy, Jonathan},
	journal = {Arkiv f\"or Mathematik},
	volume = {30},
	pages = {1-24},
	year = {1992}
}

@book{Atiyah,
	author    = {Atiyah, Michael and Macdonald, I. G.}, 
	title     = {{Introduction to Commutative Algebra}},
	publisher = {Westview Press},
	year      = 1969,
	isbn 	  = {0-201-40751-5}
}

@misc{stacks-project,
	shorthand    = {Stacks},
	author       = {The {Stacks Project Authors}},
	title        = {\textit{Stacks Project}},
	howpublished = {\url{https://stacks.math.columbia.edu}},
	year         = {2018},
}

@book{Matsumura,
author    = {Matsumura, Hideyuki}, 
title     = {{Commutative Algebra}},
publisher = {W.A. Benjamin, Inc.},
year      = 1970,
isbn 	  = {8053-7025-0}
}

@article{Jacobi,
	title = {{De resolutione aequationum per series infinitas}},
	author = {Jacobi, C.},
	journal = {Journal f\"ur die reine und angewandte Mathematik},
	volume = {5},
	pages = {257-286},
	year = {1830}
}

@article{Good,
	title = {{Generalizations to several variables of Lagrange's expansion, with applications to stochastic processes}},
	author = {Good, Irving J.},
	journal = {Mathematical Proceedings of the Cambridge Philosophical Society},
	volume = {56},
	issue = {4},
	pages = {367-380},
	year = {1960}
}

@article{Gessel_multi,
	title = {{A Combinatorial Proof of the Multivariable Lagrange Inversion Formula}},
	author = {Gessel, Ira M.},
	journal = {Journal of Combinatorial Theory, Series A},
	volume = {45},
	pages = {178-195},
	year = {1987}
}

@article{Gessel,
	title = {{Lagrange inversion}},
	author = {Gessel, Ira M.},
	journal = {Journal of Combinatorial Theory, Series A},
	volume = {144},
	pages = {212-249},
	year = {2016}
}

@book{Faraut,
	author    = {Faraut, Jacques}, 
	title     = {{Analysis on Lie Groups: An Introduction}},
	publisher = {Cambridge University Press},
	year      = 2008,
	isbn 	  = {978-0-511-42296-6}
}

@book{Goodman,
	author    = {Goodman, Roe and Wallach, Nolan R.}, 
	title     = {{Symmetry, Representations and Invariants}},
	publisher = {Springer New York},
	year      = 2009,
	isbn 	  = {978-0-387-79852-3}
}

@book{Riehl,
	author    = {Riehl, Emily}, 
	title     = {{Category Theory in Context}},
	publisher = {Dover Pulications Inc., New York},
	year      = 2016,
	isbn 	  = {978-0-486-80903-8}
}
	
\end{document}